\documentclass[11pt,twoside]{article}
\newcommand{\ifims}[2]{#1} 
\newcommand{\ifAMS}[2]{#1}   
\newcommand{\ifau}[4]{#2}  
\newcommand{\ifbook}[2]{#1}   

\usepackage{ifthen}
\usepackage[ampersand]{easylist}
\usepackage{amsmath,amssymb,amsthm}
\usepackage{mathtools}
\usepackage{natbib}
\usepackage{epsfig,graphicx}
\usepackage{comment}
\usepackage{color}
\usepackage{srcltx}
\usepackage[mathscr]{eucal}
\usepackage[math]{easyeqn}
\usepackage{etoolbox}
\usepackage{hyperref}
\hypersetup{
            colorlinks,
            linkcolor=hookersgreen,
            linktoc=blue,
            citecolor=ultramarine,
            urlcolor=black,
            filecolor=black
            }
\usepackage{nameref}
\usepackage{multirow}

\ifims{
\textheight=22cm
\textwidth=14.8cm
\topmargin=0pt
\oddsidemargin=1.0cm
\evensidemargin=1.0cm
\linespread{1.3}

}{ 
}

\numberwithin{equation}{section}
\numberwithin{figure}{section}
\newcounter{example}[section]
\numberwithin{example}{section}
\newcounter{remark}[section]
\numberwithin{remark}{section}
\newtheorem{theorem}{Theorem}[section]
\newtheorem{proposition}[theorem]{Proposition}
\newtheorem{lemma}[theorem]{Lemma}
\newtheorem{corollary}[theorem]{Corollary}

\newtheorem{exmp}[example]{Example}
\newtheorem{rmrk}[remark]{Remark}
\newenvironment{example}{\begin{exmp}\rm}{\end{exmp}}
\newenvironment{remark}{\begin{rmrk}\rm}{\end{rmrk}}

\ifbook{

  }
  {

  }

\renewcommand{\(}{$\,}
\renewcommand{\)}{\,$}

\def\nquad{\hspace{-1cm}}
\def\eqdef{\stackrel{\operatorname{def}}{=}}

\def\tow{\stackrel{w}{\longrightarrow}}

\newcommand{\cc}[1]{\mathscr{#1}}
\newcommand{\bb}[1]{\boldsymbol{#1}}

\renewcommand{\bar}[1]{\overline{#1}}

\renewcommand{\tilde}[1]{\widetilde{#1}}

\newcommand{\thankstitle}[1]{\ifthenelse{\equal{#1}{}}{}{\thanks{#1}}}
\newcommand{\thanksau}[1]{\ifthenelse{\equal{#1}{}}{}{\thanks{#1}}}

\newcommand{\aua}[6]
{\def\authora{#1}
\def\runauthora{#2}
\def\addressa{#3}
\def\emaila{#4}
\def\affiliationa{#5}
\def\thanksa{#6}}

\newcommand{\aub}[6]
{\def\authorb{#1}
\def\runauthorb{#2}
\def\addressb{#3}
\def\emailb{#4}
\def\affiliationb{#5}
\def\thanksb{#6}}

\def\theauthors{
\ifau{ 
  \author{
    \authora
    \thanksau{\thanksa}
    \\[5.pt]
    \addressa \\
    \texttt{ \emaila}
  }
}
{  
  \author{
    \authora
    \thanksau{\thanksa}
    \\[5.pt]
    \addressa \\
    \texttt{ \emaila}
    \and
    \authorb
    \thanksau{\thanksb}
    \\[5.pt]
    \addressb \\
    \texttt{ \emailb}
  }
}
{   
  \author{
    \authora
    \thanksau{\thanksa}
    \\[5.pt]
    \addressa \\
    \texttt{ \emaila}
    \and
    \authorb
    \thanksau{\thanksb}
    \\[5.pt]
    \addressb \\
    \texttt{ \emailb}
    \and
    \authorc
    \thanksau{\thanksc}
    \\[5.pt]
    \addressc \\
    \texttt{ \emailc}
  }
} {   
  \author{
    \authora
    \thanksau{\thanksa}
    \\[5.pt]
    \addressa \\
    \texttt{ \emaila}
    \and
    \authorb
    \thanksau{\thanksb}
    \\[5.pt]
    \addressb \\
    \texttt{ \emailb}
    \and
    \authorc
    \thanksau{\thanksc}
    \\[5.pt]
    \addressc \\
    \texttt{ \emailc}
    \and
    \authord
    \thanksau{\thanksd}
    \\[5.pt]
    \addressd \\
    \texttt{ \emaild}
  }
}
}
\def\therunauthors{ \hfill \textsc{\small
 \ifau{\runauthora}
      {\runauthora and \runauthorb}
      {\runauthora, \runauthorb, and \runauthorc}
      {\runauthora, \runauthorb, \runauthorc, and \runauthord}
 }
 \hfill
 }

\renewcommand{\Gamma}{\varGamma}
\renewcommand{\Pi}{\varPi}
\renewcommand{\Sigma}{\varSigma}
\renewcommand{\Delta}{\varDelta}
\renewcommand{\Lambda}{\varLambda}
\renewcommand{\Psi}{\varPsi}
\renewcommand{\Phi}{\varPhi}
\renewcommand{\Theta}{\varTheta}
\renewcommand{\Omega}{\varOmega}
\renewcommand{\Xi}{\varXi}
\renewcommand{\Upsilon}{\varUpsilon}

\def\arginf{\operatornamewithlimits{arginf}}

\def\argmax{\operatornamewithlimits{argmax}}

\def\av{\bb{a}}

\def\uv{\bb{u}}

\def\wv{\bb{w}}
\def\xv{\bb{x}}

\def\zv{\bb{z}}

\def\Uv{\bb{U}}

\def\Xv{\bb{X}}
\def\Yv{\bb{Y}}

\def\etav{\bb{\eta}}
\def\gammav{\bb{\gamma}}

\def\lambdav{\bb{\lambda}}

\def\thetav{\bb{\theta}}

\def\xiv{\bb{\xi}}
\def\zetav{\bb{\zeta}}
\def\Psiv{\bb{\Psi}}

\def\sumi{\sum_{i=1}^{n}}

\usepackage{color}

\definecolor{blue(pigment)}{rgb}{0.2, 0.2, 0.6}
\definecolor{ultramarine}{rgb}{0.07, 0.04, 0.56}
\definecolor{darkspringgreen}{rgb}{0.09, 0.45, 0.27}
\definecolor{hookersgreen}{rgb}{0.0, 0.44, 0.0}
\definecolor{plum(traditional)}{rgb}{0.56, 0.27, 0.52}
\definecolor{purple(html/css)}{rgb}{0.5, 0.0, 0.5}
\definecolor{magenta(dye)}{rgb}{0.79, 0.08, 0.48}

\def\II{\mathcal{I}}

\def\rrbias{\rr_{b}}

\def\R{I\!\!R}
\def\E{I\!\!E}

\def\P{I\!\!P}

\def\kappa{\varkappa}

\def\Frobg{\Lambda}

\def\diag{\operatorname{diag}}

\def\Fr{\operatorname{Fr}}

\def\ND{\mathcal{N}}

\def\oper{\operatorname{op}}

\def\Var{\operatorname{Var}}

\def\T{\top}
\def\tr{\operatorname{tr}}

\def\Sphere{\cc{S}}
\def\CONST{\mathtt{C} \hspace{0.1em}}
\def\cond{\, \big| \,}

\def\nsize{{n}}

\def\sumi{\sum_{i=1}^{\nsize}}

\def\ex{\mathrm{e}}

\def\Id{I\!\!\!I}
\def\Ind{\operatorname{1}\hspace{-4.3pt}\operatorname{I}}

\def\alp{\alpha}

\def\AAla{\TAU_{\zev}}

\def\AP{A}

\def\BB{I\!\!B}     
\def\BB{B}

\def\BBB{\cc{B}}

\def\cdens{\phi}

\def\CA{\cc{A}}

\def\CAt{\tilde{\CA}}
\def\CAGP{\CA_{\GP}}

\def\CONSTPsi{\CONST_{\Psi}}

\def\CONSTfour{\kappa}
\def\CONSTru{\CONST_{0}}

\def\CS{\cc{E}}

\def\DP{D}
\def\DPc{\DP}

\def\DPGP{\DP_{\GP}}

\def\DPt{\tilde{\DPc}}
\def\DPGPt{\DPt_{\GP}}

\def\DPt{\tilde{\DP}}

\def\dimp{p}

\def\dimA{\mathtt{p}}
\def\dimAA{\dimA_{\tau}}

\def\dimG{\dimA_{\GP}}

\def\dimAla{\dimA_{\lambdav}}

\def\dimt{\tilde{\dimA}}

\def\dens{f}

\def\err{\diamondsuit}
\def\eps{\epsilon}			
\def\eps{\varepsilon}

\def\fs{f}

\def\gaussv{\bb{\gauss}}
\def\gauss{\gamma}

\def\gm{\mathtt{g}}
\def\gmc{\gm_{c}}
\def\gmb{\gm}

\def\gp{g}

\def\GP{G}

\def\HM{\mathbb{H}}
\def\HP{\mathsf{H}}

\def\hyper{\kappa}

\def\IF{\Bbb{F}}

\def\LT{L}
\def\LGP{\LT_{\GP}}

\def\LL{\cc{L}}

\def\lambdav{\bb{\lambda}}

\def\ldens{\ell}

\def\mm{m}

\def\mms{\mm^{*}}

\def\muc{\mu_{c}}

\def\MM{\cc{M}}

\def\nunu{\nu_{0}}

\def\PAAr{\Phi}

\def\PAArg{\PAAr_{\lambdav}}

\def\Pone{P}
\def\Pdom{\mu_{0}}
\def\PDOM{\bb{\mu}_{0}}

\def\PEF{\cc{P}}

\def\PG{\P'}

\def\Psimean{\bar{\Psi}}
\def\Proj{\Pi}

\def\Prior{\Pi}

\def\QP{Q}
\def\QPGP{\QP | \GP}

\def\qq{q}

\def\rhot{t}

\def\rr{\mathtt{r}}

\def\rups{\rr_{0}}

\def\Spsi{S}

\def\supA{\lambda}

\def\thetav{\bb{\theta}}
\def\thetavs{\thetav^{*}}
\def\thetavc{\thetav'}
\def\thetavd{\thetav^{\circ}}

\def\TAU{\mathcal{T}}

\def\Thetad{\Theta^{\circ}}

\def\uvd{\uv^{\circ}}

\def\UV{\mathcal{U}}
\def\UVd{\UV^{\circ}}

\def\vA{\mathtt{v}}

\def\VP{V}
\def\VPc{\VP_{0}} 	
\def\VPc{\VP}

\def\VPGP{\VP | \GP}

\def\VV{\mathbb{V}}

\def\vtheta{\vartheta}
\def\vthetav{\bb{\vtheta}}
\def\vthetavb{\bar{\vthetav}}

\def\wv{\bb{w}}
\def\wvd{\wv^{\circ}}

\def\xb{\bar{\xv}}

\def\xvd{\xv^{\circ}}
\def\xx{\mathtt{x}}
\def\xxc{\xx_{c}}

\def\xvd{\xv^{\circ}}

\def\Xs{X^{*}}

\def\XX{\cc{X}}

\def\zev{\zv}

\def\zq{z}
\def\zqc{\zq_{c}}

\def\zz{\mathfrak{z}}

\def\thetitle{Accuracy of Gaussian approximation in nonparametric Bernstein -- von Mises Theorem}
\def\theruntitle{Accuracy of Gaussian approximation in BvM Theorem}

\def\theabstract{
The prominent Bernstein -- von Mises (BvM) result claims that the posterior distribution after centering by the efficient estimator and standardizing by the square root of the total Fisher information is nearly standard normal.
In particular, the prior completely washes out from the asymptotic posterior distribution.
This fact is fundamental and justifies the Bayes approach from the frequentist viewpoint.
In the nonparametric setup the situation changes dramatically and the impact of prior becomes essential even for the contraction of the posterior; see~\cite{vdV2008}, \cite{Bo2011}, \cite{CaNi2013,CaNi2014} for different models like Gaussian regression or i.i.d. model in different weak topologies.
This paper offers another non-asymptotic approach to studying the behavior of the posterior 
for a special but rather popular and useful class of statistical models and for Gaussian priors. 
First we derive tight finite sample bounds on posterior contraction in terms of the so-called effective dimension of the parameter space.
Our main results describe the accuracy of Gaussian approximation of the posterior.
In particular, we show that restricting to the class of all centrally symmetric credible sets around the penalized
maximum likelihood estimator (pMLE) allows to get Gaussian approximation up to order \( n^{-1} \).
We also show that the posterior distribution mimics well the distribution of the pMLE and reduce the question of reliability of credible sets 
to consistency of the pMLE-based confidence sets. 
The obtained results are specified for nonparametric log-density estimation and generalized regression.
}

\def\kwdp{62F15}
\def\kwds{62F25}

\def\thekeywords{posterior, concentration, contraction Gaussian approximation}

\def\thethankstitle{The research was supported by the Russian Science Foundation grant 19-71-30020.}

\aub{Maxim Panov}
{Panov, M. }
{Skolkovo Institute \\ of Science and Technology \\
    121205 Moscow, Russia, 
    }
{m.panov@skoltech.ru}
{Skolkovo Institute of Science and Technology, }
{}

\aua
{Vladimir Spokoiny}
{Spokoiny, V. }
{
    Weierstrass Institute and HU Berlin,  \\
    HSE University, Skoltech and IITP RAS,
    \\
    Mohrenstr. 39, 10117 Berlin, Germany
    }
{spokoiny@wias-berlin.de}
{Weierstrass-Institute and Humboldt University Berlin}
{
   Financial support by the German Research Foundation (DFG) through the Collaborative Research Center 1294 ``Data assimilation'' is gratefully acknowledged.
} 

\def\GPm{\GP_{1}}
\def\rsmall{\varrho}
\def\rrbias{\rr_{b}}

\def\hhh{v}
\def\Xs{\cc{X}}

\def\CONSTi{\mathtt{C}}
\def\CONSTfour{\CONSTi_{\fs}}
\def\CONSTphi{\CONSTi_{\cdens}}
\def\CONSTPsi{\CONSTi_{\Psi}}

\def\CONSTpsi{\CONSTi_{\psi}}
\def\CONSTIF{\CONSTi_{\IF}}

\def\CGP{w}
\def\Psid{\Psi_{\thetav}}

\def\detHP{q_{0}}
\def\CAt{{\CS}}

\def\smp{s}
\def\mux{\mu_{\xx}}

\def\dltwu{\tau}
\def\dltw{\delta}
\def\HPc{H}

\begin{document}
\thispagestyle{empty}

\title{\thetitle\thankstitle{\thethankstitle}}

\theauthors

\maketitle
\begin{abstract}
  \theabstract
\end{abstract}

\ifAMS
    {\par\noindent\emph{AMS 2010 Subject Classification:} Primary \kwdp. Secondary \kwds}
    {\par\noindent\emph{JEL codes}: \kwdp}


\tableofcontents


\section{Introduction}
\label{Swilksint}
Bernstein -- von Mises (BvM) Theorem is one of most prominent results in statistical inference.
It claims that the posterior measure is asymptotically normal with the mean close to the maximum likelihood estimator (MLE) and the variance close to the variance of the MLE.
This explains why this result is often considered as the Bayesian counterpart of the frequentist Fisher Theorem about asymptotic normality of the MLE.
The BvM result provides a theoretical background for different Bayesian procedures.
In particularly, one can use Bayesian computations for evaluation of the MLE and its variance.
Also one can build elliptic credible sets using the first two moments of the posterior.
The main questions to address by studying the behavior of a nonparametric Bayes procedure are 
\smallskip
\begin{easylist}[itemize]
	& concentration: find possibly small concentration sets of the posterior distribution;
	& asymptotic normality or any other asymptotic approximation of the posterior;
	& covering: whether one can use credible sets as frequentist confidence sets.
\end{easylist}
\smallskip
The classical versions of the BvM Theorem claim that the posterior concentrates on a root-n vicinity of the true 
parameter, after proper centering and scaling it is root-n standard normal, and credible sets can be well used 
in place of classical confidence sets.
However, these results require a fixed finite dimensional parameter set, 
correct model specification, and large samples.
We refer to~\cite{VW1996,vdV00} for a detailed historical overview.

Any extension of the BvM approach to the case of a large or infinite dimensional parameter space appears to be very involved, in particular, more involved than the expansions of the maximum likelihood estimate.
The first problem is related to the posterior concentration. 
Such a result requires to bound the integral of the likelihood process in the complement of the local vicinity and this is a hard task in the nonparametric setup.
The second problem is due to fact that a standard Gaussian measure on \( \R^{\infty} \) is only defined in a weak sense.
In particular, it does not concentrate on any \( \ell_{2} \) ball in \( \R^{\infty} \).
This makes it difficult to study asymptotically the total variation distance between the scaled posterior and the Gaussian law. 
We refer to \cite{CaNi2013,CaNi2014}, and \cite{GhVa2017} for a more discussion. 
Our approach can be called \emph{preasymptotic}: we fix the sample and study 
the distance between the posterior and an accompanying Gaussian distribution.
Similar approach was recently used in \cite{yano2020} for a special \( \ell_{\infty} \) topology
and it applies to the case of a moderate or high parameter dimension \( \dimp \) with the rate 
\( \bigl( \log \dimp \bigr)^{a} \nsize^{-1/6} \).
In our approach the parameter dimension can be arbitrary, however, the so called effective dimension 
has to be relatively small.
Restricting to centrally symmetric credible sets allows to get the accuracy of order \( \nsize^{-1} \) 
for Gaussian approximation.
One more crucial issue is an inconsistency problem: in some situations, 
 Bayesian credible sets do not contain the true parameter with the probability close to one; cf.~\cite{Cox1993,Fr1999}, or~\cite{KlVa2006,KlVa2012}.
%
It appears that the posterior in the case of a Gaussian prior is nearly normal but its first two moments mimic the \emph{penalized maximum likelihood estimator} (pMLE) with the quadratic penalization coming from the prior distribution.
It is well known that the penalization yields some estimation bias.
If the squared bias exceeds the variance of the penalized MLE, 
the Bayesian credible sets become unreliable.

The main results of this paper describe the properties of the posterior distribution for Gaussian priors 
in a high-dimensional or nonparametric setups.
In particular, we establish a nonasymptotic upper bound on concentration and on the error of Gaussian approximation for the posterior in
total variation distance in terms of efficient dimension of the problem. 
We also show that the latter bound can be dramatically improved when restricting to the class of centrally symmetric sets around pMLE.
Our assumptions include two important conditions. 
The first one requires that the stochastic part of the log-likelihood is linear in the target parameter, 
while the second one is about concavity of the expected log-likelihood. 
These two conditions are automatically fulfilled in a number of popular models like Gaussian, Poissonian, Binary or Generalized Linear regression, log-density estimation, linear diffusion, etc.
Under these assumptions we manage to state and prove our results in a concise way and avoid 
the machinery of the empirical process theory. 
A forthcoming paper \cite{Sp2019NLI} explains how the approach and the results can be extended to much more general setups
including nonlinear (generalized) regression and nonlinear inverse problems with noisy observations.
%
The main contributions of the paper are \emph{finite samples} results with \emph{accuracy guarantees}
 including
\smallskip
\begin{easylist}[itemize]
	& sharp bounds on concentration of pMLE and of the posterior distribution;
	& Gaussian approximation of the posterior with an explicit error term for the total variation distance and for the class of centrally symmetric sets around pMLE;
	& systematic use of an \emph{effective dimension} in place of the total parameter dimension;
	& addressing frequentist validity of Bayesian credible sets; 
	& specification of the results to log-density estimation and generalized regression.
\end{easylist}
\smallskip
The whole approach is \emph{coordinate free}, we do not use any spectral decomposition and/or any basis representation for the target parameter and penalization.
In this paper we suppose the prior to be given and do not address the question of prior selection.
A number of studies explain how an empirical or hierarchical Bayes approach can be used for building 
adaptive confidence sets;
see e.g.
\cite{KnSzVa2016,NiSz2016,SnVa2015}.
We, however, indicate below how our approach can be used to reduce the original problem of Bayes model selection to the well studied Gaussian case using uniform Gaussian approximation; see Section~\ref{Spriorfamily}.

The paper is structured as follows. 
Section~\ref{SnonBvM} describes our setup, presents the main conditions and states the main results 
about the properties of the pMLE and of the posterior.
Section~\ref{Scondgllo} collects our conditions and main notations.
The central notion of effective dimension is discussed in Section~\ref{Seffdim}.
Properties of the pMLE \( \tilde{\thetav}_{\GP} \) are described in Section~\ref{StiteGpro}.
Section~\ref{Spostcontraction} presents our main results about 
Gaussian approximation of the posterior; see  Theorem~\ref{TnonparBvMm} and its Corollary~\ref{CnonparBvMm}.
We also address the issues of contraction and coverage of Bayesian credible sets.
Section~\ref{SGBvM} comments how the result can be applied to the case of the Bayesian nonparametric log-density estimation,
while Section~\ref{SEFcposter} discusses generalized regression estimation. 
Some more results are given in the Supplement.
Section~\ref{Ssomeextensions} presents some extensions including 
the use of posterior mean in place of the pMLE in the construction of credible sets,
or the use of a general prior with a log-concave density in place of a Gaussian one.
It also 
addresses the important question of prior impact.
%
The proofs and auxiliary results are collected in the Appendix.

\section{Nonparametric BvM Theorem}
\label{SnonBvM}
This section discusses an extension of the BvM result for a a class of models with a high-dimensional or infinite dimensional parameter set and for a Gaussian prior.
Compared to existing literature, our results provide finite sample bounds on posterior concentration and on accuracy 
of Gaussian approximation for the posterior. Moreover, we show that the quality of Gaussian approximation 
can be gradually improved up to order \( \nsize^{-1} \) if we only consider credible sets which are centrally symmetric around 
the pMLE.

Below \( \R^{\dimp} \) means a \( \dimp \)-dimensional Euclidean space equipped with the norm
\( \| \cdot \| \), \( \dimp \leq \infty \). 
Scalar product in \( \R^{\dimp} \) is denoted by \( \langle \cdot,\cdot \rangle \).
For a linear operator \( \BB \) in \( \R^{\dimp} \), the norm \( \| \BB \| \) means 
the largest eigenvalue of \( \BB \).
   
First we specify our setup. 
Let \( \Yv \) denote the observed data and \( \P \) mean their distribution.
A general parametric assumption (PA) means that \( \P \) belongs to infinite-dimensional family 
\( (\P_{\thetav}, \thetav \in \Theta \subseteq \R^{\dimp}) \) dominated by a measure \( \PDOM \).
This family yields the log-likelihood function 
\( L(\thetav) = L(\Yv,\thetav) \eqdef \log \frac{d\P_{\thetav}}{d\PDOM}(\Yv) \).
The PA can be misspecified, so, in general, \( L(\thetav) \) is a \emph{quasi log-likelihood}.
The classical maximum likelihood principle suggests to estimate  \( \thetav \) by
maximizing the function \( L(\thetav) \):
\begin{EQA}[c]
	\tilde{\thetav}
	\eqdef
	\argmax_{\thetav \in \Theta} L(\thetav) .
\label{tthetamkGP}
\end{EQA}
If 
\( \P \not\in \bigl( \P_{\thetav} \bigr) \), 
then the estimate \( \tilde{\thetav} \) from~\eqref{tthetamkGP} is still meaningful and it appears to be an estimate of the value \( \thetavs \) defined by maximizing the expected value of \( L(\thetav) \)
w.r.t. \( \P \):
\begin{EQA}[c]
	\thetavs
	\eqdef
	\argmax_{\thetav \in \Theta} \E L(\thetav) .
\label{thetavsd}
\end{EQA}
Such a value \( \thetavs \) is the true parameter under correct model specification and 
it can be viewed as the parameter of the best parametric fit in the general case.
In the Bayes setup, the parameter \( \vthetav \) is a random element following a prior measure \( \Prior \) on the parameter set \( \Theta \).
The posterior describes the conditional distribution of \( \vthetav \) given \( \Yv \) obtained by normalization of the product \( \exp \bigl\{ L(\thetav) \bigr\} \Prior(d\thetav) \).
This relation is usually written as
\begin{EQA}
	\vthetav \cond \Yv
	& \propto &
	\exp \bigl\{ L(\thetav) \bigr\} \, \Prior(d\thetav).
\label{Bayesform}
\end{EQA}
Below we focus on the case of a Gaussian prior.
Without loss of generality, a Gaussian prior \( \Prior(\thetav) \) will be assumed to be centered at zero.
By \( \GP^{-2} \) we denote its covariance matrix, so that, \( \Prior \sim \ND(0, \GP^{-2}) \).
The main question studied below is to understand under which conditions on the prior covariance 
\( \GP^{-2} \) and the model, the BvM-type result holds and what is the error term in the BvM approximation.
For a Gaussian likelihood, the posterior is Gaussian as well and its properties can be studied directly; see e.g.~\cite{Bo2011,Le2011}.   
For the case when the log-likelihood function is not quadratic in \( \thetav \), the study is more involved.
The posterior is obtained by normalizing the product density \( \exp\bigl\{\LGP(\thetav)\bigr\} \) with 
\begin{EQA}
    \LGP(\thetav)
    &=&
    L(\thetav) - \bigl\| \GP \thetav \bigr\|^{2}/2 \, ,
  \label{LGPBvMG}
\end{EQA}
where \( \| \cdot \| \) is the Euclidean norm in \( \R^{\dimp} \).
This expression arises in penalized maximum likelihood estimation, one can treat the prior term \( \bigl\| \GP \thetav \bigr\|^{2}/2 \) as roughness penalty.
Define
\begin{EQA}
    \tilde{\thetav}_{\GP} 
    = 
    \argmax_{\thetav \in \Theta} \LGP(\thetav) ,
    &\quad &
    \thetavs_{\GP} 
    = 
    \argmax_{\thetav \in \Theta} \E \LGP(\thetav).
  \label{ttGatLGttsGELG}
\end{EQA}

\subsection{Conditions} 
\label{Scondgllo} 
This section collects the conditions which are systematically used in the text. 
We mainly require that the stochastic part of the log-likelihood process \( L(\thetav) \) is linear in \( \thetav \),
while its expectation is a smooth concave function of \( \thetav \).
We also implicitly assume that the parameter set \( \Theta \) 
is an open subset of \( \R^{\dimp} \) where \( \dimp \) is typically equal to infinity.
The model and complexity reduction will be done via the the prior structure in terms of the so called effective dimension.

\begin{description}
    \item[\label{LLref} \( \bb{(\LL)} \)]
      \textit{The set \( \Theta \) is open and convex in \( \R^{\dimp} \).
        The function \( \E L(\thetav) \) is concave in \( \thetav \in \Theta \).
      }
\end{description}

\begin{description}
    \item[\label{Eref} \( \bb{(E)} \)]
      \textit{The stochastic component \( \zeta(\thetav) = L(\thetav) - \E L(\thetav) \) of the process \( L(\thetav) \) is linear in \( \thetav \). We denote by \( \nabla \zeta \equiv \nabla \zeta(\thetav) \) its gradient
      and by \( \VP^{2} = \Var(\nabla \zeta) \) its covariance}.
  \end{description}

\begin{description}
	\item[\( \bb{(E\VP)} \)\label{ED0ref}]
  \emph{There exist a positive self-adjoint operator \( \HPc \) with \( \HPc^{2} \geq \VP^{2} \),
	and constant \( \nunu \ge 1 \) such that \( \nabla \zeta \) fulfills
  }
\begin{EQA}[c]
\label{expzetaclocGP}
    \sup_{\uv \in \R^{\dimp}} 
    \log \E \exp\biggl\{
      \lambda \frac{\langle \uv, \HPc^{-1} \nabla \zeta \rangle}
                   {\| \uv \|}
    \biggr\} \le
    \frac{\nunu^{2} \lambda^{2}}{2} \, .
\end{EQA}
\end{description}
Condition~\nameref{ED0ref} basically requires that the normalized score 
\( \xiv = \HPc^{-1} \nabla \zeta \)
is a sub-Gaussian random vector.
One can relax this condition to finite exponential moments for \( |\lambda| \leq \gmb \)
with \( \gmb \) sufficiently large; see \eqref{expgamgm} of Section~\ref{SdevboundnonGauss}. 
In fact, \nameref{ED0ref} is only used to establish the deviation bounds for quadratic forms of
\( \nabla \zeta \); see e.g. \eqref{uTzDGvTxm12z}.
One can directly operate with the quantiles of the corresponding distribution.
In the finite dimensional case \( \dimp < \infty \), one can often take \( \HPc = \VP \); see Section~\ref{SexamplesBvM} for more examples.

Apart the basic conditions \nameref{LLref}, \nameref{Eref}, \nameref{ED0ref}  
we need some local properties of the expected log-likelihood \( \E L(\thetav) \).
Let \( \Thetad \) be a local set. 
It is required that this set contains the concentration set 
\( \CAGP(\rr_{\GP}) \) of the estimate \( \tilde{\thetav}_{\GP} \); see Proposition~\ref{TconstGrDG} below.
%
Define 
\begin{EQA}
    \IF(\thetav)
    & \eqdef &
    - \nabla^{2} \E L(\thetav) , 
    \\
    \IF_{\GP}(\thetav)
    & \eqdef &
    - \nabla^{2} \E \LGP(\thetav) 
    = 
    - \nabla^{2} \E L(\thetav) + \GP^{2} 
    =
    \IF(\thetav) + \GP^{2} ,
  \label{DP2na2Ltn2}
\end{EQA}

\begin{description}
\item[\( \bb{(\HPc\GP)} \) \label{HGref}] 
For all \( \thetav \in \Thetad \), it holds \( \HPc^{2} \leq \IF_{\GP}(\thetav) \) and
\( \HPc \, \IF_{\GP}^{-1}(\thetav) \, \HPc \) is a trace operator:
\begin{EQA}
	\dimA_{\GP}(\thetav)
	& \eqdef &
	\tr \bigl\{ \HPc \, \IF_{\GP}^{-1}(\thetav) \, \HPc \bigr\}
	< \infty .
\label{pHpHm2Gm2ity}
\end{EQA}
\end{description}

Also we require that the function \( \E L(\thetav) \) is four times differentiable on \( \Thetad \).
Define for each \( \thetav \in \Thetad \), and any \( \uv \in \R^{\dimp} \), 
the directional 
Gateaux derivative
\begin{EQA}
    \dltw_{k}(\thetav,\uv)
    & \eqdef &
    \frac{1}{k!}
      \frac{d^{k}}{d\rhot^{k}} \E L(\thetav + \rhot \uv) \biggr|_{\rhot=0} \, ,
    \quad 
    k = 3, 4.
\label{d3utd4utd4dt4c}
\end{EQA}
Clearly \( \dltw_{k}(\thetav,\uv) \) is proportional to \( \| \uv \|^{k} \).
%
Later we need a uniform bound on \( \dltw_{k}(\thetav,\uv) \).
\begin{description}
    \item[\( \bb{(\LL_{0})} \)\label{LL0ref}]
       \textit{
		It holds with \( \HPc \) from \nameref{ED0ref} for \( k=3,4 \) 
		and some \( \rr \) sufficiently large
		}
\begin{EQA}
	\dltwu_{k,\HPc}(\rr)
	\eqdef
	\sup_{\thetav \in \Thetad} \sup_{\| \HPc \uv \| \leq \rr} \rr^{-k} \, \dltw_{k}(\thetav, \uv)
	& < &
	\infty .
\label{dmtuCmnumddL0}
\end{EQA}		
		
\end{description}

In what follows we consider the situation with \( \nsize \) independent observations.
The resulted expected log-likelihood \( \E L(\thetav) \) is of order \( \nsize \) as well as its derivatives. 
Also we will assume that \( \| \HPc^{-2} \| \leq \CONST / \nsize \).
Then for some fixed \( \CONSTi_{3,\dltw}, \CONSTi_{4,\dltw} \)
\begin{EQA}
	\dltwu_{k,\HPc}(\rr)
	& \leq &
	\CONSTi_{k,\dltw} \, \rr^{-k} \nsize \, (\rr \nsize^{-1/2})^{k} 
	=
	\CONSTi_{k,\dltw} \, \nsize^{1-k/2}  
	,
	\qquad
	k=3,4.
\label{dmtuCmnumdd}
\end{EQA}


\subsection{Effective dimension}
\label{Seffdim}
This section discusses the central notion of \emph{effective dimension}.
Define with \( \HPc^{2} \) from \nameref{ED0ref} and \( \IF_{\GP} = \IF_{\GP}(\thetavs_{\GP}) \)
\begin{EQA}[rcccccl]
	\BB_{\GP} 
	& \eqdef &
	\HPc \, \IF_{\GP}^{-1} \, \HPc \,  ,
	\quad
    \dimA_{\GP} 
    & \eqdef &
    \tr \BB_{\GP} \, , 
    \qquad
    \lambda_{\GP}
    & \eqdef & 
    \| \BB_{\GP} \|   \, . 
\label{lamGPDPVPfis}
\end{EQA}
Here \( \| \BB \| \) means the operator norm or the maximal eigenvalue of \( \BB \).
The values \( \tr \BB_{\GP} \) and \( \lambda_{\GP} \) are important because they enter in the definition of the upper quantile function \( \zq(\BB_{\GP}, \xx) \) 
for \( \| \IF_{\GP}^{-1/2} \nabla \zeta \| \); see \eqref{uTzDGvTxm12z} below.
In \cite{SP2013_rough} the quantity \( \dimA_{\GP} \) was called 
the effective dimension.
%
Below we consider the \emph{local effective dimension} \( \dimG(\thetav) \) at a point 
\( \thetav \in \Thetad \) given by
\begin{EQA}
    \dimA_{\GP}(\thetav)
    & \eqdef &
    \tr \bigl\{ \HPc^{2} \, \IF_{\GP}^{-1}(\thetav) \bigr\} 
    =
    \tr \bigl\{ \HPc^{2} \bigl( \IF(\thetav) + \GP^{2} \bigr)^{-1} \bigr\}.
\label{dimGtdef}
\end{EQA}
Our results involve two particular values:
\begin{EQA}
	\dimA_{\GP}
	\eqdef
	\dimA_{\GP}(\thetavs_{\GP})
	=
	\dimA_{\GP} \, ,
	& \quad &
	\dimt_{\GP}
	\eqdef
	\dimA_{\GP}(\tilde{\thetav}_{\GP}) .
\label{dGtsttttsGte}
\end{EQA}
Usually all the values \( \dimA_{\GP} \) and \( \dimA_{\GP}(\thetav) \) for \( \thetav \) close to 
\( \thetavs_{\GP} \) are close to each other; see \eqref{DPGPm1Cd3rG} below. 
If \( \GP^{2} \equiv 0 \) and \( \HPc^{2} \approx \IF(\thetavs) \), one obviously gets
\( \dimA_{\GP}(\thetav) \asymp \dimp \).

In what follows we will consider two non-trivial examples of priors.
A \emph{truncation prior} assumes that a growing sequence of nested linear subspaces \( \VV_{\mm} \) of dimension 
\( m \) is given, and the prior distribution is supported to a \( \mm \)-dimensional subspace 
\( \VV_{\mm} \) of \( \R^{\dimp} \).
This formally corresponds to a covariance operator \( \GP_{\mm}^{-2} \) with 
\( \GP_{\mm}^{-2} \bigl( \Id - \Proj_{\VV_{m}} \bigr) = 0 \), where \( \Proj_{\VV} \) means 
the projector on \( \VV \).

\begin{lemma}
\label{Ltruncprior}
Suppose that \( \VV_{\mm} \) is a linear subspace
of \( \R^{\dimp} \) with \( \dim(\VV_{\mm}) = \mm \) and the prior \( \ND(0,\GP_{\mm}^{-2}) \)
is supported to \( \VV_{\mm} \).
Then \( \dimA_{\GP}(\thetav) \leq \mm \).
Moreover, if
\begin{EQA}
	\| \GP_{\mm} \uv \|^{2}
	& \leq &
	\gp_{\mm}^{2} \| \uv \|^{2}, 
	\qquad 
	\uv \in \VV_{\mm} \, ,
\label{uinVmgm2Gu2}
\end{EQA}
and if \( \IF(\thetav) \) fulfills 
\begin{EQ}[rccclcl]
	\CONST_{1,\IF} \, \nsize \| \uv \|^{2} 
	& \leq &
	\uv^{\T} \IF(\thetav) \uv
	& \leq &
	\CONST_{2,\IF} \, \nsize \| \uv \|^{2},
	& \qquad &
	\uv \in \VV_{\mm} \, ,
\label{uVmGu2C1C2}
\end{EQ}
and \( \HPc^{2} \geq \IF(\thetav) \), then 
\begin{EQA}
	\frac{\CONST_{1,\IF} \, n}{\CONST_{2,\IF} \, \nsize + \gp_{\mm}^{2}} \, \mm
	& \leq &
	\dimA_{\GP}(\thetav)
	\leq 
	\mm .
\label{fC1nC2ngm2mpGm}
\end{EQA}
If \( \gp_{\mm}^{2} \ll \nsize \), then \( \dimA_{\GP}(\thetav) \asymp \mm \).
\end{lemma}

A \emph{smoothing prior} corresponds to the situation when 
\( \| \GP \uv \|^{2} \) becomes large for \( \uv \not\in \VV_{\mm} \) and \( \mm \) large.
Usually one assumes that \( \VV_{\mm} \) is spanned by the eigenvectors
of \( \GP^{2} \) corresponding to its smallest eigenvalues 
\( \gp_{1}^{2} \leq \gp_{2}^{2} \leq \ldots \leq \gp_{\mm}^{2} \).
One can write this condition in the form
\begin{EQ}[rcl]
	\| \GP \uv \|^{2}
	& \leq &
	\gp_{\mm}^{2} \| \uv \|^{2} ,
	\qquad 
	\uv \in \VV_{\mm} \, ,
	\\
	\| \GP \uv \|^{2}
	& \geq &
	\gp_{\mm}^{2} \| \uv \|^{2} \, ,
	\qquad 
	\uv \in \VV_{\mm}^{c} \, ,
\label{uniVmgm2u2uiVm}
\end{EQ}
where \( \VV_{\mm}^{c} \) is the orthogonal complement of \( \VV_{\mm} \).
Further we assume that \( \gp_{j}^{2} \) grow polynomially yielding
for some \( \CONST_{1,\gp} , \CONST_{2,\gp} \) and each \( J \)
\begin{EQA}
	\CONST_{1,\gp}
	\leq 
	\frac{1}{J \gp_{J}^{-2}} \sum_{j \geq J} \gp_{j}^{-2}
	& \leq &
	\CONST_{2,\gp} \,  .
\label{sumjJgjm2C}
\end{EQA}	
A typical example is given by \( \gp_{j}^{2} = \CGP j^{2\smp} \) for \( \smp > 1/2 \)
and some window parameter \( \CGP \).

\begin{lemma}
\label{Lsmoothprior}
Let the matrices \( \IF(\thetav) \) and \( \HPc^{2} \) satisfy \eqref{uVmGu2C1C2} and let 
the prior precision matrix \( \GP^{2} \) follow \eqref{uniVmgm2u2uiVm} and \eqref{sumjJgjm2C}.
Define \( \mm \) by \( \gp_{\mm}^{2} \approx \nsize \).
Then 
\begin{EQA}
	\CONST_{3} \, \mm
	\leq 
	\dimA_{\GP}(\thetav)
	& \leq &
	\CONST_{4} \, \mm ,
\label{C4mC3mpGAtm}
\end{EQA}
where \( \CONST_{3} \) and \( \CONST_{4} \) only depend on \( \CONST_{1,\IF}, \CONST_{2,\IF} \),
and \( \CONST_{1,\gp}, \CONST_{2,\gp} \).
In particular, if \( \gp_{j}^{2} = \CGP j^{2\smp} \) with \( \smp > 1/2 \), then the effective 
dimension can be obtained as \( \dimA_{\GP}(\thetav) \asymp (n/\CGP)^{1/(2\smp)} \).
\end{lemma}

The next technical result explains how the prior can be linked to ``smoothness'' of the unknown 
vector \( \thetav \).

\begin{lemma}
\label{Lsmoothnessts}
Let \( \IF(\thetav) \) and \( \HPc^{2} \) satisfy \eqref{uVmGu2C1C2} for all \( \thetav \in \Theta \).
Suppose that \( \Theta \) is a subset of a Sobolev ball \( \BBB_{\smp}(1) \)
\begin{EQA}
	\BBB_{\smp}(1)
	& \eqdef &
	\Bigl\{ \thetav = (\theta_{j}) \colon \sum_{j \geq 1} j^{2\smp} \theta_{j}^{2} \leq 1 \Bigr\} .
\label{BBBsttsjj1sj2}
\end{EQA}
If \( \GP^{2} = \diag\bigl( \gp_{j}^{2} \bigr) \) with \( \gp_{j}^{2} = \CGP j^{2\smp} \) and \( \smp > 1/2 \), 
then \( \| \GP \thetav \|^{2} \leq \CGP \).
Define \( \CGP \) by the bias-variance trade-off
\begin{EQA}
	(n/\CGP)^{1/(2\smp)}
	\asymp 
	\CGP
	& \text{ or } &
	\CGP^{2\smp+1}
	\asymp 
	n .
\label{mgsm2nn2}
\end{EQA}
Then \( \CGP \asymp \nsize^{1/(2\smp+1)} \),
\( \dimA_{\GP}(\thetav) \lesssim \nsize^{1/(2\smp+1)} \) and 
\( \| \GP \thetav \|^{2} \lesssim \nsize^{1/(2\smp+1)} \)
for all \( \thetav \in \BBB_{\smp}(1) \).
\end{lemma}

Here and below ``\( a \lesssim b \)'' means \( a \leq \CONST b \) with an 
absolute constant \( \CONST \).
%
In what follows we implicitly 
assume that each value \( \dimA_{\GP}(\thetav) \) is much smaller than the full dimension \( \dimp \) which can be even infinite.
Most of our results requires \( \dimA_{\GP}(\thetav) \ll \nsize^{1/3} \),
that is, \( \smp > 1 \) in case of Lemma~\ref{Lsmoothnessts}.

\subsection{Properties of the pMLE \( \tilde{\thetav}_{\GP} \)}
\label{StiteGpro}
This section briefly reviews some properties of the penalized MLE \( \tilde{\thetav}_{\GP} = \argmax \LGP(\thetav) \).
Our results are based on conditions \nameref{LLref}, \nameref{Eref}, \nameref{ED0ref}, \nameref{HGref}, and \nameref{LL0ref} even if not mentioned explicitly.
In particular, we systematically use that the stochastic term in the log-likelihood only linearly depends on \( \thetav \)
and that the expected log-likelihood is concave in \( \thetav \).
We state two results, the first one claims a kind of local concentration of the penalized MLE \( \tilde{\thetav}_{\GP} \),
while the second one describes some useful expansions for the estimator \( \tilde{\thetav}_{\GP} \) and for the fitted log-likelihood \( \LGP(\tilde{\thetav}_{\GP}) \).
The presented results substantially improve similar statements in \cite{SP2013_rough}.
  
Remind that \( \thetavs_{\GP} = \arginf_{\thetav} \E \LGP(\thetav) \) and \( \DPGP^{2} =  \IF(\thetavs_{\GP}) + \GP^{2} \).
Below we show that the penalized MLE \( \tilde{\thetav}_{\GP} \) concentrates with a high probability on the elliptic set 
\begin{EQA}
    \CAGP(\rr)
    & \eqdef &
    \bigl\{ \thetav \colon \| \DPGP (\thetav - \thetavs_{\GP}) \| \leq \rr \bigr\}
\label{DG2Dts2G2Trr}
\end{EQA}
under a proper choice of \( \rr \).

As the stochastic component of \( \LGP(\thetav) \) is linear in \( \thetav \),
the gradient \( \nabla \zeta = \nabla \bigl\{ \LGP(\thetav) - \E \LGP(\thetav) \bigr\} \)
does not depend on \( \thetav \).  
Under condition \nameref{ED0ref}, there exists a random set \( \Omega(\xx) \) with 
\( \P\bigl( \Omega(\xx) \bigr) \geq 1 - \CONST \ex^{-\xx} \) such that on this set
\( \bigl\| \DPGP^{-1} \nabla \zeta \bigr\| \leq \zq(\BB_{\GP},\xx) \):
\begin{EQA}
    \bigl\| \DPGP^{-1} \nabla \zeta \bigr\|
    & \leq &
    \zq(\BB_{\GP},\xx) 
    \quad
    \text{on }
    \Omega(\xx)
    \text{ with }
    \,
    \P\bigl( \Omega(\xx) \bigr) \geq 1 - \CONST \ex^{-\xx},
\label{uTzDGvTxm12z}
\end{EQA}
where \( \BB_{\GP} = \HPc \, \DPGP^{-2} \, \HPc \) and \( \zq(\BB_{\GP},\xx) \) is given by \eqref{zzxxppdBlroB};
see Theorem~\ref{TdevboundGauss} with \( \xiv = \nabla \zeta \).
One can use the simplified bound 
\begin{EQA}
	\zq(\BB_{\GP},\xx) 
	& \leq &
	\sqrt{\tr (\BB_{\GP})} + \sqrt{2\xx \, \| \BB_{\GP} \|} .
\label{zqBVGlstBVG2}
\end{EQA}
It is worth mentioning that this deviation bound is the only place where the stochastic nature of the log-likelihood \( L(\thetav) \) is accounted for. 
In the rest, we only use the condition \nameref{Eref} about linearity the stochastic component \( \zeta(\thetav) \) in 
\( \thetav \).

Our first result describes the concentration properties of the penalized MLE \( \tilde{\thetav}_{\GP} \).

\begin{proposition}
\label{TconstGrDG}
Assume \eqref{uTzDGvTxm12z}.
Let also \( \rr_{\GP} \) be such that \( \Thetad \) contain the set 
\( \CAGP(\rr_{\GP}) \eqdef 
\bigl\{ \thetav \colon \| \DPGP (\thetav - \thetavs_{\GP}) \| \leq \rr_{\GP} \bigr\} \), and with
\( \dltwu_{3,\HPc} = \dltwu_{3,\HPc}(\rr_{\GP}) \) from \eqref{dmtuCmnumddL0}, it holds
\begin{EQ}[rcccl]
    3 \rr_{\GP} \, \dltwu_{3,\HPc}
    & \leq &
    \rho \leq 1/2, 
    \qquad
    (1 - \rho) \rr_{\GP}
    & \geq &
    \zq(\BB_{\GP},\xx) .
\label{3rGm2d3rG12}
\end{EQ}
Then on \( \Omega(\xx) \), the estimate 
\( \tilde{\thetav}_{\GP} \) belongs to this set \( \CAGP(\rr_{\GP}) \) as well, that is, 
\begin{EQA}
    \bigl\| \DPGP \bigl( \tilde{\thetav}_{\GP} - \thetavs_{\GP} \bigr) \bigr\|
    & \leq &
    \rr_{\GP} \, .
\label{DPttGtsGrG}
\end{EQA}
\end{proposition}
  
\begin{remark}
In words, \eqref{DPttGtsGrG} means that \( \tilde{\thetav}_{\GP} \) belongs with a high probability
to the vicinity \( \CAGP(\rr_{\GP}) \) from \eqref{DG2Dts2G2Trr} with \( \rr_{\GP} \leq 2 \zq(\BB_{\GP},\xx) \).
Under \eqref{dmtuCmnumdd} \( \dltwu_{3,\HPc} \lesssim \nsize^{-1/2} \)
while 
\( \zq^{2}(\BB_{\GP},\xx) \asymp \dimA_{\GP} = \tr (\BB_{\GP}) \);
see \eqref{zqBVGlstBVG2} and examples in Section~\ref{SexamplesBvM}.
Therefore, \( \rr_{\GP} \, \dltwu_{3,\HPc} \asymp (\dimA_{\GP}/ \nsize)^{1/2} \),
and conditions \eqref{3rGm2d3rG12} require only that the value
\( \dimA_{\GP} \) is smaller in order
than the sample size \( \nsize \), i.e. \( \dimA_{\GP} \ll \nsize \). 
\end{remark}

Due to the concentration result of Proposition~\ref{TconstGrDG}, the estimate 
\( \tilde{\thetav}_{\GP} \) lies with a dominating probability in a local vicinity 
of the point \( \thetavs_{\GP} \).
Now one can use a quadratic approximation for the penalized log-likelihood process 
\( \LGP(\thetav) \) to establish an expansion for the penalized MLE
\( \tilde{\thetav}_{\GP} \) and for the excess \( \LGP(\tilde{\thetav}_{\GP}) - \LGP(\thetavs_{\GP}) \).

\begin{theorem}
\label{TFiWititG}
Under the conditions of Proposition~\ref{TconstGrDG}, it holds on \( \Omega(\xx) \) 
\begin{EQA}
    \bigl\| \DPGP \bigl( \tilde{\thetav}_{\GP} - \thetavs_{\GP} \bigr) - \DPGP^{-1} \nabla \zeta \bigr\|^{2}
    & \leq &
    4 \rr_{\GP}^{3} \, \dltwu_{3,\HPc} .
\label{DGttGtsGDGm13rG}
    \\
    \biggl| 
    \LGP(\tilde{\thetav}_{\GP}) - \LGP(\thetavs_{\GP}) 
    - \frac{1}{2} \bigl\| \DPGP^{-1} \nabla \zeta \bigr\|^{2}
    \biggr|
    & \leq &
    \rr_{\GP}^{3} \, \dltwu_{3,\HPc} .
\label{3d3Af12DGttG}
\end{EQA}
Also 
\begin{EQA}
    \biggl| 
    	\LGP(\tilde{\thetav}_{\GP}) - \LGP(\thetavs_{\GP}) 
    	- \frac{1}{2} \bigl\| \DPGP (\tilde{\thetav}_{\GP} - \thetavs_{\GP}) \bigr\|^{2}
    \biggr|
    & \leq &
    \rr_{\GP}^{3} \, \dltwu_{3,\HPc}  ,
    \\
    \sup_{\thetav \in \CAGP(\rr_{\GP})}
    \biggl| 
    	\LGP(\tilde{\thetav}_{\GP}) - \LGP(\thetav) 
    	- \frac{1}{2} \bigl\| \DPGPt (\tilde{\thetav}_{\GP} - \thetav) \bigr\|^{2}
    \biggr|
    & \leq &
    \rr_{\GP}^{3} \, \dltwu_{3,\HPc} , 
\label{3d3Af12DGttt}
\end{EQA}
where the random matrix \( \DPGPt^{2} = \IF_{\GP}(\tilde{\thetav}_{\GP}) \) fulfills
on \( \Omega(\xx) \) for some universal constant \( \CONST \)
\begin{EQA}
    \bigl\| \DPGP^{-1} \bigl( \DPGPt^{2} - \DPGP^{2} \bigr) \DPGP^{-1} \bigr\|
    & \leq &
    \CONST \rr_{\GP} \, \dltwu_{3,\HPc} .
\label{DPGPm1Cd3rG}
\end{EQA}
\end{theorem}

\begin{remark}
The results of Theorem~\ref{TFiWititG} can be viewed as finite sample nonparametric analogs
of classical asymptotic parametric results such as Fisher expansion for the MLE and Wilks 
phenomenon. 
In fact, all the mentioned classical results can be easily derived from \eqref{DGttGtsGDGm13rG} through \eqref{3d3Af12DGttt} provided asymptotic normality of the normalized score 
\( \DPGP^{-1} \nabla \zeta \); see Theorem~\ref{TestlosspMLE} below.
In particular, one can state root-\( \nsize \) normality of the pMLE and a generalized chi-squared limit
distribution for the excess \( \LGP(\tilde{\thetav}_{\GP}) - \LGP(\thetavs_{\GP}) \).
\end{remark}

\begin{remark}
Similarly to Proposition~\ref{TconstGrDG}, the results of Theorem~\ref{TFiWititG} are 
meaningful if \( \dimA_{\GP} \) is significantly smaller than \( \nsize \).
\end{remark}

The concentration set \( \CAGP(\rr_{\GP}) \) becomes smaller when \( \GP^{2} \) increases. 
In particular, if \( \GP^{2} \) is large then \( \tilde{\thetav}_{\GP} \) concentrates 
on a small vicinity of \( \thetavs_{\GP} \).
At the same time, 
penalization \( \bigl\| \GP \thetav \bigr\|^{2} \) yields some estimation bias measured by 
\( \E \LGP(\thetavs_{\GP}) - \E \LGP(\thetavs) \) and \( \thetavs_{\GP} - \thetavs \).
The bias is not critical if the true value \( \thetavs \) is ``smooth'', that is, \( \| \GP \thetavs \|^{2} \) is not too big.
The next result makes these statements precise.

\begin{proposition}
\label{TbiasGP}
It holds
\begin{EQA}
	\E \LGP(\thetavs_{\GP}) - \E \LGP(\thetavs)
	& \leq &
	\frac{1}{2} \| \GP \thetavs \|^{2} .
\label{122GttsELG}
\end{EQA}
If, in addition, 
\( \| \GP \thetavs \|^{2} \leq \rrbias^{2}/2 \) for some \( \rrbias \) such that 
\( \rrbias \, \dltwu_{3,\HPc}(\rrbias) \leq 1/2 \), then 
\begin{EQA}
\label{fr33GrGd3rELG21}
	\Bigl| 
		\E \LGP(\thetavs_{\GP}) - \E \LGP(\thetavs) - \bigl\| \DPGP (\thetavs_{\GP} - \thetavs) \bigr\|^{2}/2 
	\Bigr|
	& \leq &
	\rrbias^{3} \, \dltwu_{3,\HPc}(\rrbias),
	\\
	\bigl\| \DPGP (\thetavs - \thetavs_{\GP}) - \DPGP^{-1} \GP^{2} \thetavs \bigr\|^{2}
	& \leq &
	4 \rrbias^{3} \, \dltwu_{3,\HPc}(\rrbias) \, .
	\qquad
\label{fr33GrGd3rELG22}
\end{EQA}
Moreover, for any linear mapping \( \QP \) in \( \R^{\dimp} \) 
it holds
\begin{EQA}
	\| \QP (\thetavs_{\GP} - \thetavs) \|
	& \leq &
	\| \QP \DPGP^{-2} \GP^{2} \thetavs \|
	+ 2 \sqrt{\| \QP \DPGP^{-2} \QP^{\T} \| \, \rrbias^{3} \, \dltwu_{3,\HPc}(\rrbias)} \, .
\label{22GtsQDGm2}
\end{EQA}
\end{proposition}

%

Finally, we can combine all the previous results together to bound the loss 
\( \| \QP (\tilde{\thetav}_{\GP} - \thetavs) \| = \| \QP (\tilde{\thetav}_{\GP} - \thetavs) \| 
= \| \QP (\tilde{\thetav}_{\GP} - \thetavs_{\GP}) + \QP (\thetav_{\GP} - \thetavs) \| \). 
One can apply \( \QP = \sqrt{\IF(\thetavs)} \) for prediction and \( \QP = \sqrt{n} \, \Id_{\dimp} \)
for estimation.
%
The first bound 
can be viewed as analog of classical bias-variance decomposition of the loss 
\( \bigl\| \QP (\tilde{\thetav}_{\GP} - \thetavs) \bigr\| \).
The second one is asymptotic and it corresponds to the case of 
``small bias'' or ``undersmoothing''. 
Below \( o(1) \) means a small asymptotically vanishing value,
and  ``\( a \lesssim b \)'' means \( a \leq \CONST b \) with an 
absolute constant \( \CONST \) that possibly depends on the constants from our conditions.
We also ignore the small error term in \eqref{22GtsQDGm2}.

\begin{theorem}
\label{TestlosspMLE}
Given \( \QP \) with \( \QP^{\T} \QP \leq \DPGP^{2} \), 
define \( \BB_{\QPGP} = \QP \DPGP^{-2} \HPc^{2} \DPGP^{-2} \QP^{\T} \). \\
(1) On a random set \( \Omega_{1}(\xx) \) with \( \P\bigl( \Omega_{1}(\xx) \bigr) \leq 2 \ex^{-\xx} \),
it holds
\begin{EQA}
	\| \QP (\tilde{\thetav}_{\GP} - \thetavs) \|
	& \leq &
	\| \QP \DPGP^{-2} \GP^{2} \thetavs \|
	+ \zq(\BB_{\QPGP},\xx) \, ,
\label{22GtsQDGm2loss}
\end{EQA}
where \( \zq(\BB,\xx) \leq \sqrt{\tr \BB} + \sqrt{2 \| \BB \| \, \xx} \); see \eqref{zqBVGlstBVG2}.
Under the bias-variance trade-off
\begin{EQA}
	\| \QP \DPGP^{-2} \GP^{2} \thetavs \|^{2}
	& \lesssim &
	\tr\bigl( \BB_{\QPGP} \bigr) ,
\label{QDGm2G2ts2trV2}
\end{EQA}
one obtains on \( \Omega(\xx) \)
\begin{EQA}
	\| \QP (\tilde{\thetav}_{\GP} - \thetavs) \|^{2}
	& \lesssim &
	\tr (\BB_{\QPGP}) + \| \BB_{\QPGP} \| \, \xx \, .
\label{QttGts2trVQGx}
\end{EQA}
(2) Let \( \nabla \zeta \) be nearly normal in the sense that 
\begin{EQA}
	\sup_{\av \in \R^{\dimp}} \, \sup_{\zq > 0} \Bigl| 
		\P\bigl( \bigl\| \QP \DP_{\GP}^{-2} \nabla \zeta - \av \bigr\| \leq \zq \bigr)
		- \P\bigl( \bigl\| \VP_{\QPGP} \gammav - \av \bigr\| \leq \zq \bigr)
	\Bigr|
	&=&
	o(1)
\label{szQDGm2nzzVQG}
\end{EQA}
with \( \gammav \in \R^{\dimp} \) standard normal and 
\( \VP_{\QPGP}^{2} = \Var\bigl( \QP \DP_{\GP}^{-2} \nabla \zeta \bigr) 
= \QP \DPGP^{-2} \VP^{2} \DPGP^{-2} \QP^{\T} \).
Assume also the ``small bias'' condition 
\begin{EQA}
	\frac{\| \QP \DPGP^{-2} \GP^{2} \thetavs \|^{2}}{\tr\bigl( \VP_{\QPGP}^{2} \bigr)}
	&=&
	o(1) .
\label{QDGm2G2ts2V2}
\end{EQA} 
Then it holds
\begin{EQA}
	\sup_{\zq > 0} \Bigl| 
		\P\bigl( \bigl\| \QP (\tilde{\thetav}_{\GP} - \thetavs) \bigr\| \leq \zq \bigr)
		- \P\bigl( \bigl\| \VP_{\QPGP} \gammav \bigr\| \leq \zq \bigr)
	\Bigr|
	&=&
	o(1) .
\label{01zPPQVzzts}
\end{EQA}
\end{theorem}

\begin{remark}
Given \( \alp \), define \( \zq_{\alp} \) by 
\( \P( \| \VP_{\QPGP} \gammav \| \geq \zq_{\alp} ) = \alp \).
It follows from \eqref{01zPPQVzzts} that \( \CS_{\QPGP}(\zq_{\alp}) \eqdef 
	\bigl\{ \thetav \colon \| \QP (\tilde{\thetav}_{\GP} - \thetav) \| \leq \zq_{\alp} \bigr\} \) is asymptotically valid confidence set for \( \thetavs \).
\end{remark}

\begin{remark}
\label{RTeffdimML}
The result \eqref{QttGts2trVQGx} with \( \QP = \DPGP \) yields on \( \Omega(\xx) \) a bound
\( \| \QP (\tilde{\thetav}_{\GP} - \thetavs) \|^{2} \lesssim \dimA_{\GP} \) as in \cite{SP2013_rough}.
\end{remark}

\subsection{Posterior concentration and Gaussian approximation}
\label{Spostcontraction}
Now we turn to the properties of the posterior \( \vthetav_{\GP} \cond \Yv \).
Our first result shows that the posterior concentrates on the elliptic set 
\( \CAt_{\GP}(\rups) = \bigl\{\uv \colon \| \HPc (\uv - \tilde{\thetav}_{\GP})\| \leq \rups \bigr\} \) for a proper value 
\( \rups \geq \CONST \sqrt{\dimt_{\GP}} + \CONST \sqrt{\xx} \) for 
\( \dimt_{\GP} = \dimA_{\GP} (\tilde{\thetav}_{\GP}) \). 
For this we bound from above the random quantity
\begin{EQA}
	\rho(\rups)
	& \eqdef &
	\frac{\int \Ind\bigl( \|\HPc \uv\| > \rups \bigr)
      		\exp \bigl\{ \LGP(\tilde{\thetav}_{\GP}+\uv) \bigr\} d \uv}
   		 {\int \Ind\bigl( \|\HPc \uv\| \leq \rups \bigr) \exp \bigl\{ \LGP(\tilde{\thetav}_{\GP}+\uv) \bigr\} d \uv} \, .
\label{rhopipoprDGP}
\end{EQA}
Obviously \( \P\bigl( \| \HPc (\vthetav_{\GP} - \tilde{\thetav}_{\GP}) \| > \rups \cond \Yv \bigr) \le \rho(\rups) \).
Therefore, small values of \( \rho(\rups) \) indicate a concentration of 
\( \vthetav_{\GP} \cond \Yv \) on the set \( \CAt_{\GP}(\rups) \).



\begin{proposition}
\label{PrhoQPBvM}
Let conditions of Proposition~\ref{TconstGrDG} be satisfied. 
Assume that for some fixed values \( \rups \) and \( \xx > 0 \), it holds
\begin{EQA}
\label{LmgfquadELGP}
    \err(\rups)
    \eqdef 
    4 \rups^{6} \, \dltwu_{3,\HPc}^{2} + 4 \rups^{4} \, \dltwu_{4,\HPc}
    \leq 
    1/2,
    & \quad &
    \CONSTru
    \eqdef 
    1 - 3 \rups \, \dltwu_{3,\HPc} 
    \geq 1/2,
\end{EQA}
and also 
\begin{EQA}
    \CONSTru \rups 
    & \geq &
    2 \sqrt{\dimA_{\GP}(\thetav)} + \sqrt{\xx} \, ,
    \qquad
    \thetav \in \Thetad \, . 
\label{CONSTruAxx1}
\end{EQA}
Then,  on the random set \( \Omega(\xx) \) from \eqref{uTzDGvTxm12z}, 
\( \rho(\rups) \) from \eqref{rhopipoprDGP} fulfills with
\( \dimt_{\GP} = \dimA_{\GP}(\tilde{\thetav}_{\GP}) \)
\begin{EQA}
    \rho(\rups)
    & \leq &
    \frac{1}{1 - \err(\rups)} \,\,
    \frac{\exp\bigl\{- (\dimt_{\GP} + \xx)/2\bigr\}}
    	 {1 - \exp\{ -(\dimt_{\GP} + \xx)/2 \}} \, .
\label{rhoQPubnBvm}
\end{EQA}
\end{proposition}

\begin{remark}
Conditions \eqref{LmgfquadELGP} and \eqref{CONSTruAxx1} can be spelled out as follows:
the value \( \rups^{2} \) is of order at least the effective dimension \( \dimA_{\GP}(\thetav) \),
while the values \( \rups^{3} \, \dltwu_{3,\HPc} \) and \( \rups^{4} \, \dltwu_{4,\HPc} \) 
should be small.
If \( \dltwu_{3,\HPc} \asymp \nsize^{-1/2} \), \( \dltwu_{4,\HPc} \asymp \nsize^{-1} \), then 
\eqref{LmgfquadELGP} requires \( \dimA_{\GP}^{3}(\thetav) \ll \nsize \), \( \thetav \in \Thetad \). 
\end{remark}

\begin{remark}
It is of interest to compare the concentration result of Proposition~\ref{TconstGrDG} for the pMLE
\( \tilde{\thetav}_{\GP} \)
and the concentration bound \eqref{rhoQPubnBvm} for the posterior \( \vthetav_{\GP} \cond \Yv \).
Suppose that \( \dimA_{\GP}(\thetav) \asymp \dimA_{\GP} \) for \( \thetav \in \Thetad \).
Then also \( \rups \asymp \rr_{\GP} \).
Therefore, the concentration results for the penalized MLE \( \tilde{\thetav}_{\GP} \) 
and for the posterior \( \vthetav_{\GP} \cond \Yv \) look similar, but there is one essential difference. 
The properly shifted MLE \( \tilde{\thetav}_{\GP} \) well concentrates on a rather small elliptic set 
\( \bigl\{ \uv \colon \| \DPGP (\uv - \thetavs_{\GP}) \| \leq \rr_{\GP} \bigr\} \)
centered at \( \tilde{\thetav}_{\GP} \).
In other words, \( \DPGP \bigl( \tilde{\thetav}_{\GP} - \thetavs_{\GP} \bigr) \) belongs to the ball in \( \R^{\dimp} \) 
of radius \( \rr_{\GP} \) with a high probability.
This holds true even if \( \dimp = \infty \).
Our result of Proposition~\ref{PrhoQPBvM} claims concentration 
of \( \vthetav_{\GP} \cond \Yv \)
on a larger set \( \CAt_{\GP}(\rups) \), also with an elliptic shape, but 
centered at \( \tilde{\thetav}_{\GP} \), and with 
larger axes corresponding to \( \HPc^{-1} \) instead of \( \DPGP^{-1} \).
Later we will see that 
\( \DPGP \bigl( \vthetav_{\GP} - \tilde{\thetav}_{\GP} \bigr) \cond \Yv \) is close to the standard Gaussian measure and it does not concentrate on a ball in \( \R^{\infty} \) for any radius \( \rr \). 
\end{remark}

Our main result claims that the posterior can be well approximated
by a Gaussian distribution \( \ND\bigl( \tilde{\thetav}_{\GP}, \DPGPt^{-2} \bigr) \).
By \( \PG \) we denote a standard normal distribution of a random vector 
\( \gammav \in \R^{\dimp} \) given \( \DPGPt = \DPGP(\tilde{\thetav}_{\GP}) \). 
%
We distinguish between the class \( \cc{B}_{\smp}(\R^{\dimp}) \) of centrally symmetric Borel sets and the class \( \cc{B}(\R^{\dimp}) \) of all Borel sets in \( \R^{\dimp} \).

\begin{theorem}
\label{TnonparBvMm}
Let the conditions of Proposition~\ref{PrhoQPBvM} hold and \( \rho(\rups) \) satisfy~\eqref{rhoQPubnBvm}.
It holds on the set \( \Omega(\xx) \) from \eqref{uTzDGvTxm12z} 
for any centrally symmetric Borel set \( A \in \cc{B}_{\smp}(\R^{\dimp}) \) 
\begin{EQA}
    \P\bigl( \vthetav_{\GP} - \tilde{\thetav}_{\GP} \in A \cond \Yv \bigr)
    & \geq &
    \frac{1 - \err(\rups)}
     	 {\bigl\{ 1 + \err(\rups) + \rho(\rups) \bigr\}} \,
    \PG\bigl( \DPGPt^{-1} \gammav \in A \bigr) - \rho(\rups) \, ,
    \\
    \P\bigl( \vthetav_{\GP} - \tilde{\thetav}_{\GP} \in A \cond \Yv \bigr)
    & \leq & 
    \frac{1 + \err(\rups)}
     {\bigl\{ 1 - \err(\rups) \bigr\} \bigl( 1 - \ex^{-\xx} \bigr)} 
    \PG\bigl( \DPGPt^{-1} \gammav \in A \bigr)
     + \rho(\rups) \, .
\label{1p1m1a1emx}
\end{EQA}
For any measurable set \( A \in \cc{B}(\R^{\dimp}) \), similar bounds 
hold with \( \dltwu_{3,\HPc} \) in place of \( \err(\rups) \).
\end{theorem}

\noindent
The first result of the theorem states for any symmetric set \( A \in \cc{B}_{\smp}(\R^{\dimp}) \) 
\begin{EQA}
    \left| 
    {\P\bigl( \vthetav_{\GP} - \tilde{\thetav}_{\GP} \in A \cond \Yv \bigr)}
    - {\PG\bigl( \DPGPt^{-1} \gammav \in A \bigr)} 
    \right|
    & \lesssim &
    {\PG\bigl( \DPGPt^{-1} \gammav \in A \bigr)}
    \bigl\{ \err(\rups) + \ex^{-\xx} \bigr\} + \rho(\rups).
\label{PvtGttGAcY}
\end{EQA}
The second statement applies to any \( A \in \cc{B}(\R^{\dimp}) \) and hence, 
it bounds the distance in total variation between the posterior and its Gaussian approximation
\( \DPGPt^{-1} \gammav \).

\begin{corollary}
\label{CnonparBvMm}
Suppose that \( \rups \) satisfies the conditions \eqref{LmgfquadELGP} and \eqref{CONSTruAxx1} with 
\( \xx = 2 \log \nsize \).
It holds on \( \Omega(\xx) \) 
\begin{EQA}
\label{1p1m1a1emxrr}
	\sup_{A \in \cc{B}_{\smp}(\R^{\dimp})}
	\left| 
		\P\bigl( \vthetav_{\GP} - \tilde{\thetav}_{\GP} \in A \cond \Yv \bigr)
		- \PG\bigl( \DPGPt^{-1} \gammav \in A \bigr)  
	\right|
	& \lesssim &
	\err(\rups) 
	+ 1/n ,
	\\
\label{1p1m1a1emxrras}
	\sup_{A \in \cc{B}(\R^{\dimp})}
	\left| 
		\P\bigl( \vthetav_{\GP} - \tilde{\thetav}_{\GP} \in A \cond \Yv \bigr)
		- \PG\bigl( \DPGPt^{-1} \gammav \in A \bigr)  
	\right|
	& \lesssim &
	\rups^{3} \, \dltwu_{3,\HPc} 
	+ 1/n  .
\end{EQA}
\end{corollary}

Comparison of two bounds of Corollary~\ref{CnonparBvMm} reveals
that the use of symmetric credible sets improves the accuracy of Gaussian approximation
from \( \rups^{3} \, \dltwu_{3,\HPc} \) to 
\( \err(\rups) \asymp \rups^{6} \, \dltwu_{3,\HPc}^{2} + \rups^{4} \, \dltwu_{4,\HPc} \).
In particular, under \eqref{dmtuCmnumdd} and \eqref{uVmGu2C1C2},
\( \dltwu_{3,\HPc} \asymp \rups^{3} \, \nsize^{-1/2} \) while 
\( \err(\rups) \asymp \rups^{6}/n \).
The choice \( \xx = 2 \log \nsize \) and 
\( \rups = \CONST \bigl( \sqrt{\dimA_{\GP}} + \sqrt{\log n} \bigr) \) yields 
\( \rho(\rups) \leq 1/n \),
and the leading term in the error of Gaussian approximation \eqref{1p1m1a1emxrr} is 
\( \err(\rups) \asymp \dimA_{\GP}^{6}/n \).
The bound \eqref{1p1m1a1emxrras} in TV-distance ensures an error of order
\( \sqrt{\dimA_{\GP}^{3}/n} \).


\smallskip
Now we turn to \emph{contraction property}
describing the distance between the support of the posterior and the true value \( \thetavs \).
For a given linear mapping from \( \R^{\dimp} \) satisfying \( \QP^{\T} \QP \leq \DPGP^{2} \), 
we would like to describe a minimal radius \( \rr \) ensuring 
\begin{EQA}
	\P\bigl( \| \QP (\vthetav_{\GP} - \thetavs) \| > \rr \cond \Yv \bigr)
	&=&
	o(1) . 
\label{PQvtGtslrYo1}
\end{EQA}

\begin{theorem}
\label{Ccontactionrate}
Assume that \( \QP^{\T} \QP \leq \DPGP^{2} \), \( \HPc^{2} \lesssim \DPGP^{2} \), and
\begin{EQA}
	\| \QP \DPGP^{-2} \GP^{2} \thetavs \|^{2} 
	& \lesssim &
	\tr \bigl( \QP \DPGP^{-2} \QP^{\T} \bigr) .
\label{trQDGm2G22Q2DG2}
\end{EQA}
Then it holds on \( \Omega(\xx) \) for some fixed 
\( \CONST, \CONST_{1}, \CONST_{2} \)
\begin{EQA}
	\P\Bigl( 
		\bigl\| \QP (\vthetav_{\GP} - \thetavs) \bigr\|^{2} 
		\geq \CONST_{1} \tr \bigl( \QP \DPGP^{-2} \QP^{\T} \bigr) 
		+ \CONST_{2} \log \nsize \, \| \QP \DPGP^{-2} \QP^{\T} \| 
		\, \cond \Yv 
	\Bigr)
	& \leq &
	\CONST \, \nsize^{-1} .
	\qquad
\label{PDvtGtsCrGY}
\end{EQA}
\end{theorem}

\begin{remark}
If \( \QP = \HPc \), then \( \tr \bigl( \QP \DPGP^{-2} \QP^{\T} \bigr) = \dimA_{\GP} \), 
and the contraction radius squared is proportional to the effective dimension. 
The relation \eqref{trQDGm2G22Q2DG2} is called ``bias-variance trade-off'',
and it means that the squared bias \( \bigl\| \QP \bigl( \thetavs_{\GP} - \thetavs \bigr) \bigr\|^{2} \) is not larger in order than the trace of the variance of the posterior
\( \QP \vthetav_{\GP} \cond \Yv \). 
The bound \eqref{PDvtGtsCrGY} on posterior contraction is sharp in the sense that it cannot be improved even in the Gaussian case.
\end{remark}

Our results help to address the issue of \emph{frequentist coverage} of Bayesian credible sets.
Corollary~\ref{CnonparBvMm} suggests to consider credible sets of the form
\begin{EQA}
	\CAt_{\QPGP}(\rr)
	& \eqdef &
	\bigl\{ \thetav \colon \bigl\| \QP (\tilde{\thetav}_{\GP} - \thetav) \bigr\| \leq \rr \bigr\} ,
\label{SAQPrtcQttr}
\end{EQA}
where \( \QP^{\T} \QP \leq \DPGP^{2} \) and \( \rr = \rr_{\alp} \) is fixed to ensure 
\begin{EQA}
	\PG\bigl( \bigl\| \QP \DPGPt^{-1} \gammav \bigr\| > \rr_{\alp} \bigr)
	& = &
	\alp ,
	\qquad
	\gammav \sim \ND(0,\Id_{\dimp}) .
\label{alPGQDtm1gra}
\end{EQA}
It appears that frequentist validity of such credible sets can be stated under the same 
``small bias'' condition \eqref{QDGm2G2ts2V2} as for pMLE-based confidence sets; see Theorem~\ref{TestlosspMLE}.

\begin{theorem}
\label{CThonestCS}
Suppose that 
\( \VP^{2} \leq \DPGP^{2} \) and moreover, 
the score \( \nabla \zeta \) is nearly normal as in \eqref{szQDGm2nzzVQG}.
Assume also ``small bias'' condition \eqref{QDGm2G2ts2V2}. 
Then 
\begin{EQA}
	\P\bigl( \thetavs \not\in \CAt_{\QPGP}(\rr_{\alp}) \bigr)
	& \leq &
	\alp + o(1) . 
\label{Palpo1AGra}
\end{EQA}
\end{theorem}

In general, 
the coverage probability of \( \CAt_{\QPGP}(\rr_{\alp}) \) is larger than 
the nominal probability \( \alp \).
The reason is that the prior covariance \( \DPGPt^{-2} \) is larger than the covariance 
\( \DPGP^{-2} \VPc^{2} \DPGP^{-2} \)
of the pMLE \( \tilde{\thetav}_{\GP} \).
%

\subsection{Some extensions}
\label{Ssomeextensions}
Here we list some possible straightforward extensions of the results presented above.

\subsubsection{The use of posterior mean}
The main results so far used the pMLE \( \tilde{\thetav}_{\GP} \) for centering the posterior.
In practical applications, computing the pMLE or equivalently, MAP \( \tilde{\thetav}_{\GP} \)
could be a hard task.
A natural question here is whether one can instead use 
the posterior mean \( \bar{\vthetav}_{\GP} = \E \bigl( \vthetav_{\GP} \cond \Yv \bigr) \) 
which can be efficiently computed by Monte Carlo simulations. 
One can easily show that the posterior mean \( \bar{\vthetav}_{\GP} \) is very close to the pMLE.
Unfortunately, our most advanced results on \( 1/n \) accuracy of Gaussian approximation
of the posterior would not apply after centering by \( \bar{\vthetav}_{\GP} \) because
the symmetricity arguments do not apply.
This section states some sufficient conditions under which we still get the desirable 
\( 1/n \) accuracy even after centering by \( \bar{\vthetav}_{\GP} \).


We start 
with the bounds on the first two moments of the posterior. 
The results claim that posterior mean and variance are close to that of normal 
\( \ND\bigl( \tilde{\thetav}_{\GP}, \DPGPt^{-2} \bigr) \).
Let \( \bar{\vthetav}_{\GP} = \E \bigl( \vthetav_{\GP} \cond \Yv \bigr) \) and 
\( \av = \bar{\vthetav}_{\GP} - \tilde{\thetav}_{\GP} \). 
Note that symmetricity arguments do not apply to the posterior mean,
therefore one can expect an accuracy of order \( \| \av \| \asymp \rups^{3} \, \delta_{3,\HPc} \asymp \sqrt{\dimA_{\GP}^{3}/n} \).
In the contrary, the posterior covariance can be estimated with a higher accuracy of order 
\( \err(\rups) \), again by symmetricity arguments. 
The next results describes the moments of 
\( \QP \bigl( \vthetav_{\GP} - \tilde{\thetav}_{\GP} \bigr) \) for
an arbitrary linear operator \( \QP \) with \( \QP^{\T} \QP \leq \DPGP^{2} \).

\begin{theorem}
\label{Tpostmeanv}
Let the conditions of Theorem~\ref{PrhoQPBvM} hold.
Let \nameref{LL0ref} be fulfilled with \( \rups \) 
satisfying \eqref{CONSTruAxx1} with \( \sqrt{\dimA_{\GP}(\thetav) + 1} \)
in place of \( \sqrt{\dimA_{\GP}(\thetav)} \), and \( \err(\rups) \) be given by \eqref{LmgfquadELGP}
and \eqref{CONSTruAxx1} with \( \xx = 2 \log n \).
Then \( \bar{\vthetav}_{\GP} = \E \bigl( \vthetav_{\GP} \cond \Yv \bigr) \)
fulfills on the set \( \Omega(\xx) \) for any linear operator \( \QP \) 
\begin{EQA}
    \bigl\| 
    	\QP \bigl( \bar{\vthetav}_{\GP} - \tilde{\thetav}_{\GP} \bigr)   
    \bigr\|
    & \lesssim &
    \delta_{3}(\rups) \sqrt{\dimt_{\QPGP}} 
    + n^{-1} \, \dimt_{\QPGP} \, ,
\label{DGbvtGttGd3r0}
\end{EQA}
where 
\( \dimt_{\QPGP} = \tr (\QP \DPGPt^{-2} \QP^{\T}) \).
The posterior variance \( \Var \bigl( \vthetav_{\GP} \cond \Yv \bigr) \) fulfills
on \( \Omega(\xx) \)
\begin{EQA}
    \bigl\| \Id - \DPGPt \Var \bigl( \vthetav_{\GP} \cond \Yv \bigr) \DPGPt \bigr\|
    =	
    \sup_{\zv \in \Sphere_{\dimp}}
    \biggl| 
    \E \Bigl[ 
      \bigl\langle \zv, \DPGPt \bigl( \vthetav_{\GP} - \tilde{\thetav}_{\GP} \bigr) \bigr\rangle^{2} 
    \cond \Yv 
    \Bigr] - 1 
    \biggr|
    & \lesssim &
    \err(\rups) ,
    \\
    \biggl| 
    \E \Bigl(
    \bigl\|  
      \QP \bigl( \vthetav_{\GP} - \tilde{\thetav}_{\GP} \bigr) 
    \bigr \|^{2} 
    \cond \Yv 
    \Bigr) - \dimt_{\QPGP} 
    \biggr|
    & \lesssim &
    \err(\rups) \, \dimt_{\QPGP} \, .
\label{bBbDGtvtGtG2u2}
\end{EQA}
\end{theorem}

\begin{remark}
\label{Rpostmeanv}
The result \eqref{DGbvtGttGd3r0} of Theorem~\ref{Tpostmeanv} does not apply to \( \QP = \DPGPt \)
for \( \dimp = \infty \) because \( \DPGPt \bigl( \vthetav_{\GP} - \tilde{\thetav}_{\GP} \bigr) \) 
is nearly standard normal given \( \Yv \) and \( \dimt_{\QPGP} = \infty \).
However, it well applies to \( \QP = \DPt \) yielding 
\begin{EQA}
	\bigl\| \DPt \bigl( \bar{\vthetav}_{\GP} - \tilde{\thetav}_{\GP} \bigr) \bigr\| 
	& \lesssim &
	\delta_{3}(\rups) \sqrt{\dimt_{\GP}}  \, ,
 	\\
	\E \Bigl(
        \bigl\|  
          \DPt \bigl( \vthetav_{\GP} - \tilde{\thetav}_{\GP} \bigr) 
        \bigr \|^{2} 
        \cond \Yv 
      \Bigr) 
     & \approx &
     \dimt_{\GP} .
\label{DtbvttGd3r0spG}
\end{EQA}
Another typical choice of \( \QP \) is \( \QP = \Pi_{m} \DPGPt \), where \( \Pi_{m} \) is 
the projector on the first \( m \) eigenvectors of \( \DPGPt \).
Then \( \dimt_{\QPGP} = m \).
\end{remark}


Now we state the result about the Gaussian approximation of 
\( \vthetav_{\GP} - \bar{\vthetav}_{\GP} \cond \Yv \).

\begin{theorem}
\label{TbarthetavG}
Let the conditions of Theorem~\ref{Tpostmeanv} be satisfied, 
let \( \QP = \QP \Pi \) where \( \Pi \) is a projector in \( \R^{\dimp} \).
Then it holds on the set \( \Omega(\xx) \) shown in Proposition~\ref{TconstGrDG}
\begin{EQA}
    \sup_{\rr > 0}
    \left| 
    \P\Bigl( 
    	\bigl\| \QP \bigl( \vthetav_{\GP} - \bar{\vthetav}_{\GP} \bigr) \bigr\| \leq \rr \cond \Yv 
    \Bigr)
    - \PG\Bigl( \| \QP \DPGPt^{-1} \gammav \| \leq \rr \Bigr) 
    \right|
    & \lesssim &
    \err(\rups) \sqrt{\dimt_{\Pi}} 
    + n^{-1} \dimt_{\Pi} 
\label{PvtGttGAcYb}
\end{EQA}
with \( \dimt_{\Pi} \eqdef \tr \bigl( \Pi \DPGPt^{2} \Pi \DPGPt^{-2} \Pi \bigr) \).
\end{theorem}

\begin{remark}
If \( \DPGPt \Pi = \Pi \DPGPt \), then \( \dimt_{\Pi} = \dim(\Pi) \).
In particular, if \( \QP = \DPGPt \Pi_{m} \) for the eigenprojector \( \Pi_{m} \) as in Remark~\ref{Rpostmeanv},
then \( \dimt_{\GP} = \dimt_{\Pi_{m}} = m \).
One can see that the use of the posterior mean instead of the penalized MLE \( \tilde{\thetav}_{\GP} \)
is justified under a bit stronger condition ``\( \err(\rups) \sqrt{\dimt_{\Pi}} + n^{-1} \dimt_{\Pi} \) is small''
compared to ``\( \err(\rups) \) is small''.
\end{remark}


\subsubsection{Prior comparison and prior impact}
\label{Spriorimpact}
The classical BvM result claims that the prior impact asymptotically washes out,
as the sample size increases.
The posterior becomes close to the normal distribution with the same distribution as for the MLE \( \tilde{\thetav} \),
namely, to \( \ND(\tilde{\thetav},\DPc^{-2}) \).
It is well understood that a general BvM result is impossible in a infinite dimensional nonparametric set-up 
whatever sample size is. 
In this section we want to quantify the accuracy of the BvM approximation using the obtained bounds 
on the Gaussian approximation of the posterior and the results on Gaussian comparison. 
As in previous sections we show that restricting to the class of elliptic sets helps to improve the bounds.
We slightly change the statement of the problem and consider it as a problem of prior impact.
Let \( \GP^{-2} \) and \( \GPm^{-2} \) be prior covariance matrices for two different priors.
We want to compare their posteriors. 
A special case of interest is when \( \GPm \) is large and corresponds to a model of low complexity,
which can be measured by the effective dimension.
The second prior has a smaller precision matrix \( \GP^{2} \) and hence a more complex model.
By Theorem~\ref{TnonparBvMm} the posterior \( \vthetav_{\GP} \cond \Yv \) is nearly Gaussian
\( \ND(\tilde{\thetav}_{\GP},\DPGPt^{-2}) \).
Similarly the posterior \( \vthetav_{\GPm} \cond \Yv \) is close to \( \ND(\tilde{\thetav}_{\GPm},\DPt_{\GPm}^{-2}) \).
We are interested if elliptic credible sets calibrated for the simple \( \GPm \)-prior
can be used for the more complex \( \GP \)-prior.
The Gaussian approximation reduces this question to Gaussian comparison; see \cite{GNSUl2017}
or Theorem~\ref{Tgaussiancomparison3}.
Motivated by the above discussion we assume \( \GP^{2} \leq \GPm^{2} \).

\begin{theorem}
\label{TnonparamBvM}
Let the conditions of Theorem~\ref{Tpostmeanv} be satisfied for two priors \( \ND(0,\GP^{-2}) \) and 
\( \ND(0,\GPm^{-2}) \) with \( \GP^{2} \leq \GPm^{2} \). 
Then it holds on a set \( \Omega(\xx) \) with \( \P\bigl( \Omega(\xx) \bigr) \geq 1 - 2/n \)        
\begin{EQA}
    && \nquad  
    \sup_{\rr}
    \left| 
		\P\Bigl( 
	  		\bigl\| \QP \bigl( \vthetav_{\GP} - \tilde{\thetav}_{\GPm} \bigr) \bigr\| \leq \rr \cond \Yv 
		\Bigr)    
		- \P\Bigl( 
	  		\bigl\| \QP \bigl( \vthetav_{\GPm} - \tilde{\thetav}_{\GPm} \bigr) \bigr\| \leq \rr \cond \Yv 
		\Bigr)    
	\right|
	\\
	& \lesssim &
	\delta_{3}(\rups) + n^{-1} + \frac{1}{\| \QP \DPGPt^{-2} \QP^{\T} \|_{\Fr}} 
    \Bigl\{ \tr\bigl(\QP (\DPGPt^{-2} - \DPt_{\GPm}^{-2}) \QP^{\T}\bigr) 
    + \bigl\| \QP \bigl( \tilde{\thetav}_{\GP} - \tilde{\thetav}_{\GPm} \bigr) \bigr\|^{2} \Bigr\}	.
    \qquad
\label{PvtGttGAcY2}
\end{EQA}
\end{theorem}

\begin{remark}
The last term in the bound \eqref{PvtGttGAcY2} comes from the Gaussian comparison result
of Theorem~\ref{Tgaussiancomparison3}.
It includes the ``variance'' part that relates two covariance operators 
\( \QP \DPGPt^{-2} \QP^{\T} \) and
\( \QP \DPt_{\GPm}^{-2} \QP^{\T} \), and the ``squared bias'' term 
\( \| \QP \bigl( \tilde{\thetav}_{\GP} - \tilde{\thetav}_{\GPm} \bigr) \|^{2} \).
Applicability of the prior \( \GPm \) in place of \( \GP \) is justified under ``small bias'' condition 
\( \| \QP \bigl( \tilde{\thetav}_{\GP} - \tilde{\thetav}_{\GPm} \bigr) \|^{2} \ll \| \QP \DPGPt^{-2} \QP^{\T} \|_{\Fr} \), and under the ``variance'' condition
\( \dimt_{\GP} - \dimt_{\GPm} \ll \| \QP \DPGPt^{-2} \QP^{\T} \|_{\Fr} \).
\end{remark}

\subsubsection{A family of priors and uniform Gaussian approximation}
\label{Spriorfamily}
Consider a more general situation when a family of Gaussian priors \( \ND(0,\GP_{\hyper}^{-2}) \),
\( \hyper \in \MM \), is given.
For each of them, under appropriate conditions, one can state the Gaussian approximation result as in 
Corollary~\ref{CnonparBvMm} or Theorem~\ref{TnonparBvMm}.
The choice of a prior by empirical or full Bayes approaches requires to state this approximation
uniformly in \( \hyper \in \MM \).
Surprisingly, in the contrary to the classical frequentist model selection, 
such a uniform approximation can be stated in a straightforward way. 
In fact, all the results about posterior distribution are stated conditionally on the data after
restricting to the random set \( \Omega(\xx) \) on which 
a deviation bound \( \| \DPGP^{-1} \nabla \zetav \| \leq \zq(\BB_{\VPGP},\xx) \) holds.
All we need is a uniformly in \( \hyper \) version of this bound.
Note that the conditions \nameref{LLref}, \nameref{Eref}, \nameref{ED0ref},
and \nameref{LL0ref} do not involve any prior.
The only condition \nameref{HGref} 
has to be fulfilled uniformly in \( \GP \in \bigl\{ \GP_{\hyper} \bigr\} \). 
In the case when the family \( \bigl\{ \GP_{\hyper} \bigr\} \) contains the smallest covariance 
\( \GP_{\min}^{2} \leq \GP_{\hyper}^{2} \), it suffices to check 
\nameref{HGref} for \( \GP_{\min}^{2} \) only.
This is trivially fulfilled for the family of smooth \( (\smp,\CGP) \)-priors for different 
\( \smp > 1 \).  
For \( \mm \)-truncation priors, one can take any \( \mm \leq \CONST \nsize^{1/(2\smp_{\min}+1)} \).

The next result describes the uniform properties of the estimators 
\( \tilde{\thetav}_{\hyper} \) and the posteriors \( \vthetav_{\hyper} \cond \Yv \).
Everywhere we write the subindex \( \hyper \) in place of \( \GP_{\hyper} \).
In particular, 
\( \DPc_{\hyper}^{2} = \IF_{\hyper}(\thetavs_{\hyper}) 
= \IF(\thetavs_{\hyper}) + \GP_{\hyper}^{2} \) with
\( \IF(\thetav) = - \nabla^{2} \E L(\thetav) \).

\begin{theorem}
\label{TGARunif}
Let the conditions \nameref{LLref}, \nameref{Eref}, \nameref{ED0ref}, and \nameref{LL0ref} be fulfilled. 
Let also \( \bigl\{ \ND(0,\GP_{\hyper}^{-2}) \, , \hyper \in \MM \bigr\} \) be a family of Gaussian priors, and let 
\nameref{HGref} be satisfied uniformly in \( \hyper \in \MM \).
Let also there exist a random set \( \Omega(\xx) \) with \( \P\bigl( \Omega(\xx) \bigr) \geq 1 - \ex^{-\xx} \)
and an upper function \( \zq(\BB_{\VPGP_{\hyper}},\xx) \) of \( \hyper \)
such that it holds on \( \Omega(\xx) \)
\begin{EQA}
	\| \DP_{\hyper}^{-1} \, \nabla \zetav \| 
	& \leq &
	\zq(\BB_{\VPGP_{\hyper}},\xx) ,
	\qquad
	\forall \hyper \in \MM.
\label{DGkm1nzzVGkx}
\end{EQA}
Then all the statements of Proposition~\ref{TconstGrDG} through \ref{TnonparamBvM} are fulfilled on 
\( \Omega(\xx) \) uniformly in \( \hyper \in \MM \).
In particular, each of \( \tilde{\thetav}_{\hyper} \) concentrates on the elliptic set
\begin{EQA}
	\CA_{\hyper}(\rr_{\hyper})
	& = &
	\bigl\{ \thetav \colon \| \DPc_{\hyper} (\thetav - \thetavs_{\hyper}) \| \leq \rr_{\hyper} \bigr\} ,
\label{CAhyrhytcDcht}
\end{EQA}
with the center at \( \thetavs_{\hyper} \),
while the posterior \( \vthetav_{\hyper} \cond \Yv \) concentrates on a larger vicinity 
\begin{EQA}
	\CAt_{\hyper}(\rr_{\hyper}) 
	&=& 
	\bigl\{ \thetav \colon 
	  	\bigl\| \HPc (\thetav - \tilde{\thetav}_{\hyper}) \bigr\| 
		\leq \rr_{\hyper} 
	\bigr\}
\label{CAthyrhI12ttthyr}
\end{EQA}
centered at \( \tilde{\thetav}_{\hyper} \).
Each posterior \( \vthetav_{\hyper} \cond \Yv \) can be approximated 
by the Gaussian \( \ND(\tilde{\thetav}_{\hyper},\DPt_{\hyper}^{-2}) \), the error of approximation is given by Theorem~\ref{TnonparBvMm} or Corollary~\ref{CnonparBvMm}.
\end{theorem}

This result follows from the fact that after restricting to the set \( \Omega(\xx) \) we only operate 
with deterministic function \( \E L(\thetav) \) and use its local smoothness properties from 
\nameref{LL0ref}.

\begin{remark}
The probabilistic bound 
\( \P\bigl( \| \DP_{\hyper}^{-1} \nabla \zetav \| > \zq(\BB_{\VPGP_{\hyper}},\xx) \bigr) \leq \ex^{-\xx} \) follows from \nameref{ED0ref} for each \( \hyper \in \MM \).
If \( \MM \) is a finite set and \( |\MM| \) is its cardinality then a uniform version of this bound can be 
easily obtained by an increase of \( \xx \) to \( \xx + \log|\MM| \).
\end{remark}

\subsubsection{Non-Gaussian priors} 
The assumption of a Gaussian prior was essential in the stated results. 
However, the applied technique allows to extend the approach to a broad class of non-Gaussian priors under some mild assumptions.
For simplicity we restrict ourselves to the case of Gaussian likelihood \( L(\thetav) \) which is a quadratic function of \( \thetav \), 
so that \( -\nabla^{2} L(\thetav) = \IF \) for a symmetric operator \( \IF \) in \( \R^{\dimp} \). 
Now, let the prior has a log-concave density \( \Prior(\thetav) \) which we also assume to be a sufficiently smooth.
More precisely, we suppose that  
\begin{EQA}
	\GP^{2}(\thetav)
	& \eqdef &
	- \nabla^{2} \log \Prior(\thetav)
	\geq 
	0.
\label{G2tmn2lP}
\end{EQA}
Define 
\begin{EQA}
	\IF_{\GP}(\thetav)
	& \eqdef &
	\IF + \GP^{2}(\thetav)
\label{DG2tD2G2t}
\end{EQA}
and redefine \( \delta_{m}(\thetav,\uv) \) in~\eqref{d3utd4utd4dt4c} using \( \log \Prior(\thetav) \) in place of \( \E L(\thetav) \).
It is rather straightforward to see that with this exchange, all the previous results continue to apply without any change.  


\section{Examples}
\label{SexamplesBvM}
In this section we illustrate the general results of the Section~\ref{SnonBvM} by applying to nonparametric density estimation and generalized regression.
Log-density model is a popular example in statistical literature related to BvM Theorem and nonparametric Bayes study. 
We mention~\cite{CaNi2014}, \cite{CaRo2015} among many others.
Generalized regression model includes the logit model for binary response or classification problems, Poisson and Cox regression, several reliability models and so on. 
The related BvM results can be found e.g. in~\cite{CaNi2014}, \cite{GhVa2017} and references therein.
The results on Gaussian approximation of the posterior are typically asymptotic and do not provide any accuracy guarantees for this approximation.
Our results are stated for finite samples and deliver the quantitative and tight bounds on the accuracy of this approximation in terms of effective dimension of the problem.

\subsection{Nonparametric log-density estimation}
\label{SGBvM}
Suppose we are given a random sample \( X_{1},\ldots,X_{\nsize} \) in \( \R^{d} \).
The i.i.d. model assumption means that all these random variables are independent 
identically distributed from some measure \( \Pone \) with a density 
\( \dens(x) \) with respect to a \( \sigma \)-finite measure \( \Pdom \) in \( \R^{d} \).
This density function is the target of estimation.
By definition, the function \( \dens \) is non-negative, measurable, and integrates to one:
\( \int \dens(x) \, \Pdom(dx) = 1 \).
%
Here and below, the integral \( \int \) without limits means the integral over the whole space \( \R^{d} \).
If \( \dens(\cdot) \) has a smaller support \( \XX \), one can restrict integration to this set.
%
Below we parametrize the model by a linear decomposition of the log-density function.
Let \( \bigl\{ \psi_{j}(x), \, j=1,\ldots,\dimp \bigr\} \) with \( \dimp \leq \infty \)
be a collection of functions in \( \R^{d} \) (a dictionary).
%
For each \( \thetav = (\theta_{j}) \in \R^{\dimp} \), define   
\begin{EQA}
    \ldens(x,\thetav)
    &\eqdef &
    \sum_{j=1}^{\dimp} \theta_{j} \psi_{j}(x) - \cdens(\thetav)
    =
    \bigl\langle \Psi(x), \thetav \bigr\rangle - \cdens(\thetav),
\label{logdenssumj}
\end{EQA}
where \( \cdens(\thetav) \) is given by 
\begin{EQA}
    \cdens(\thetav)
    & \eqdef &
    \log \int \exp\biggl\{ \sum_{j=1}^{\dimp} \theta_{j} \psi_{j}(x) \biggr\} \, \Pdom(dx)
    =
    \log \int \ex^{\langle \Psi(x), \thetav \rangle} \, \Pdom(dx) .
\label{intRdsumjtpj}
\label{gtdelointTPxm0}
\end{EQA}
Here \( \Psi(x) \) is a vector with components \( \psi_{j}(x) \).
%
Linear log-density modeling assumes  
\begin{EQA}
  	\log \dens(x)
  	&=&
  	\ldens(x,\thetavs)
	=
	\bigl\langle \Psi(x), \thetavs \bigr\rangle - \cdens(\thetavs)
\label{dnesdxtsl}
\end{EQA}  
for some \( \thetavs \in \Theta \subseteq \R^{\dimp} \).
A nice feature of such representation is that the function \( \log \dens(x) \) in the contrary to the density itself does not need to be non-negative.
One more important benefit of using the log-density is that the stochastic part of the corresponding log-likelihood is \emph{linear} w.r.t. the parameter \( \thetav \).
The log-likelihood \( L(\thetav) \) reads as
\begin{EQA}[c]
    L(\thetav) 
    = 
    \sum \ell(X_{i},\thetav) 
    = 
    \sum \bigl\langle \Psi(X_{i}), \thetav \bigr\rangle - \nsize \cdens(\thetav) 
    = 
    \langle \Spsi, \thetav \rangle - \nsize \cdens(\thetav) ,
    \quad
    \Spsi
    = \sum \Psi(X_{i}),
\label{density_likelyhood}
\end{EQA}
where  \( \sum \) means \( \sum_{i=1}^{n} \).
%
For applying the general results of Section~\ref{SnonBvM}, it suffices to check
the general conditions of Section~\ref{Scondgllo} for the log-density model.
First note that the generalized linear structure of the model yields automatically 
conditions \nameref{LLref} and \nameref{Eref}.
Indeed, convexity of \( \cdens(\cdot) \) implies that 
\( \E L(\thetav) = \langle \E \Spsi, \thetav \rangle - \nsize \cdens(\thetav) \) is concave. 
Further, for the stochastic component \( \zeta(\thetav) = L(\thetav) - \E L(\thetav) \), it holds 
\begin{EQA}
\label{stochDensDef}
    \nabla \zeta(\thetav)
    &=&
    \nabla \zeta
    =
    \Spsi - \E \Spsi
    = 
    \sum \bigl[\Psi(X_{i}) - \E \, \Psi(X_{i})\bigr],
\label{ztLtEltnztSPXi}
\end{EQA}
and \nameref{Eref} follows.
Further,
the representation \( \E L(\thetav) = \langle \E \Spsi, \thetav \rangle - \nsize \cdens(\thetav) \) implies
\begin{EQA}
	\IF(\thetav)
	&=&
	- \nabla^{2} \E L(\thetav)
	=
	n \nabla^{2} \cdens(\thetav) .
\label{IFttn2ELtnn2cdt}
\end{EQA}
To simplify our presentation, we assume that \( X_{1},\ldots,X_{\nsize} \) are indeed i.i.d.
This can be easily extended to non i.i.d. r.v.'s at cost of more complicated notations.
Then
\begin{EQA}
    \E \Spsi
    &=&
    \sum \E \, \Psi(X_{i}) = \nsize \, \E\Psi(X_{1}) 
    = \nsize \, \Psimean
\label{EPsi1nEPXi}
\end{EQA}
with \( \Psimean = \E \, \Psi(X_{1}) \).
%
We further assume that the underlying density \( \dens(x) \)
can be represented in the form \eqref{dnesdxtsl} for some parameter vector
\( \thetavs \). 
It also holds
\begin{EQA}
    \thetavs
    & = &
    \argmax_{\thetav \in \Theta} \E L(\thetav)
    =
    \argmax_{\thetav \in \Theta} 
      	\bigl\{ \langle \E \Spsi, \thetav \rangle - \nsize \cdens(\thetav) \bigr\} 
	=
    \argmax_{\thetav \in \Theta} 
      	\bigl\{ \langle \Psimean, \thetav \rangle - \cdens(\thetav) \bigr\} .
\label{ttsT1dedef}
\end{EQA}
Now we switch to the Bayesian framework and restrict ourselves to \( \mm \)-truncation or 
\( (\smp,\CGP) \)-smoothing priors;
see Section~\ref{Seffdim}.
In the second case, we still truncate the prior at the fixed level 
\( \mms = \nsize^{1/3}/\log \nsize \).
This means that the prior is supported to the subspace \( \R^{\mms} \) of \( \R^{\dimp} \).
For a given penalty operator \( \GP^{2} \), the corresponding penalized MLE 
\( \tilde{\thetav}_{\GP} \), and the target \( \thetavs_{\GP} \) are
\begin{EQ}[rcl]
    \tilde{\thetav}_{\GP}
    &=&
    \argmax_{\thetav \in \Theta} L_{\GP}(\thetav)
    =
    \argmax_{\thetav \in \Theta} 
      \Bigl\{ \langle \thetav, \Spsi \rangle - \nsize \cdens(\thetav) 
    - \frac{1}{2} \| \GP \thetav \|^{2} \Bigr\} ,
    \\
    \thetavs_{\GP}
    &=&
    \argmax_{\thetav \in \Theta} \E L_{\GP}(\thetav)
    =
    \argmax_{\thetav \in \Theta} \Bigl\{ 
    	\langle \thetav, \E \Spsi \rangle - \nsize \cdens(\thetav) 
    	- \frac{1}{2} \| \GP \thetav \|^{2} 
    \Bigr\} .
\label{tsatTqLGttT122}
\end{EQ}
We write \( \DPGP^{2} = \IF_{\GP}(\thetavs_{\GP}) \) and \( \dimA_{\GP} = \dimA_{\GP}(\thetavs_{\GP}) \).
Below we assume:
\begin{description}
\item[\label{Thetaref} \( \bb{(\Theta)} \)]
\( \Theta = \BBB_{\smp}(1) \) for \( \smp > 1 \) yielding 
\( \sum_{j} \theta_{j}^{2} j^{2\smp} \leq 1 \) for all \( \thetav = (\theta_{j}) \).

\end{description}
This condition is standard in of log-density estimation; cf. \cite{CaNi2014}.

\begin{description}
\item[\label{Psithetaref} \( \bb{(\psi_{j}) }\)]
Define \( q_{j}^{-2} = j \log^{2} (j) \). Then 
\begin{EQA}
	\sup_{x \in \XX}  
		\sum_{j \geq 1} \psi_{j}^{2}(x) \, q_{j}^{2} 
	& \leq &
	\CONSTpsi^{2}
	 \, .
\label{ggsjpsjxhjuj}
\end{EQA}

\end{description}

One can check that this condition is fulfilled in two important special case:
(1) all the basis functions \( \psi_{j}(x) \) are uniformly bounded by a constant \( \CONSTpsi \),
e.g. Fourier or cosine basis;
(2) \( (\psi_{j}) = \bigl( \psi_{kl} \bigr)_{k \geq 0,l \in \II_{k}} \) is a double indexed set of wavelet functions 
satisfying 
\begin{EQA}
	\sup_{x \in \XX} \sum_{l \in \II_{k}} \psi_{kl}^{2}(x)
	& \leq &
	\CONST 2^{k} .
\label{CiksxXXslIIk}
\end{EQA}

In what follows, for some small but fixed \( \rsmall \), we use \( \Thetad \) of the form
\begin{EQA}
	\Thetad
	=
	\Theta_{\rsmall}
	&=&
	\bigl\{ \thetav = \thetavs + \uv \colon 
		\bigl\langle \nabla^{2} \cdens(\thetav) \uv,\uv \bigr\rangle \leq \rsmall^{2} 
	\bigr\} \, .
\label{TTrrna22uur2}
\end{EQA}

\begin{description}
	\item[\label{phiref}\( \bb{(\nabla^{2}\phi)} \)]
For some \( \CONSTi_{\phi,1},\CONSTi_{\phi,2} \geq 1 \) and for all \( \thetav \in \Thetad \)
and all \( \uv \in \R^{\mms} \)
\begin{EQA}
	\CONSTi_{\phi,1}^{-2} 
	\leq 
	\| \uv \|^{-2} \bigl\langle \nabla^{2} \cdens(\thetav) \uv, \uv \bigr\rangle
	& \leq &
	\CONSTi_{\phi,2}^{2}  .
\label{infuUdln2cduu2}
\end{EQA}
\end{description}
This is an identifiability condition ensuring that different features can
be well identified from the data.

Define for any \( \thetav \in \Theta \) a measure \( P_{\thetav} \) by the relation:
\begin{EQA}
    \frac{d P_{\thetav}}{d\Pdom}(x)
    &=&
    \exp\bigl\{ \bigl\langle \Psi(x), \thetav \bigr\rangle - \cdens(\thetav) 
    \bigr\} .
\label{dPdlxeLtSn}
\end{EQA}
The identity \eqref{gtdelointTPxm0} ensures that \( P_{\thetav} \) is a probabilistic measure.
Moreover, due to \nameref{Thetaref} and \nameref{Psithetaref} all such measures are equivalent in the sense 
the ratio \( d\P_{\thetav}/dP_{\thetavd} \) is bounded by a constant for all \( \thetav, \thetavd \)
in \( \Theta \) because \( \sup_{x \in \XX}| \langle \thetav - \thetavd, \Psi(x) \rangle | < \infty \).

\begin{description}
\item[\label{Psiuref} \( \bb{(\Psi\uv)} \)]
  \textit{
  There exists a constant \( \CONSTfour \geq 1 \) such that 
  it holds for all \( \thetav \in \Thetad \), all \( \uv \in \R^{\mms} \), and \( k=3,4 \)} 
\begin{EQA}
    && \nquad
    E_{\thetav} \bigl| 
    	\bigl\langle \Psi(X_{1}) - E_{\thetav} \Psi(X_{1}), \uv \bigr\rangle 
    \bigr|^{k}
    \leq 
    \bigl\{ 
    	\CONSTfour^{2} E_{\thetav} \langle \Psi(X_{1}) - E_{\thetav} \Psi(X_{1}), \uv \rangle^{2} 
    \bigr\}^{k/2} .
    \qquad
\label{EtuTPsXimEtuT}
\end{EQA}

\end{description}

\begin{remark}
In fact, in view of \eqref{ggsjpsjxhjuj}, it suffices to check 
\eqref{infuUdln2cduu2} and \eqref{EtuTPsXimEtuT} 
for one point \( \thetav \in \Thetad \), in particular, for the true point \( \thetavs \)
corresponding to the underlying measure \( P \).
Condition \nameref{Psiuref} means that
each measure \( P_{\thetav} \) for \( \thetav \) in a vicinity \( \Thetad \) of \( \thetavs \) satisfies 
a kind of Khinchin's inequality which relates the fourth and the second directional moments 
of \( \Psi(X_{1}) \).
Conditions like \eqref{EtuTPsXimEtuT} are often used in high dimensional probability. 
It is automatically fulfilled for Gaussian measures and for the case when the individual features 
\( \psi_{j}(x) \) of the vector \( \Psi(x) \) are independent and have fourth moments. 

\end{remark}

%
%
Below by \( \CONST_{0}, \CONST_{1}, \CONST_{2}, \ldots \) we denote some fixed constants 
which possibly depending on \( \CONSTfour \) and \( \CONSTphi \) from \nameref{Psiuref}
and \nameref{phiref}.

\begin{theorem}
\label{Tconcdens}
Assume \nameref{Thetaref}, \nameref{Psithetaref}, 
\nameref{phiref}, and
\nameref{Psiuref}.
Let \( \GP^{2} \) correspond to a \( \mm \)-truncation prior with \( \mm \leq \mms \)
yielding \( \dimA_{\GP} \leq \mm \), 
or to a \( \mms \)-truncated \( (\smp,\CGP) \)-smooth prior with \( \CGP \gg \nsize^{1 - 2\smp/3} \)
and \( \dimA_{\GP} \lesssim (\nsize/\CGP)^{1/(2\smp)} \).
%
Define 
\begin{EQA}
	\zq_{\GP} 
	& \eqdef & 
	\sqrt{\dimA_{\GP}} \, + \sqrt{2 \log \nsize} \, . 
\label{zqGPdefpG2n}
\end{EQA}
Let finally \( \thetavs_{\GP} \in \Thetad \).
Then 
it holds for \( \nsize \geq \CONST_{0} \zq_{\GP}^{2} \)
on a set \( \Omega_{\nsize} \) with \( \P(\Omega_{\nsize}) \geq 1 - 3/\nsize \) 
\begin{EQ}[rcl]
    \| \DPGP (\tilde{\thetav}_{\GP} - \thetavs_{\GP}) \| 
    & \leq &
    2 \zq_{\GP} \, ,
    \\
    \bigl\| \DPGP \bigl( \tilde{\thetav}_{\GP} - \thetavs_{\GP} \bigr) - \DPGP^{-1} \nabla \zeta \bigr\|
    & \leq &
    \CONST_{1} \, \zq_{\GP}^{2} \, \nsize^{-1/2} ,
\label{DGttGtsde}
    \\
    \biggl| 
    \LGP(\tilde{\thetav}_{\GP}) - \LGP(\thetavs_{\GP}) 
    - \frac{1}{2} \bigl\| \DPGP^{-1} \nabla \zeta \bigr\|^{2}
    \biggr|
    & \leq &
    \CONST_{1} \, \zq_{\GP}^{3} \, \nsize^{-1/2} .
\label{3d3Af12DGtde}
\end{EQ}
\end{theorem}

\noindent

Now we continue with the properties of the posterior \( \vthetav_{\GP} \cond \Xv \) for a Gaussian prior \( \ND(0,\GP^{2}) \).
Remind 
 \( \DPGPt^{2} = \IF_{\GP}(\tilde{\thetav}_{\GP}) = \IF(\tilde{\thetav}_{\GP}) + \GP^{2} \).

\begin{theorem}
\label{Tconcdensex}
Suppose that the conditions of Theorem~\ref{Tconcdens} hold.
Then Theorem~\ref{PrhoQPBvM} through \ref{CThonestCS} 
continue to apply to the posterior \( \vthetav_{\GP} \cond \Xv \).
In particular, on a set \( \Omega_{\nsize} \) with \( \P\bigl( \Omega_{\nsize} \bigr) \geq 1 - 3/\nsize \), it holds with \( \zq_{\GP} \) from \eqref{zqGPdefpG2n}
\begin{EQA}
\label{1p1m1a1emxrrde}
	\sup_{A \in \cc{B}_{\smp}(\R^{\dimp})}
	\left| 
		\P\bigl( \vthetav_{\GP} - \tilde{\thetav}_{\GP} \in A \cond \Xv \bigr)
		- \PG\bigl( \DPGPt^{-1} \gammav \in A \bigr)  
	\right|
	& \leq &
	\CONST_{4}\, \zq_{\GP}^{6} \,  \nsize^{-1} \, ,
	\\
\label{1p1m1a1emxrrasde}
	\sup_{A \in \cc{B}(\R^{\dimp})}
	\left| 
		\P\bigl( \vthetav_{\GP} - \tilde{\thetav}_{\GP} \in A \cond \Xv \bigr)
		- \PG\bigl( \DPGPt^{-1} \gammav \in A \bigr)  
	\right|
	& \leq &
	\CONST_{4} \, \zq_{\GP}^{3} \, \nsize^{-1/2} .
\end{EQA}
\end{theorem}

%
%
Under a proper choice \( \mm = \CGP = \nsize^{1/(2\smp+1)} \) of the prior parameters
one can derive finite sample versions of 
the standard nonparametric \emph{rate optimal} results about concentration of the pMLE
and posterior contraction; cf \cite{CaNi2014}, \cite{CaRo2015}.
We use that \( \DPGPt^{2} \geq \nsize \, \CONSTphi^{2} \, \Id_{\dimp} \)
and \( \nsize^{-1} \zq_{\GP}^{2} \asymp \nsize^{-1} \mm \asymp \nsize^{-2\smp/(2\smp+1)} \).
The error term in \eqref{1p1m1a1emxrrde} is of order \( \nsize^{(2 - 2 \smp) / (2\smp + 1)} 
\to 0 \) as \( \nsize \to \infty \) because \( \smp > 1 \).

\begin{corollary}
\label{CTconcdens}
Assume conditions of Theorem~\ref{Tconcdens}.
Define \( \mm = \CGP = \nsize^{1/(2\smp+1)} \).
For the \( \mm \)-truncation prior 
or for the \( (\smp,\CGP) \) smooth prior, 
it holds on the same set \( \Omega_{\nsize} \) 
\begin{EQA}
	\| \tilde{\thetav}_{\GP} - \thetavs \|^{2}
    & \leq &
    \CONST \nsize^{-2\smp/(2\smp+1)} \, ,
    \\
    \P\Bigl( \| \vthetav_{\GP} - \thetavs \|^{2}  
    		&>& \CONST \nsize^{-2\smp/(2\smp+1)}
    	\cond \Xv 
    \Bigr)
    \leq 
    n^{-1} \, ,
\label{DGttGtsCzGns}
\end{EQA}
for \( \nsize \) sufficiently large.
Moreover, under ``undersmoothing'' \( \thetavs \in \BBB_{\smp^{*}}(1) \) with \( \smp^{*} > \smp \),
Bayesian credible sets 
\( \CAt_{\GP}(\rr_{\alp}) =
\bigl\{ \thetav \colon \bigl\| \sqrt{\nsize} (\tilde{\thetav}_{\GP} - \thetav) \bigr\| \leq \rr_{\alp} \bigr\} \) with \( \rr_{\alp} \) satisfying
\( \PG\bigl( \sqrt{\nsize} \bigl\| \DPGPt^{-1} \gammav \bigr\| > \rr_{\alp} \bigr) = \alp \),
\( \gammav \sim \ND(0,\Id_{\dimp}) \) are asymptotic valid; see \eqref{Palpo1AGra}.
\end{corollary}


\subsection{Generalized regression}
\label{SEFcposter}
Now we discuss how the general results apply to generalized regression.
Suppose we are given independent data \( Y_{1},\ldots,Y_{\nsize} \) which follow the model 
\begin{EQA}
	Y_{i}
	& \sim &
	P_{\upsilon_{i}}
	\in \PEF,
	\qquad
	i=1,\ldots,n,
\label{YiPvii1n}
\end{EQA}
where \( \PEF = (P_{\upsilon}, \upsilon \in \Upsilon) \) be a univariate exponential family with a canonical parameter.
The latter means that \( \PEF \) is dominated by a \( \sigma \)-finite measure \( \mu \) and 
\begin{EQA}
	\log \frac{d P_{\upsilon}}{d\mu}(y)
	&=&
	\upsilon y - \cdens(\upsilon) + \ell(y)
\label{ldPvdmyuyCy}
\end{EQA}
for a convex function \( \cdens(\upsilon) \) of a univariate parameter \( \upsilon \).
A typical example is given by the logistic regression with binary observations \( Y_{i} \).
Then \( \cdens(\upsilon) = \log(1 + \ex^{\upsilon}) \).
%
%
The model \eqref{YiPvii1n} yields \( \E Y_{i} = \cdens'(\upsilon_{i}) \) and \( \Var(Y_{i}) = \cdens''(\upsilon_{i}) \).
We, however, do not assume that the model is correct.
The value \( \upsilon_{i} \) is just defined by the canonical link \( \E Y_{i} = \cdens'(\upsilon_{i}) \).

Generalized regression assumes that the \( \upsilon_{i} \)'s in \eqref{YiPvii1n} are values of a  function \( \fs(X_{i}) \) at deterministic design points \( X_{1},\ldots,X_{\nsize} \).
A linear basis expansion \( \fs(x) = \sum_{j} \theta_{j} \psi_{j}(x) \) leads to a generalized linear model
\begin{EQA}
	Y_{i}
	& \sim &
	P_{\langle \Psiv_{i},\thetav \rangle}
\label{YisPtTPiv}
\end{EQA}
with \( \thetav = (\theta_{1},\ldots,\theta_{\dimp})^{\T} \in \R^{\dimp} \) and
\( \Psiv_{i} = \bigl( \psi_{1}(X_{i}),\ldots,\psi_{\dimp}(X_{i}) \bigr)^{\T} \in \cc{X} \subset \R^{\dimp} \) for \( \dimp \leq \infty \).
The model \eqref{YisPtTPiv} yields 
\begin{EQA}
	L(\thetav)
	&=&
	\sum Y_{i} \, \langle \Psiv_{i},\thetav \rangle 
		- \cdens\bigl(\langle \Psiv_{i},\thetav \rangle \bigr) .
\label{LtsYiPiTcPiT}
\end{EQA}
Here again \( \sum = \sum_{i=1}^{n} \).
The corresponding Fisher information operator reads as
\begin{EQA}
	\IF(\thetav)
	=
	- \nabla^{2} L(\thetav)
	&=&
	\sum \cdens''\bigl(\langle \Psiv_{i},\thetav \rangle \bigr) \, \Psiv_{i} \otimes \Psiv_{i}
	\geq 
	0, 
\label{n2LtsiPiPiT}
\end{EQA}
because \( \cdens \) is strictly convex. 
This yields \nameref{LLref}.
Define \( \eps_{i} = Y_{i} - \E Y_{i} \).
Then the stochastic component of the log-likelihood is linear in \( \thetav \) and~\nameref{Eref} is fulfilled with 
\begin{EQA}
	\zeta(\thetav)
	&=&
	L(\thetav) - \E L(\thetav)
	=
	\sum \eps_{i} \, \langle \Psiv_{i},\thetav \rangle,
	\qquad
	\nabla \zeta
	=
	\sum \eps_{i} \, \Psiv_{i} \, .
\label{ztLtELtsinPi}
\end{EQA}
We assume a number of regularity conditions similar to the log-density case.

\begin{description}
\item[\label{PsithetaGRref} \( \bb{(\psi_{\infty}) }\)]
It holds \( \| \psi_{j} \Ind_{\XX} \|_{\infty} \leq \CONSTpsi \, j^{1/2} \) for \( j \leq \dimp \).
\smallskip

\item[\label{ThetaGRref} \( \bb{(\Theta)} \)]
\( \Theta \subseteq \BBB_{\smp}(1) \) for \( \smp > 1 \). 
This condition and \nameref{PsithetaGRref} yield by the Cauchy-Schwarz inequality, 
for any \( x \in \XX \) and \( \thetav \in \Theta \)
\begin{EQA}
	\bigl| \langle \Psiv(x),\thetav \rangle \bigr|
	=
	\biggl( \sum_{j=1}^{\dimp} \theta_{j} \psi_{j}(x) \biggr)^{2}
	& \leq &
	\biggl( \sum_{j=1}^{\dimp} \theta_{j}^{2} j^{2\smp} \biggr)
	\biggl( \sum_{j=1}^{\dimp} \CONST_{\psi}^{2} \, j^{-2\smp+1} \biggr)
	\leq 
	\frac{\CONST_{\psi}^{2}}{2\smp - 2}
	\eqdef
	\CONSTPsi \, .
\label{tTsjpjxj2sj2}
\end{EQA}
 
\item[\label{phi4ref} \( \bb{(\phi^{(k)})} \)]
\( \cdens(\cdot) \) is a smooth function with a continues fourth
derivative on the interval \( [-\CONSTPsi,\CONSTPsi] \) and the second derivative
\( \cdens''(\cdot) \) fulfills for some \( \CONST_{1,\cdens} \leq 1 \leq \CONST_{2,\cdens} \)
\begin{EQA}
	\CONST_{1,\phi} \, \cdens''(0)
	& \leq &
	\cdens''(t)
	\leq 
	\CONST_{2,\cdens} \, \cdens''(0), 
	\quad
	|t| \leq \CONSTPsi \, .
\label{CphiCh2t1t1Cp}
\end{EQA} 

\item[\label{epsiref}\( \bb{(\eps_{i})} \)]
There are \( \varrho > 0 \), \( \sigma_{\max} \), and \( \nunu \geq 1 \) such that \( \eps_{i} = Y_{i} - \E Y_{i} \) satisfy with \( \sigma_{i}^{2} = \E \eps_{i}^{2} \)
\begin{EQA}
	\max_{i \leq n} \sigma_{i}^{2}
	\leq 
	\sigma_{\max}^{2},
	& \quad &
	\sup_{|\lambda| \leq \varrho} \max_{i \leq n} \log \E \exp\bigl( \lambda \eps_{i} \bigr)
	\leq 
	\frac{\nunu^{2} \lambda^{2} \sigma_{\max}^{2}}{2}. 
\label{maisi2lCs2}
\end{EQA}

\item[\label{PsiuGRref}\( \bb{(\Psi\uv)}\)]
It holds for all \( \uv \in \R^{\mms} \) 
and for a fixed constant \( \CONSTfour \)
\begin{EQA}
    \frac{1}{\nsize} \sumi \langle \Psiv_{i},\uv \rangle^{4}
    & \leq &
    \biggl\{ \frac{\CONSTfour^{2}}{\nsize} \sumi \langle \Psiv_{i} ,\uv \rangle^{2} 
    \biggr\}^{2} \, .
\label{tuTPsXimtuT}
\end{EQA}

\item[\label{IFtref}\( \bb{(\IF)} \)]
For some \( \CONSTIF > 0 \)
\begin{EQA}
	\inf_{\uv \in \R^{\mms}} 
		\frac{\bigl\langle \IF(\thetavs) \uv,\uv \bigr\rangle}{\nsize \| \uv \|^{2}}
	& \geq &
	\CONSTIF^{-2} \, .
\label{infuUdln2cduu2GLM}
\end{EQA}

\end{description}

Further, define 
\begin{EQA}
	\VPc^{2}
	& \eqdef &
	\Var\bigl( \nabla \zeta \bigr)
	=
	\sum \sigma_{i}^{2} \, \Psiv_{i} \otimes \Psiv_{i} \, ,
	\\
	\HPc^{2}
	&=&
	\sum \sigma_{\max}^{2} \, \Psiv_{i} \otimes \Psiv_{i} \, ,
\label{VP2Cnzi2PP}
\end{EQA}
Under the correct model specification \( Y_{i} \sim P_{\langle \Psiv_{i},\thetavs \rangle} \), it holds 
\( \sigma_{i}^{2} = \cdens''(\langle \Psiv_{i},\thetavs \rangle) \) and 
\( \VPc^{2} = \IF(\thetavs) \).
Let us fix some Gaussian prior \( \ND(0,\GP^{-2}) \).
As previously, we focus on a \( \mm \)-truncation or \( \mms \)-truncated \( (\smp,\CGP) \)-smooth prior with \( \mms = \nsize^{1/3}/\log \nsize \). 
In the first case, \( \mm \leq \mms \).
For any \( \thetav \in \Theta \), conditions \nameref{Thetaref} 
and \nameref{IFtref} imply
\begin{EQA}
	\dimA_{\GP}(\thetav)
	& = &
	\tr\bigl\{ \bigl( \IF(\thetav) + \GP^{2} \bigr)^{-1} \HPc^{2} \bigr\}
	\lesssim
	\mm ,
\label{mltrItG2nm1Itp}
\end{EQA}
where \( \mm = (\nsize/\CGP)^{1/(2\smp)} \) for the \( (\smp,\CGP) \)-prior.

\begin{theorem}
\label{TGLMBvM}
Suppose \nameref{PsithetaGRref},
\nameref{ThetaGRref}, \nameref{phi4ref}, \nameref{epsiref}, \nameref{PsiuGRref}, and \nameref{IFtref} for the model \eqref{YisPtTPiv}.
Then all the expansions in \eqref{3d3Af12DGtde} for the penalized MLE 
\( \tilde{\thetav}_{\GP} \) and the properties of the posterior \( \vthetav_{\GP} \cond \Yv \)
listed in Theorem~\ref{Tconcdensex} continue to hold provided that \( \mm^{3} \ll \nsize \).
Under the prior choice with \( \mm = \CGP = \nsize^{1/(2\smp + 1)} \),
the results of Corollary~\ref{CTconcdens} hold as well.
\end{theorem}

\appendix
\section{Tools}

\subsection{Gaussian comparison}
\label{SGausscomp}
Let \(\HM\) be a Hilbert space and 
\(\Sigma_{\xiv}\) be a covariance operator  of an arbitrary Gaussian random element in \(\HM\). 
By \(\{\lambda_{k\xiv}\}_{k \geq 1}\) we denote the set of its eigenvalues arranged in the non-increasing order, i.e. 
\(\lambda_{1\xiv} \geq \lambda_{2\xiv} \geq \ldots \), and let 
\(\lambdav_{\xiv} \eqdef \diag(\lambda_{j \xiv})_{j=1}^{\infty}\). 
Note that  \(\sum_{j=1}^{\infty} \lambda_{j \xiv} < \infty\).  
Introduce the following quantities
\begin{EQA}[c]
    \Frobg_{k\xiv}^{2} 
    \eqdef 
    \sum_{j=k}^{\infty} \lambda_{j \xiv}^{2}, \quad k = 1,2,
\label{Lambda def}
\end{EQA}

\begin{theorem}[\cite{GNSUl2017}]
\label{Tgaussiancomparison3}
Let \(\xiv\) and \(\etav\) be Gaussian elements in \(\HM\) with zero mean and covariance operators 
\(\Sigma_{\xiv}\) and \(\Sigma_{\etav}\) respectively.
Then for any \(\av \in \HM\)
\begin{EQA}[rcl]
    && \nquad
    \sup_{x > 0} \left|\P( \| \xiv - \av \| \leq x) - \P( \| \etav  \| \leq x) \right| 
    \\ 
    &&
    \lesssim  
    \bigg(\frac{1}{(\Frobg_{1\xiv}\Frobg_{2\xiv})^{1/2}} + \frac{1}{(\Frobg_{1\etav}\Frobg_{2\etav})^{1/2}}\bigg) 
    \bigg( \| \lambdav_{\xiv} - \lambdav_{\etav}  \|_{1} + \| \av\|^{2}\bigg).
\label{expl_gauss22}
\end{EQA}
Moreover, assume that
\begin{EQA}[c]
    3 \| \Sigma_{\xiv}\|^{2} \le \| \Sigma_{\xiv}\|_{\Fr}^{2} 
    \quad \text{ and }  \quad 
    3 \| \Sigma_{\etav}\|^{2} \le \| \Sigma_{\etav}\|_{\Fr}^{2} \, .
\label{2plusdelta}
\end{EQA}
Then for any \(\av \in \HM\) 
\begin{EQA}
    \sup_{x > 0} \left|\P(\| \xiv - \av \| \leq x) - \P(\| \etav  \| \leq x) \right| 
    & \lesssim &
    \biggl(
    \frac{1}{\| \Sigma_{\xiv}\|_{\Fr}} 
    + \frac{1}{\| \Sigma_{\etav}\|_{\Fr}}
    \biggr) 
    \biggl( \| \lambdav_{\xiv} - \lambdav_{\etav} \|_{1} + \| \av\|^{2}\biggr).
\label{expl_gauss 2}
\end{EQA}
\end{theorem}

\def\Eg{\E_{\gammav \sim \ND(0,\Id)}}
\def\Pg{\P_{\gammav \sim \ND(0,\Id)}}
\def\HPi{\QP}
\def\dimH{\dimA}
\def\vH{\vA}
\def\BBH{W}
\def\mux{\mu_{\xx}}

\subsection{Deviation bounds for Gaussian quadratic forms}
\label{SdevboundGauss}
The next result explains the concentration effect of 
\( \langle \BB \xiv, \xiv \rangle \)
for a centered Gaussian vector \( \xiv \sim \ND(0,\HPc^{2}) \) and a symmetric trace operator \( \BB \) in \( \R^{\dimp} \),
\( \dimp \leq \infty \).
We use a version from \cite{laurentmassart2000}.
For completeness, we present a simple proof of the bound.

\begin{theorem}
\label{TexpbLGA}
\label{Lxiv2LD}
\label{Cuvepsuv0}
Let \( \xiv \sim \ND(0,\HPc^{2}) \) be a Gaussian element in \( \R^{\dimp} \) and \( \BB \) be symmetric non-negative such that \( \BBH = \HPc^{-1} \BB \, \HPc \) is a trace operator in \( \R^{\dimp} \).
Then with \( \dimA = \tr(\BBH) \), \( \vH^{2} = \tr(\BBH^{2}) \), and 
\( \supA = \| \BBH \| \), it holds for each \( \xx \geq 0 \)
\begin{EQA}
\label{Pxiv2dimAvp12}
	&&
	\P\Bigl( \langle \BB \xiv, \xiv \rangle > \zq^{2}(\BBH,\xx) \Bigr)
	\leq
	\ex^{-\xx} ,
	\\
	&& \zq(\BBH,\xx)
	\eqdef
	\sqrt{\dimA + 2 \vH \, \xx^{1/2} + 2 \supA \xx} \,\, .
\label{zqdefGQF}
\end{EQA}
It also implies 
\begin{EQA}
	\P\bigl( \| \BB^{1/2} \xiv \| > \dimA^{1/2} + (2 \supA \xx)^{1/2} \bigr)
	& \leq &
	\ex^{-\xx} .
\label{Pxiv2dimAxx12}
\end{EQA}
%
If \( \BB \) is symmetric but non necessarily positive then
\begin{EQA}
	\P\bigl( \bigl| \langle \BB \xiv, \xiv \rangle - \dimA \bigr| > 2 \vH \, \xx^{1/2} + 2 \supA \xx \bigr)
	& \leq &
	2 \ex^{-\xx} .
\label{PxivTBBdimA2vp}
\end{EQA}
\end{theorem}

\begin{proof}
W.l.o.g. assume that \( \HPc = \Id_{\dimp} \) and \( \supA = \| \BBH \| = 1 \).
Then \( \xiv = \gammav \sim \ND(0,\Id_{\dimp}) \) is standard Gaussian.
We apply the exponential Chebyshev inequality: with \( \mu > 0 \)
\begin{EQA}
	\P\Bigl( \langle \BB \gaussv, \gaussv \rangle > \zq^{2} \Bigr)
	& \leq &
	\E \exp \Bigl( \mu \langle \BB \gaussv, \gaussv \rangle / 2 - \mu \zq^{2} / 2 \Bigr) \, .
\label{PBggiz2E2mz2}
\end{EQA}
Given \( \xx > 0 \), fix \( \mu < 1 \) by the equation
\begin{EQA}
	\frac{\mu}{1 - \mu} 
	&=&
	\frac{2 \sqrt{\xx}}{\vH} \, 
	\quad \text{ or } \quad
	\mu^{-1} 
	=
	1 + \frac{\vH}{2 \sqrt{\xx}} \, .
\label{1v2sxm12m1m}
\end{EQA}
Let \( \lambda_{j} \) be ordered eigenvalues of \( \BBH \), 
\( 1 = \lambda_{1} \geq \lambda_{2} \geq \lambda_{3} \ldots \).
It holds with \( \dimH = \tr \BBH \)
\begin{EQA}
	&& \nquad
	\log \E \Bigl\{ \frac{\mu}{2} \bigl( \langle \BBH \gaussv, \gaussv \rangle - \dimA \bigr) \Bigr\}
	=
	\log \det(\Id - \mu \BB)^{-1/2} - \frac{\mu \, \dimA}{2}
	=
	- \frac{1}{2} \sum_{j \leq \dimp} \bigl\{ \log(1 - \mu \lambda_{j}) + \mu \lambda_{j} \bigr\}
	\\
	&=&
	\frac{1}{2} \sum_{j=1}^{\dimp} \sum_{k=2}^{\infty} \frac{(\mu \lambda_{j})^{k}}{k}
	\leq 
	\sum_{j=1}^{\dimp} \frac{(\mu \lambda_{j})^{2}}{4} \sum_{k=0}^{\infty} (\mu \lambda_{j})^{k}
	\leq 
	\sum_{j=1}^{\dimp} \frac{\mu^{2} \lambda_{j}^{2}}{4 (1 - \mu)} 
	=
	\frac{\mu^{2} \vH^{2}}{4 (1 - \mu)} \, .
\label{m2v241m4mj1p}
\end{EQA}
It remains to check that the choice \( \mu \) by \eqref{1v2sxm12m1m} and 
\( \zq = \zq(\BBH,\xx) \) yields
\begin{EQA}
	\frac{\mu^{2} \vH^{2}}{4 (1 - \mu)} - \frac{\mu (\zq^{2} - \dimA)}{2}
	& = &
	\frac{\mu^{2} \vH^{2}}{4 (1 - \mu)} - \mu \bigl( \vH \sqrt{\xx} + \xx \bigr)
	=
	\mu \Bigl\{ \frac{\vH \sqrt{\xx}}{2} - \vH \sqrt{\xx} - \xx \Bigr\}
	=
	- \xx
	\qquad
	\qquad
\label{m2vA241muz2}
\end{EQA}
as required.
\end{proof}
As a special case, we present a bound for the chi-squared distribution 
corresponding to \( \BB = \HPc^{2} = \Id_{\dimp} \), \( \dimp < \infty \).
Then \( \tr (\BBH) = \dimp \), \( \tr(\BBH^{2}) = \dimp \) and \( \supA(\BBH) = 1 \).

\begin{corollary}
\label{Cchi2p}
Let \( \gaussv \) be a standard normal vector in \( \R^{\dimp} \).
Then for any \( \xx > 0 \)
\begin{EQA}[ccl]
\label{Pxi2pm2px}
	\P\bigl( \| \gaussv \|^{2} \geq \dimp + 2 \sqrt{\dimp \, \xx} + 2 \xx \bigr)
	& \leq &
	\ex^{-\xx},
	\\
	\P\bigl( \| \gaussv \| \,\,  \geq \sqrt{\dimp} + \sqrt{2 \xx} \bigr)
	& \leq &
	\ex^{-\xx} ,
\label{Pxi2pm2px12}
	\\
	\P\bigl( \| \gaussv \|^{2} \leq \dimp - 2 \sqrt{\dimp \, \xx} \bigr)
	& \leq &
	\ex^{-\xx}	.
\label{Pxi2pm2px22}
\end{EQA}
\end{corollary}

\subsection{Deviation bounds for non-Gaussian quadratic forms}
\label{Sprobabquad}
\label{SdevboundnonGauss}
This section collects some probability bounds for non-Gaussian quadratic forms
starting from the subgaussian case.
Then we extend the result to the case of exponential tails. 
%
Let \( \xiv \) be a random vector in \( \R^{\dimp} \), \( \dimp \leq \infty \)
satisfying \( \E \xiv = 0 \).
We suppose that there exists an operator \( \HPc \) in \( \R^{\dimp} \) such that
\begin{EQA}
	\log \E \exp \bigl( \langle \uv, \HPc^{-1} \xiv \rangle \bigr)
	& \leq &
	\frac{\| \uv \|^{2}}{2} \, ,
	\qquad 
	\uv \in \R^{\dimp} .
\label{devboundinf}
\end{EQA}
In the Gaussian case, one obviously takes \( \HPc^{2} = \Var(\xiv) \).
In general, \( \HPc^{2} \geq \Var(\xiv) \).
We consider a quadratic form \( \langle \BB \xiv, \xiv \rangle \), where 
\( \xiv \) satisfies \eqref{devboundinf} and \( \BB \) is a given symmetric non-negative 
operator in \( \R^{\dimp} \) such that \( \BB \leq \HPc^{-2} \) and
\( \BBH = \HPc \BB \HPc \) is a trace operator:
\begin{EQA}
	\dimH
	&=&
	\tr\bigl( \BBH \bigr) 
	< \infty.
\label{dimAtHm2Bpp}
\end{EQA}
Denote also 
\begin{EQA}[c]
    \vH^{2}
    \eqdef
    \tr(\BBH^{2}) .
\label{BBrddB}
\end{EQA}   
We show that under these conditions, the quadratic form \( \langle \BB \xiv, \xiv \rangle \)
follows the same deviation bound 
\( \P\bigl( \langle \BB \xiv, \xiv \rangle \ge \zq^{2}(\BBH,\xx) \bigr) \leq \ex^{-\xx} \) 
with \( \zq^{2}(\BBH,\xx) \) from \eqref{zqdefGQF}
as in the Gaussian case.

\begin{theorem}
\label{TdevboundGauss}
Suppose \eqref{devboundinf}. 
Let \( \dimH = \tr \BBH < \infty \) for \( \BBH = \HPc \BB \HPc \).
Then
\begin{EQA}
	\P\bigl( \langle \BB \xiv, \xiv \rangle > \zq(\BBH,\xx) \bigr)
	& \leq &
	\ex^{-\xx} .
\label{lBxiq}
\end{EQA}
\end{theorem}

\begin{proof}
For any \( \mu < 1 \), we use the identity
\begin{EQA}
	\E \exp\bigl( \mu \langle \BB \xiv, \xiv \rangle / 2 \bigr)
	&=&
	\E \Eg \exp\bigl( \mu^{1/2} \langle \HPc \BB^{1/2} \gammav, \HPc^{-1} \xiv \rangle \bigr)
\label{Egexmu12lHm1B12}
\end{EQA}
Application of Fubini's theorem and \eqref{devboundinf} yields
\begin{EQA}
	\E \exp\bigl( \mu \langle \BB \xiv, \xiv \rangle / 2 \bigr)
	& \leq &
	\exp \Bigl( \frac{\mu^{2} \tr \BBH^{2}}{4 (1 - \mu)} + \frac{\mu \tr \BBH}{2} \Bigr) .
\label{wBmu2vA241mmutrBi}
\end{EQA}
Further we proceed as in the Gaussian case.
\end{proof}

Now we turn to the main case of light exponential tails of \( \xiv \).
Namely, we suppose that \( \E \xiv = 0 \) and 
for some fixed \( \gmb > 0 \) 
\begin{EQA}[c]
    \log \E \exp\bigl( \langle \uv, \HPc^{-1} \xiv \rangle \bigr)
    \le
    \frac{\| \uv \|^{2}}{2} \, ,
    \qquad
    \uv \in \R^{\dimp}, \, \| \uv \| \le \gmb ,
\label{expgamgm}
\end{EQA}
for some self-adjoint operator \( \HPc \) in \( \R^{\dimp} \), \( \HPc \geq \Id_{\dimp} \).
In fact, it suffices to assume that 
\begin{EQA}
	\sup_{\| \uv \| \leq \gmb} \E \exp\bigl( \langle \uv, \HPc^{-1} \xiv \rangle \bigr)
	& \leq &
	\CONST .
\label{sgagEexlgHx}
\end{EQA}
Then one can use the fact that existence of the exponential moment \( \E \ex^{\lambda_{0} \xi} \) for a centered random variable \( \xi \) and some fixed \( \lambda_{0} \) implies that the moment generating function 
\( f_{\xi}(\lambda) \eqdef \log \E \ex^{\lambda \xi} \) is analytic in \( \lambda \in (0, \lambda_{0}) \) with 
\( f_{\xi}(0) = f'_{\xi}(0) = 0 \) and hence, it can be well majorated by a quadratic function in a smaller interval 
\( [0,\lambda_{1}] \) for \( \lambda_{1} < \lambda_{0} \); see~\cite{GolSpo2009}.

Remind \( \BBH = \HPc \BB \HPc \).
By normalization, one can easily reduce the study to the case  
\begin{EQA}
	\| \BBH \|
	&=&
	1 .
\label{BBeq1qgnG}
\end{EQA}
Let \( \dimH = \tr(\BBH) \), \( \vH^{2} = \tr(\BBH^{2}) \), and
\( \mu(\xx) \) be defined by \( \mu(\xx) = \bigl( 1 + \frac{\vH}{2 \sqrt{\xx}} \bigr)^{-1} \);
see \eqref{1v2sxm12m1m}. 
Obviously \( \mu(\xx) \) grows with \( \xx \).
Define the value \( \xxc \) by the equation
\begin{EQA}
	\frac{\gmb - \sqrt{\dimH \, \mu(\xxc)}}{\mu(\xxc)} 
	& = &
	\zq(\BBH,\xxc) + 1 .
\label{gmsqpAHmxzBx1}
\end{EQA}
The left hand-side here decreases with \( \xx \), while the right hand-side is increasing in \( \xx \).
Therefore, the solution exists and unique.
Also denote \( \muc = \mu(\xxc) \) and
\begin{EQA}
	\gmc
	&=&
	\gmb - \sqrt{\dimH \muc} \, .
\label{gcgbsqpHmc}
\end{EQA}
Obviously \( \gmb - \sqrt{\dimH} \leq \gmc \leq \gmb \).
Moreover, as \( \zq(\BBH,\xx) \leq \sqrt{\dimH} + \sqrt{2 \xx} \), we derive that 
\begin{EQA}
	\sqrt{2\xxc}
	& \geq &
	\gmc/\muc - \sqrt{\dimH} \, .
\label{sxcgbspHspA}
\end{EQA}
Our results implicitly assume that \( \xxc \) is large, that is, 
\( \gmb \gg 2 \sqrt{\dimH} \). 

\begin{theorem}
\label{LLbrevelocroB}   
Let \eqref{expgamgm} hold and let \( \BB \) be such that \( \BBH = \HPc \BB \HPc \)
satisfies \( \| \BBH \| = 1 \)
and \( \dimH = \tr(\BBH) < \infty \).
Define \( \xxc \) by \eqref{gmsqpAHmxzBx1} and \( \gmc \) by \eqref{gcgbsqpHmc}, and suppose
\( \gmc \geq 1 \).
Then for any \( \xx > 0 \)
\begin{EQA}
    \P\bigl( \langle \BB \xiv, \xiv \rangle \ge \zqc^{2}(\BBH,\xx) \bigr)
    & \le &
    2 \ex^{-\xx} + \ex^{-\xxc} \Ind(\xx < \xxc) ,
\label{PxivbzzBBroB}
\end{EQA}    
where \( \zqc(\BBH,\xx) \) is defined by
\begin{EQA}
\label{PzzxxpBroB}
    \zqc(\BBH,\xx)
    & \eqdef &
    \begin{cases}
      \sqrt{ \dimH + 2 \vH \, \xx^{1/2} + 2 \xx } \, , &  \xx \le \xxc \,, 
      \\
      \gmc/\muc + 2 (\xx - \xxc)/\gmc \, , & \xx > \xxc \, .
    \end{cases}
\label{zzxxppdBlroB}
\end{EQA}    
\end{theorem}

The value 
\( \zq(\BBH,\xx) = \sqrt{ \dimH + 2 \vH \, \xx^{1/2} + 2 \xx } \) can be upper bounded 
by \( \sqrt{\dimH} + \sqrt{2 \xx} \):
\begin{EQA}
\label{PzzxxpBroBu}
    \zqc(\BBH,\xx)
    & \leq &
    \begin{cases}
        \sqrt{\dimH} + \sqrt{2 \xx}, &  \xx \le \xxc, \\
        \zqc + 2 (\xx - \xxc)/\gmc , & \xx > \xxc.
    \end{cases}
\end{EQA}    

Depending on the value \( \xx \), we have two types of tail behavior of the 
quadratic form \( \| \xiv \|^{2} \). 
For \( \xx \le \xxc \), we have essentially the same deviation bounds as in the Gaussian case
with the extra-factor two in the deviation probability.
For \( \xx > \xxc \), we switch to the special regime driven by the exponential moment
condition \eqref{expgamgm}.
Usually \( \gmb^{2} \) is a large number (of order \( n \) in the i.i.d. setup) and 
the second term in \eqref{PxivbzzBBroB} can be simply ignored. 

As a corollary, we state the result for the norm of \( \xiv \in \R^{\dimp} \) corresponding to the case
\( \HPc^{-2} = \BB = \Id_{\dimp} \) and \( \dimp < \infty \). 
Then 
\begin{EQA}
	\dimA
	&=&
	\dimH
	=
	\vH^{2}
	=
	\dimp .
\label{pApHpvA2l1}
\end{EQA}

\begin{corollary}
\label{LLbrevelocro}   
Let \eqref{expgamgm} hold with \( \HPc = \Id_{\dimp} \) and 
\( \gm \geq 2 \sqrt{\dimp} \).
Then for each \( \xx > 0 \)
\begin{EQA}
    \P\bigl( \| \xiv \| \ge \zqc(\dimp,\xx) \bigr)
    & \le &
    2 \ex^{-\xx} + \ex^{-\xxc } \Ind(\xx < \xxc) ,
\label{PxivbzzBBro}
\end{EQA}    
where \( \zqc(\dimp,\xx) \) is defined by
\begin{EQA}
\label{PzzxxpBro}
    \zqc(\dimp,\xx)
    & \eqdef &
    \begin{cases}
      \bigl( \dimp + 2 \sqrt{\dimp \, \xx} + 2 \xx\bigr)^{1/2}, &  \xx \le \xxc  , \\
      \gmc/\muc + 2 \gmc^{-1} (\xx - \xxc)   , & \xx > \xxc .
    \end{cases}
\label{zzxxppdBlro}
\end{EQA}    
\end{corollary}

\begin{proof}[Proof of Theorem~\ref{LLbrevelocroB}]
%
First we consider the most interesting case \( \xx \leq \xxc \). 
We expect to get Gaussian type deviation bounds for such \( \xx \).
%
The main tool for the proof is the following lemma.
\begin{lemma}
\label{LGDBqfexpB}
Let \( \mu \in (0,1) \) and \( \zz(\mu) = \gmb / \mu - \sqrt{\dimH/\mu} > 0 \).
Then \eqref{expgamgm} implies
\begin{EQA}
	\E \exp\bigl( \mu \langle \BB \xiv, \xiv \rangle / 2 \bigr)
	\Ind\bigl( \| \HPc \BB \xiv \| \leq \zz(\mu) \bigr)
	& \leq &
	2 \exp \Bigl( \frac{\mu^{2} \vH^{2}}{4 (1 - \mu)} + \frac{\mu \, \dimH}{2} \Bigr) .
\label{wBmu2vA241mmutrB}
\end{EQA}
\end{lemma}

\begin{proof}
Let us fix for a moment some \( \xiv \in \R^{\dimp} \) and \( \mu < 1 \) and
define 
\begin{EQA}
	\av = \HPc^{-1} \xiv, 
	&\qquad & 
	\Sigma = \mu \BBH = \mu \HPc \BB \HPc .
\label{aSigmHm1xivmuc}
\end{EQA}
Consider the Gaussian measure \( \P_{\av,\Sigma} = \ND(\av,\Sigma^{-1}) \), and 
let \( \Uv \sim \ND(0,\Sigma^{-1}) \).
By the Girsanov formula
\begin{EQA}
	\log \frac{d\P_{\av,\Sigma}}{d\P_{0,\Sigma}}(\uv)
	&=&
	\langle \Sigma \av, \uv \rangle - \frac{1}{2} \bigl\langle \Sigma \av, \av \bigr\rangle
\label{GirslogddPaSSm1}
\end{EQA}
and for any set \( A \in \R^{\dimp} \)
\begin{EQA}
	\P_{\av,\Sigma}(A)
	&=&
	\P_{0,\Sigma}(A - \av)
	=
	\E_{0,\Sigma} \Bigl[ 
		\exp\Bigl\{ \langle \Sigma \Uv, \av \rangle 
			- \frac{1}{2} \bigl\langle \Sigma \av, \av \bigr\rangle \Bigr\} \Ind(A) 
	\Bigr] .
\label{PaSAP0SAav}
\end{EQA}
Now we select \( A = \bigl\{ \uv \colon \| \Sigma \uv \| \leq \gmb \bigr\} \).
Under \( \P_{0,\Sigma} \), one can represent 
\( \Sigma \Uv = \Sigma^{1/2} \gammav 
= \mu^{1/2} \HPc \BB^{1/2} \gammav \) with a standard Gaussian 
\( \gammav \).
Therefore,
\begin{EQA}
	\P_{0,\Sigma}(A - \av)
	&=&
	\Pg\bigl( \| \Sigma^{1/2} (\gammav - \Sigma^{1/2} \av) \| \leq \gmb \bigr)
	\\
	& \geq &
	\Pg\bigl( \| \Sigma^{1/2} \gammav \| \leq \gmb - \| \Sigma \av \| \bigr) .
\label{P0SAmaPgSm12}
\end{EQA}
We now use that \( \Pg\bigl( \| \Sigma^{1/2} \gammav \|^{2} \leq \tr (\Sigma) \bigr) \geq 1/2 \)
with \( \tr(\Sigma) = \mu \tr(\BBH ) = \mu \, \dimH \).
Therefore, the condition \( \| \Sigma \av \| + \sqrt{\mu \, \dimH} \leq \gmb \) implies
in view of \( \langle \Sigma \av, \av \rangle = \mu \langle \BB \xiv, \xiv \rangle \)
\begin{EQA}
	1/2
	\leq 
	\P_{\av,\Sigma}(A)
	&=&
	\E_{0,\Sigma} \Bigl[ 
		\exp\Bigl\{ \langle \Sigma \Uv, \HPc^{-1} \xiv \rangle
			- \mu \langle \BB \xiv, \xiv \rangle / 2 \Bigr\} 
			\Ind( \| \Sigma \Uv \| \leq \gmb) 
	\Bigr]
\label{12PavSAE0SISm1}
\end{EQA}
or
\begin{EQA}
	&& \nquad
	\exp\bigl( \mu \langle \BB \xiv, \xiv \rangle / 2 \bigr)
	\Ind\bigl( \| \Sigma \HPc^{-1} \xiv \| \leq \gmb - \sqrt{\mu \, \dimH} \bigr)
	\\
	& \leq &
	2 \E_{0,\Sigma} \Bigl[ 
		\exp\Bigl\{ \langle \Sigma \Uv, \HPc^{-1} \xiv \rangle
			\Ind( \| \Sigma \Uv \| \leq \gmb) 
	\Bigr] .
\label{emNxx22E0Sm1u}
\end{EQA}
We now take the expectation of the each side of this equation w.r.t. \( \xiv \),
change the integration order, and use \eqref{expgamgm} yielding
\begin{EQA}
	&& \nquad
	\E \exp\bigl( \mu \langle \BB \xiv, \xiv \rangle / 2 \bigr)
	\Ind\bigl( \| \Sigma \HPc^{-1} \xiv \| \leq \gmb - \sqrt{\mu \, \dimH} \bigr)
	\leq 
	2 \E_{0,\Sigma} \exp\bigl( \| \Sigma \Uv \|^{2} / 2 \bigr)
	\\
	&=&
	2 \Eg \exp\bigl( \mu \| \BBH^{1/2} \gammav \|^{2} / 2 \bigr)
	=
	2 \det\bigl( \Id - \mu \BBH \bigr)^{-1/2} .
\label{IHm1B12EgEgg}
\end{EQA}
We also use that for any \( \mu > 0 \)
\begin{EQA}
	\log \det\bigl( \Id - \mu \BBH \bigr)^{-1/2} - \frac{\mu \tr \BBH}{2}
	& \leq &
	\frac{\mu^{2} \tr \BBH^{2}}{4 (1 - \mu)} \, ;
\label{mu2v241mmiulogIm12}
\end{EQA}
see \eqref{m2v241m4mj1p}, 
and the first statement follows in view of \( \Sigma \HPc^{-1} \xiv = \mu \HPc \BB \xiv \).
\end{proof}

The use of \( \mu = \mux \) from \eqref{1v2sxm12m1m} in \eqref{wBmu2vA241mmutrB} yields 
similarly to the proof of Theorem~\ref{TexpbLGA}
\begin{EQA}
	\P\Bigl( \langle \BB \xiv, \xiv \rangle > \zq^{2}(\BBH,\xx), \,
		\| \HPc \BB \xiv \| \leq \zz(\mux)
	 \Bigr)
	 & \leq &
	 2 \ex^{-\xx} .
\label{2emxPblrHm1B}
\end{EQA}
It remains to consider the probability of large deviation
\( \P\bigl( \| \HPc \BB \xiv \| > \zz(\mux) \bigr) \).

\begin{lemma}
\label{Ldvbetagmb}
For any \( \xxc > 0 \) such that \( \zq(\BBH,\xxc) + 1 \leq \gmc/\muc \), 
it holds with \( \muc = \bigl\{ 1 + \vH/(2\sqrt{\xxc}) \bigr\}^{-1} \) 
and \( \zqc = \zz(\muc) = \gmb/\muc - \sqrt{\dimH/\muc} \)
\begin{EQA}
	\P\bigl( \| \HPc \BB \xiv \| > \zqc \bigr) 
	& \leq &
	\P\bigl( \langle \BB \xiv, \xiv \rangle > \zqc^{2} \bigr)
	\leq 
	\ex^{-\xxc} .
\label{exmxxlPHm1Bxigm}
\end{EQA}
\end{lemma}

\begin{proof}
It follows due to \eqref{1v2sxm12m1m} and \eqref{m2vA241muz2} for any \( \mu \leq \muc \)
\begin{EQA}
	\frac{\mu^{2} \vH^{2}}{4 (1 - \mu)} + \frac{\mu \, \dimH}{2}
	\leq 
	\frac{\muc^{2} \, \vH^{2}}{4 (1 - \muc)} + \frac{\muc \, \dimH}{2}
	& \leq &
	\frac{\muc \zq^{2}(\BBH,\xxc)}{2} - \xxc ,
\label{mumucmuz2Bxmx24}
\end{EQA}
where the right hand-side does not depend on \( \mu \).
We now define \( \eta^{2} = \langle \BB \xiv, \xiv \rangle \) and use that 
\( \| \HPc \BB \xiv \| \leq \| \BB^{1/2} \xiv \| = \eta \).
Then by \eqref{wBmu2vA241mmutrB}
\begin{EQA}
	\E \exp(\mu \eta^{2}/2) \Ind\bigl( \eta \leq \zz(\mu) \bigr)
	& \leq &
	\exp \bigl\{ \muc \zq^{2}(\BBH,\xxc)/2 - \xxc \bigr\} ,
\label{me22Iezmemcz2}
\end{EQA}
Define the inverse function \( \mu(\zz) = \zz^{-1}(\mu) \).
For any \( \zz \geq \zqc \), it follows from \eqref{me22Iezmemcz2} with 
\( \mu = \mu(\zz) \)
\begin{EQA}
	\E \exp\bigl\{ \mu(\zz) (\zz-1)^{2} /2 \bigr\} 
		\Ind\bigl( \zz - 1 \leq \eta \leq \zz \bigr)
	& \leq &
	\exp \bigl\{ \muc \zq^{2}(\BBH,\xxc)/2 - \xxc \bigr\} 
\label{emcz2NBxm2m1I}
\end{EQA}
yielding
\begin{EQA}
	\P\bigl( \zz - 1 \leq \eta \leq \zz \bigr)
	& \leq &
	\exp \Bigl( 
		- \mu(\zz) \, (\zz - 1)^{2}/2 + \muc \zq^{2}(\BBH,\xxc)/2 - \xxc  
	\Bigr)
\label{mzm12mcz2Bxm2x}
\end{EQA}
and hence,
\begin{EQA}
	\P\bigl( \eta > \zz \bigr)
	& \leq &
	\int_{\zz}^{\infty} \exp\bigl\{ - \mu(\zq) (\zq-1)^{2}/2 + \muc \zq^{2}(\BBH,\xxc)/2 - \xxc \bigr\}
		d\zq
\label{zzzzcm1Pezmcl22}
\end{EQA}
Further, \( \mu \, \zz(\mu) = \gmb - \sqrt{\dimH \mu} \) and
\begin{EQA}
	\gmc
	=
	\muc \, \zqc
	& \leq &
	\mu \, \zz(\mu)
	\leq 
	\gmb,
	\quad
	\mu \leq \muc .
\label{muxzxmzmglex}
\end{EQA}
This implies the same bound for the inverse function:
\begin{EQA}
	\gmc
	& \leq &
	\zz \, \mu(\zz)
	\leq 
	\gmb,
	\quad
	\zz \geq \zqc \, ,
\label{muxzxmzmglexi}
\end{EQA}
and 
\begin{EQA}
	\P\bigl( \eta > \zz \bigr)
	& \leq &
	\int_{\zz}^{\infty} 
		\exp\bigl\{ - \mu(\zq) \bigl( \zq^{2}/2 - \zq \bigr) + \muc \zq^{2}(\BBH,\xxc)/2 - \xxc \bigr\} d\zq
	\\
	& \leq &
	\int_{\zz}^{\infty} 
		\exp\bigl\{ - \gmc \, \zq/2 + \gmb + \muc \, \zq^{2}(\BBH,\xxc)/2 - \xxc \bigr\} d\zq
	\\
	& \leq &
	\frac{2}{\gmc}
		\exp\bigl\{ - \gmc \, \zz / 2 + \gmb + \muc \, \zq^{2}(\BBH,\xxc)/2 - \xxc \bigr\} .
\label{Perzxinzinmmuz2gmz2}
\end{EQA}
Conditions \( \gmc \zqc = \muc^{-1} \gmc^{2} \geq \muc \bigl\{ \zq(\BBH,\xxc) + 1 \bigr\}^{2} \)
and \( \gmc \geq 1 \) 
ensure that
\( \P\bigl( \eta > \zqc \bigr) \leq \ex^{-\xxc} \).
\end{proof}

Remind that \( \xxc \) is the largest \( \xx \)-value ensuring the condition
\( \gmc \geq \zq(\BBH,\xxc) + 1 \).
We also use that for \( \xx \leq \xxc \), it holds \( \zz(\mux) \geq \zz(\muc) = \zqc \).
Therefore, by \eqref{2emxPblrHm1B} and Lemma~\ref{Ldvbetagmb} 
\begin{EQA}
	\P\bigl( \langle \BB \xiv, \xiv \rangle \geq \zq^{2}(\BBH,\xx) \bigr)
	& \leq &
	\P\bigl( \langle \BB \xiv, \xiv \rangle \geq \zq^{2}(\BBH,\xx), \| \HPc \BB \xiv \| \leq \zz(\mux) \bigr)
	+ \P\bigl( \langle \BB \xiv, \xiv \rangle \geq \zqc^{2} \bigr)
	\\
	& \leq &
	2 \ex^{-\xx} + \ex^{-\xxc} \, .
\label{2emxpemxcPBxx}
\end{EQA}
Finally we consider \( \xx > \xxc \).
Applying \eqref{Perzxinzinmmuz2gmz2} yields by \( \zz \geq \zqc \)
\begin{EQA}
	\P\bigl( \eta > \zz \bigr)
	& \leq &
	\frac{2}{\muc \, \zqc}
		\exp\bigl\{ - \muc \, \zqc^{2} / 2 + \gmb + \muc \, \zq^{2}(\BBH,\xxc)/2 - \xxc \bigr\}
		\exp\bigl\{ - \muc \, \zqc (\zz - \zqc) / 2  \bigr\}
	\\
	& \leq &
	\ex^{- \xxc} \exp\bigl\{ - \gmc (\zz - \zqc) / 2  \bigr\} .
\label{emxcmcmzczzc2}
\end{EQA}
The choice \( \zz \) by 
\begin{EQA}
	\gmc (\zz - \zqc) / 2 
	&=&
	\xx - \xxc
\label{mzzc2xxczz}
\end{EQA}
ensures the desired bound.

%
\end{proof}

\subsection{Taylor expansions}
\label{STaylor}
Here we collect some useful bounds for various Taylor-type expansions for
a smooth function. 
Let \( f \) be a four time differentiable function on \( \R^{\dimp} \).
Here \( \dimp \leq \infty \).
By \( f^{(k)}(\xv,\uv) \) we denote the \( m \)th directional derivative at \( \xv \):
\begin{EQA}
    f^{(k)} (\xv,\uv)
    & \eqdef &
    \frac{d^{k}}{dt^{k}} f(\xv + t \uv) \bigg|_{t=0} \, .
\label{fmxudmdtmfxtu}
\end{EQA}
In particular,
\( f'(\xv,\uv) = \bigl\langle \nabla f(\xv), \uv \bigr\rangle \) and 
\( f''(\xv,\uv) = \bigl\langle \nabla^{2} f(\xv) \, \uv, \uv \bigr\rangle \).
Below we assume that some open set \( \Xs \subseteq \R^{\dimp} \) 
is fixed, and, in addition, for each \( \xv \in \Xs \), 
and a centrally symmetric convex set \( \UV(\xv) \) are fixed and
\begin{EQA}
    \frac{1}{k!} \bigl| f^{(k)}(\xv, \uv) \bigr|
    =
    \delta_{k}(\xv,\uv)
    & \leq &
    \delta_{k} \, ,
    \quad
    \xv \in \Xs, \uv \in \UV ,
    \quad
    k=3,4
\label{1mffmxum34}
\end{EQA}
for some constants \( \delta_{k} \) depending on \( \Xs \) and \( \UV \).
All bounds will be given in terms of \( \delta_{3} \) and \( \delta_{4} \).
The construction can be extended by making \( \UV \) dependent on \( \xv \in \Xs \)
at cost of more complicated notation. 

\begin{lemma}
\label{Lextftumu}
Suppose \eqref{1mffmxum34} with \( \delta_{k} \leq 1 \) for \( k=3,4 \).
Then for any point \( \xv \in \Xs \) 
\begin{EQA}
    && \nquad
    \left| \frac{1}{2} \left( \ex^{f(\xv + \uv) - f(\xv) - f'(\xv,\uv)} + \ex^{f(\xv - \uv) - f(\xv) + f'(\xv,\uv)} \right) 
    - \ex^{f''(\xv,\uv)/2} \right|
    \\
    & \leq &
    \ex^{f''(\xv,\uv)/2} \,\bigl( 4 \delta_{3}^{2} + 4 \delta_{4} \bigr) .
\label{12efxufxmufpp2}
\end{EQA}
Furthermore, 
\begin{EQA}
    \left| \ex^{f(\xv + \uv) - f(\xv) - f'(\xv,\uv)} - \ex^{f''(\xv,\uv)/2} \right|
    & \leq &
    \delta_{3} \, \ex^{f''(\xv,\uv)/2} \, .
\label{12efxufxmufpp23}
\end{EQA}    
\end{lemma}

\begin{proof}
Taylor expansions of the forth order imply 
\begin{EQA}
    f(\xv + \uv) - f(\xv) - f'(\xv, \uv) - \frac{1}{2} f''(\xv, \uv) - \frac{1}{6} f^{(3)}(\xv,\uv)  
    & = &
    \rho_{1} \, ,
    \qquad
    |\rho_{1}|
    \leq 
    \delta_{4} \, ,
    \\
    f(\xv - \uv) - f(\xv) + f'(\xv, \uv) - \frac{1}{2} f''(\xv, \uv) + \frac{1}{6} f^{(3)}(\xv,\uv) 
    & = &
    \rho_{2} \, ,
    \qquad
    |\rho_{2}|
    \leq 
    \delta_{4} \, .
\label{fu12un2fud2}
\end{EQA}
Further, define \( \kappa = f^{(3)}(\xv,\uv)/6 \), so that 
\( |\kappa| \leq \delta_{3} \leq 1 \).
Then
\begin{EQA}
    && \nquad
    \ex^{f(\xv + \uv) - f(\xv) - f'(\xv,\uv)} 
    + \ex^{f(\xv - \uv) - f(\xv) + f'(\xv,\uv)} - 2 \ex^{f''(\xv,\uv)/2}
    \\
    &=&
    \ex^{f''(\xv,\uv)/2} \left( \ex^{\kappa + \rho_{1}} + \ex^{-\kappa + \rho_{2}} - 2 \right) .
\label{12f1fm1fpp2m2}
\end{EQA}
The function 
\begin{EQA}
      g(s)
      & \eqdef &
      \frac{1}{2} \exp\bigl( s \, \kappa + \rho_{1} \bigr) 
      + \frac{1}{2} \exp\bigl( - s \, \kappa + \rho_{2} \bigr) - 1 
\label{gt12etKd1d2}
\end{EQA}
fulfills
\begin{EQA}
    |g(0)|
    &=&
    \Bigl| \frac{1}{2} \ex^{\rho_{1}} + \frac{1}{2} \ex^{\rho_{2}} - 1 \bigr|
    \leq 
    |\rho_{1}| + |\rho_{2}| , 
    \\
    |g'(0)|
    &=&
    \frac{1}{2}
    \bigl| \kappa \bigl( \ex^{\rho_{1}} - \ex^{\rho_{2}} \bigr) \bigr|
    \leq 
    |\rho_{1}| + |\rho_{2}| 
\label{g012ed12D12K}
\end{EQA}
and for any \( s \in [0,1] \) by simple algebra due to 
\( |\kappa| \leq 1 \) and \( |\rho_{k}| \leq 1 \) for
\( m=1,2 \)
\begin{EQA}
    |g''(s)|
    &=&
    \frac{1}{2} \Bigl| 
    \kappa^{2} \Bigl\{ \exp\bigl( s \, \kappa + \rho_{1} \bigr) 
    + \exp\bigl( - s \, \kappa + \rho_{2} \bigr) \Bigr\} 
    \Bigr|
    \\
    & \leq &
    \frac{|\kappa|^{2} \ex}{2} \, \bigl( \ex^{ |\kappa|} + \ex^{-|\kappa|} \bigr)
    < 
    8 | \kappa |^{2} ,
\label{gpptK2gt12eD}
\end{EQA}
and thus
\begin{EQA}
    \bigl| g(1) \bigr|
    & \leq &
    \sup_{s \in [0,1]} \bigl| g(0) + g'(0) + \frac{1}{2} g''(s) \bigr|
    \leq 
    4 | \kappa |^{2} + 2 |\rho_{1}| + 2 |\rho_{2}| ,
\label{g1supt01g0gp0}
\end{EQA}
and \eqref{12efxufxmufpp2} follows. 
The bound \eqref{12efxufxmufpp23} can be obtained in a similar way 
using the Taylor expansion of the third order.
\end{proof}

Now we study the modulus of continuity for the gradient \( \nabla f(\xv) \) and the Hessian
\( \nabla^{2} f(\xv) \).

\begin{lemma}
\label{LHessTay}
Suppose \eqref{1mffmxum34} with \( \delta_{3} \leq 1 \).
Let \( \xv \in \Xs \) and \( \uv \in \UV \) be such that
\( \xv + \uv \in \Xs \). 
Then, for any \( \wv \in \UV \)
\begin{EQA}
\label{wTnfxun2fu}
    \Bigl| \bigl\langle  \wv, \nabla f(\xv+\uv) - \nabla f(\xv) - \nabla^{2}f(\xv) \uv \bigr \rangle 
    \Bigr|
    & \leq &
    \CONST \delta_{3} \, ,
    \\
    \Bigl| \bigl\langle \wv, \bigl\{ \nabla^{2} f(\xv+\uv) - \nabla^{2}f(\xv) \bigr\} \wv \bigr\rangle \Bigr|
    & \leq &
    \CONST \delta_{3} \, .
\label{wTn2fxun2fw}
\end{EQA}
\end{lemma}

\begin{proof}
Let us fix any \( \xvd \in \Xs \) and \( \wvd \in \UV \) and define the function 
\begin{EQA}
	g(t)
	& \eqdef &
	f(\xvd + t \wvd) + f(\xvd - t \wvd) - 2 f(\xvd) - t^{2} f''(\xvd,\wvd) .
\label{fxdtwmtv222}
\end{EQA}
The Taylor expansion of the third order yields  
\begin{EQA}
    \bigl| g(1) \bigr|
    =
    \Bigl| f(\xvd+\wvd) + f(\xvd-\wvd) - 2 f(\xvd) - f''(\xvd,\wvd) \Bigr| 
    & \leq &
    2 \delta_{3}(\xvd,\wvd) \, .
\label{fxw2fxmw2fp2}
\end{EQA}
We apply this bound for \( \xvd = \xv \) and \( \xvd = \xv + \uv \) and take the difference 
between them. 
This implies
\begin{EQA}
    &&
    \bigl| f''(\xv,\wvd) - f''(\xv + \uv,\wvd) \bigr|
    \leq 
    \bigl| f(\xv+\wvd) + f(\xv-\wvd) - 2 f(\xv) 
    \\
    &&
    \qquad
    - \, f(\xv + \uv + \wvd) - f(\xv + \uv - \wvd) + 2 f(\xv + \uv) \bigr|
    + 2 \delta_{3}(\xv,\wvd) + 2 \delta_{3}(\xv+\uv,\wvd) \, . 
    \qquad \qquad
\label{12fxw224}
\end{EQA}
For given \( \xv, \uv, \wv \), and \( \xb = \xv + \uv/2 \), define 
\begin{EQA}
    g(t)
    & \eqdef &
    f\bigl( \xb + t(\uv+\wv) \bigr) - f\bigl( \xb - t(\uv+\wv) \bigr)
    \\
    &&
    + \, f\bigl( \xb + t(\uv-\wv) \bigr) - f\bigl( \xb - t(\uv-\wv) \bigr)
    - 2 f\bigl( \xb + t\uv \bigr) + 2 f\bigl( \xb - t\uv \bigr) .
\label{gtfxtuw22fxtu2}
\end{EQA}
It is straightforward to see that \( g(0) = g'(0) = g''(0) = 0 \).
Moreover, in view of \( \uv \in \UV \) and
\( (\uv \pm \wv)/2 \in \UV \), it holds \( \delta_{3}(\xb,\uv/2) = \delta_{3}(\xb,\uv)/8 \) and
for any \( |t| \leq 1/2 \)
\begin{EQA}
    \frac{1}{6} \bigl| g^{(3)}(t) \bigr|
    & \leq &
    \frac{5 \delta_{3}}{2} \, .
\label{124g46d4}
\end{EQA}
By Taylor expansion of the third order we derive
\begin{EQA}
    \bigl| g(1/2) \bigr|
    & \leq &
    \sup_{t \in [0,1]} \frac{1}{6} \bigl| g^{(3)}(t) \bigr|
    \leq 
    \frac{5 \delta_{3}}{2} \, . 
\label{g1t01124g4t6}
\end{EQA}
Note that \( g(1/2) \) is exactly the expression in the right hand-side of \eqref{12fxw224} with \( \wvd = \wv/2 \).
The use of \( \delta_{3}(\xvd,\wvd) = \delta_{3}(\xvd,\wv)/8 \)
together with \eqref{12fxw224} yields \eqref{wTn2fxun2fw} with \( \CONST = 3 \).
\end{proof}

Now we specify the result to the case of an elliptic set \( \UV \) of the form
\begin{EQA}
    \UV 
    &=& 
    \bigl\{ \uv \colon \| \QP \uv \| \leq \rr \bigr\} 
\label{UVdefQPur}
\end{EQA}
for a positive invertible operator \( \QP \) and \( \rr > 0 \).

\begin{lemma}
\label{LellUVD2w}
Let \( \UV \) be given by \eqref{UVdefQPur} with \( \QP > 0 \), and let \( \xv \in \Xs \) and 
\( \uv \in \UV \) be such that \( \xv + \uv \in \Xs \).
Then
\begin{EQA}
\label{Qm1n2fxn2txurm1}
    \bigl\| 
    \QP^{-1} \bigl\{ \nabla f(\xv + \uv) - \nabla f(\xv) - \nabla^{2} f(\xv) \uv \bigr\} 
    \bigr\|
    & \leq &
    \CONST \rr^{-1} \delta_{3} \, ,
    \\
    \bigl\| \QP^{-1} \bigl\{ \nabla^{2} f(\xv) - \nabla^{2} f(\xv + \uv) \bigr\} \QP^{-1} \bigr\|
    & \leq &
    \CONST \rr^{-2} \delta_{3} \, .
\label{Qm1n2fxn2txurm2}
\end{EQA}
\end{lemma}

\begin{proof}
For any \( \wv \in \UV \), it holds by Lemma~\ref{LHessTay}
\begin{EQA}
    \Bigl| 
    	\bigl\langle \wv, \bigl\{ \nabla^{2} f(\xv+\uv) - \nabla^{2} f(\xv) \bigr\} \wv \bigr\rangle 
    \Bigr|
    & = &
    \Bigl| \bigl\langle \QP \wv, \QP^{-1} \bigl\{ \nabla^{2} f(\xv+\uv) - \nabla^{2} f(\xv) \bigr\} 
    \QP^{-1} (\QP \wv) \bigr\rangle \Bigr|
    \leq 
    \CONST \delta_{3} \, .
\label{wTn2fxun2fw3}
\end{EQA}
As this bound holds for all \( \wv \in \UV \) with \( \| \QP \wv \| \leq \rr \), the result follows.
\end{proof}

The result of Lemma~\ref{Lextftumu} can be extended to the integral of \( \ex^{f(\xv + \uv)} \)
over \( \uv \in \UV \).

\begin{lemma}
\label{Lintfxupp2}
Let \( \UV \) be a subset in \( \R^{\dimp} \).
Suppose \eqref{1mffmxum34} with \( \delta_{k} \leq 1 \) for \( k=3,4 \).
Then for any point \( \xv \in \Xs \) and any centrally symmetric set \( A \subset \UV \)
\begin{EQA}
\label{4d324d4efppm}      
    \biggl| 
    	\int_{A} \ex^{f(\xv+\uv) - f(\xv) - f'(\xv,\uv)} \, d\uv 
    	- \int_{A} \ex^{f''(\xv,\uv)/2} \, d\uv \biggr|
    & \leq &
    \err \int_{A} \ex^{f''(\xv,\uv)/2} \, d\uv \, 
    \qquad
\end{EQA}
with \( \err = 4 \delta_{3}^{2} + 4 \delta_{4} \) and for any vector \( \zv \)
\begin{EQA}
    && \nquad
    \biggl| 
    	\int_{A} \langle \zv,\uv \rangle^{2} \ex^{f(\xv+\uv) - f(\xv) - f'(\xv,\uv)} \, d\uv 
    		- \int_{A} \langle \zv,\uv \rangle^{2}\ex^{f''(\xv,\uv)/2} \, d\uv 
    \biggr|
    \\
    & \leq &
    \err \int_{A} \langle \zv,\uv \rangle^{2} \ex^{f''(\xv,\uv)/2} \, d\uv \, .
\label{4d324d4efppmz2}      
\end{EQA}
If \( A \) is not centrally symmetric then
\begin{EQA}
\label{d324d4efppm3}      
    \biggl| \int_{A} \ex^{f(\xv+\uv) - f(\xv) - f'(\xv,\uv)} \, d\uv - \int_{A} \ex^{f''(\xv,\uv)/2} \, d\uv \biggr|
    & \leq &
    \delta_{3} \int_{A} \ex^{f''(\xv,\uv)/2} \, d\uv \, .
    \qquad
\end{EQA}
\end{lemma}

\begin{proof}
By symmetricity of \( \UV \), it holds 
\begin{EQA}
    \int_{A} \ex^{f(\xv+\uv) - f(\xv) - f'(\xv,\uv)} \, d\uv
    &=&
    \frac{1}{2} \int_{A}\left( \ex^{f(\xv+\uv) - f(\xv) - f'(\xv,\uv)} 
    + \ex^{f(\xv-\uv) - f(\xv) + f'(\xv,\uv)} \right) d\uv ,
\label{12efuAefmuA}
\end{EQA}
and the first result is proved by \eqref{12efxufxmufpp2}.
The same symmetricity arguments apply to \eqref{4d324d4efppmz2}.
The final bound for any \( A \) follows from \eqref{12efxufxmufpp23}.
\end{proof}

The bound \eqref{4d324d4efpp2du} can be specified to the case of a massive set \( \UV \).
%
We assume that \( f \) is concave and \( \HP^{2} \eqdef - \nabla^{2} f(x) \geq 0 \).

\begin{lemma}
\label{Lintfxupp2r}
Let \( f(\cdot) \) be strictly concave with \( \HP^{2} = - \nabla^{2} f(\xv) > 0 \). 
Suppose \eqref{1mffmxum34} with \( \delta_{3} \leq 1 \).
For a linear operator \( \QP \), it holds
\begin{EQA}
    \biggl\| 
    	\int_{\UV} \QP \uv \, \ex^{f(\xv+\uv) - f(\xv) - f'(\xv,\uv)} \, d\uv 
    \biggr\|
    & \leq &
    \delta_{3} \int_{\UV} \| \QP \uv \|\ex^{f''(\xv,\uv)/2} \, d\uv \, ,
\label{4d324d4efpp2du}
\end{EQA}    
Let also \( \UV \) be massive in the sense that 
\begin{EQA}
	\P\bigl( \HP^{-1} \gammav \in \UV \bigr) 
	& \geq &
	1/2 
\label{PHm1gU2m12}
\end{EQA}
with \( \gammav \) standard normal in \( \R^{\dimp} \).
    Then 
for any linear operator \( \QP \), it holds
\begin{EQA}
    \biggl\| 
    	\frac{\int_{\UV} \QP \uv \, \ex^{f(\xv+\uv) - f(\xv) - f'(\xv,\uv)} \, d\uv}
    			 {\int_{\UV} \ex^{f''(\xv,\uv)/2} \, d\uv}
    \biggr\|
    & \leq &
    2 \delta_{3} \, \E \| \QP \, \HP^{-1} \gammav \| \, .
\label{4d324d4efpp2dur}
\end{EQA}    
  \end{lemma}

\begin{proof}
The bound \eqref{4d324d4efpp2du} follows in a way similar to the proof
of Lemma~\ref{Lintfxupp2} using \eqref{12efxufxmufpp23} instead of \eqref{12efxufxmufpp2}.
We apply \eqref{4d324d4efpp2du} 
yielding in view of \( f''(\xv,\uv) = - \| \HP \uv \|^{2} \)
\begin{EQA}
    \biggl\| 
    	\frac{\int_{\UV} \QP \uv \, \ex^{f(\xv+\uv) - f(\xv) - f'(\xv,\uv)} \, d\uv}
    			 {\int_{\UV} \ex^{f''(\xv,\uv)/2} \, d\uv}
    \biggr\|
    & \leq &
    \delta_{3}  
    	\frac{\int_{\UV} \| \QP \uv \| \, \ex^{f''(\xv,\uv)/2} \, d\uv}
    			 {\int_{\UV} \ex^{f''(\xv,\uv)/2} \, d\uv} 
    \\
    & \leq &
    2 \delta_{3}  
    	\frac{\int_{\UV} \| \QP \uv \| \, \ex^{f''(\xv,\uv)/2} \, d\uv}
    			 {\int \ex^{f''(\xv,\uv)/2} \, d\uv} 
    \leq 
    2 \delta_{3} \E \| \QP \, \HP^{-1} \gammav \|
\label{d3Ageg2dg}
\end{EQA}
as required. 
\end{proof}

All the bounds presented above assume that \( \UV \) is a symmetric subset
of \( \R^{\dimp} \), in particular, an ellipsoid centred at zero.
Now we check what happens under a small departure from symmetricity.

\begin{lemma}
\label{Lintavshift}
Let \( f \) be concave with \( \HP^{2} = - \nabla^{2} f(\xv) \).
Let also \( \UV \) be centrally symmetric massive set; see \eqref{PHm1gU2m12}
and let \( \av \in \UV \) be fixed and 
\( \UV + \av \subset \CS(\rups) = \bigl\{ \uv \colon \| \HP \uv \| \leq \rups \bigr\} \). 
Suppose \eqref{1mffmxum34} with \( \delta_{k} = \delta_{k}(\rups) \leq 1 \) for \( k=3,4 \).
Then
\begin{EQA}
	&& \nquad
	\biggl| 
		{\int_{\UV} \ex^{f(\xv+\uv+\av ) - f(\xv) - f'(\xv,\uv + \av)} \, d\uv
				- \int_{\UV} \ex^{f''(\xv,\uv+\av)/2} \, d\uv} 
	\biggr|
	\\
	& \lesssim &
	\Bigl\{ 
		\err(\rups) +
		\bigl( \| \HP \av \| + \| \HP \av \|^{2} \bigr) \delta_{3}(\rups) 
	\Bigr\} 
	{\int_{\UV} \ex^{f''(\xv,\uv)/2} \, d\uv} \, .
\label{UVafxavfpau}
\end{EQA}
\end{lemma}

\begin{proof}
Define
\begin{EQA}
	h(t)
	& \eqdef &
	\detHP
	\int_{\UV} \ex^{f(\xv + \uv + t \av) - f(\xv) - f'(\xv,\uv + t \av)} \, d\uv \, 
	\\
	\hhh(t)
	& \eqdef &
	\detHP
	\int_{\UV} \ex^{f''(\xv,\uv + t \av)/2} \, d\uv \, 
\label{hadeefxavfxu}
\end{EQA} 
with 
\begin{EQA}
	\detHP
	& \eqdef &
	\left( \int_{\UV} \ex^{f''(\xv,\uv)/2} \, d\uv \right)^{-1} .
\label{detHPdefq0}
\end{EQA}
Then we have to bound the difference \( h(t) - \hhh(t) \) for \( t \leq 1 \).
For this, we bound the first two derivatives of \( h(t) \).
By \eqref{4d324d4efppm} Lemma~\ref{Lintfxupp2}, it holds \( |h(0) - \hhh(0)| \leq \err(\rups) \).
Further, in view of \( f''(\xv,\uv + t \av) = - \| \HP (\uv + t \av) \|^{2} \)
\begin{EQA}
\label{ddthtvt0efxhhh}
	h'(0) 
	&=&
	\detHP
	\int_{\UV} \bigl\langle \av, \nabla f(\xv + \uv) - \nabla f(\xv) \bigr\rangle
		\, \ex^{f(\xv+\uv) - f(\xv) - f'(\xv,\uv)} \, d\uv \, ,
	\\
	\hhh'(0) 
	&=&
	\detHP
	\int_{\UV} \bigl\langle \HP \av, \HP \uv \bigr\rangle
		\, \ex^{f''(\xv,\uv)/2} \, d\uv = 0 \, 
\label{ddthtvt0efx}
\end{EQA}
because \( \UV \) is centrally symmetric.
By \eqref{Qm1n2fxn2txurm1} of Lemma~\ref{LellUVD2w} and \eqref{4d324d4efppm} of Lemma~\ref{Lintfxupp2}
\begin{EQA}
	&& \nquad
	\left| 
	  \int_{\UV} \bigl\langle \av, \nabla f(\xv + \uv) - \nabla f(\xv) - \nabla^{2} f(\xv) \uv \bigr\rangle
		\, \ex^{f(\xv+\uv) - f(\xv) - f'(\xv,\uv)} \, d\uv 
	\right|
	\\
	& \leq &
	\CONST \, \| \HP \av \| \rups^{-1} \, \delta_{3}(\rups) \int_{\UV} \ex^{f(\xv+\uv) - f(\xv) - f'(\xv,\uv)} \, d\uv 
	\\
	& \leq &
	\CONST \, \| \HP \av \| \rups^{-1} \, \delta_{3}(\rups) \bigl( 1 + \err(\rups) \bigr) \int_{\UV} \ex^{f''(\xv,\uv)/2} \, d\uv
\label{d4d31e4fuva}
\end{EQA}
for \( \err (\rups) = 4 \delta_{3}^{2}(\rups) + 4 \delta_{4}(\rups) \).
Further we use \( f''(\xv,\uv) = - \| \HP \uv \|^{2} \), 
\( - \bigl\langle \av, \nabla^{2} f(\xv) \uv \bigr\rangle = \bigl\langle \HP \av, \HP \uv \bigr\rangle \),
and it follows by Lemma~\ref{Lintfxupp2r}  that
\begin{EQA}
	\left| 
		\detHP
		\int_{\UV} \bigl\langle \av, \nabla^{2} f(\xv) \uv \bigr\rangle
		\ex^{f(\xv+\uv) - f(\xv) - f'(\xv,\uv)} \, d\uv \right|
	& \leq &
	2 \delta_{3}(\rups) \, \E \bigl| \bigl\langle \HP \av, \gammav \bigr\rangle \bigr|
	\leq 
	2 \delta_{3}(\rups) \| \HP \av \| .
\label{rd3iUefppxu2}
\end{EQA}
Putting this together with \eqref{ddthtvt0efxhhh}, \eqref{d4d31e4fuva} yields
\begin{EQA}
	\bigl| h'(0) \bigr|
	& \lesssim &
	\delta_{3}(\rups) \| \HP \av \| .
\label{hp0d3r0Hal}
\end{EQA}
For the second derivative, 
\begin{EQA}
	h''(t) 
	&=&
	\detHP
	\int_{\UV} \Bigl\{ \bigl\langle \av, \nabla f(\xv + \uv + t \av) - \nabla f(\xv) \bigr\rangle^{2}
		 + \bigl\langle \av, \nabla^{2} f(\xv + \uv + t \av) \av \bigr\rangle
		 \Bigr\} 
	\\
	&& \qquad \qquad \qquad
		\times \, \ex^{f(\xv + \uv + t \av) - f(\xv) - f'(\xv,\uv + t \av)} \, d\uv \, .
\label{ddthtvt0efx2}
\end{EQA}
Similarly, by the use of \( f''(\xv,\uv) = - \| \HP \uv \|^{2} \) we derive
\begin{EQA}
	\hhh''(t) 
	&=&
	\detHP
	\int_{\UV} \Bigl\{ \bigl\langle \av, \HP^{2}(\uv + t \av) \bigr\rangle^{2}
		 - \bigl\langle \av, \HP^{2} \av \bigr\rangle
		 \Bigr\} \, \ex^{- \| \HP (\uv + t \av ) \|^{2}/2} \, d\uv \, .
\label{ddthtvt0efx2}
\end{EQA}
Now by \eqref{4d324d4efpp2du} of Lemma~\ref{Lintfxupp2} 
\begin{EQA}
	&& \nquad
	\bigl| \bigl\langle \av, \nabla^{2} f(\xv + \uv + t \av) \av \bigr\rangle 
	- \bigl\langle \av, \nabla^{2} f(\xv) \av \bigr\rangle \bigr| 
	\\
	&=&
	\bigl| \bigl\langle \HP \av, \HP^{-1} \bigl\{ \nabla^{2} f(\xv + \uv + t \av) - \nabla^{2} f(\xv) \bigr\}\HP^{-1} \HP \av \bigr\rangle \bigr|
	\lesssim 
	\| \HP \av \|^{2} \rups^{-2} \delta_{3}(\rups) \, .
\label{avn2fxutaHHpr02}
\end{EQA}
Similarly
\begin{EQA}
	\bigl| \bigl\langle \av, \nabla f(\xv + \uv + t \av) - \nabla f(\xv) + \HP^{2} (\uv + t \av) \bigr\rangle \bigr|
	& \lesssim &
	\| \HP \av \| \rups^{-1} \delta_{3}(\rups)
\label{anfxutanfxH2}
\end{EQA}
and by \( \| \HP (\uv + t \av) \| \leq \rups \) it holds 
\( \bigl| \bigl\langle \av, \HP^{2} (\uv + t \av) \bigr\rangle \bigr| \leq \| \HP \av \| \rups \) and 
\begin{EQA}
	&& \nquad
	\bigl\langle \av, \nabla f(\xv + \uv + t \av) - \nabla f(\xv) \bigr\rangle^{2}
	- \bigl\langle \av, \HP^{2} (\uv + t \av) \bigr\rangle^{2}
	\\
	& \lesssim &  
	\| \HP \av \|^{2} \, \rups^{-2} \delta_{3}^{2}(\rups) + 2 \| \HP \av \| \, \rups^{-1} \delta_{3}(\rups) \| \HP \av \| \rups
	\\	
	& \lesssim &
	\| \HP \av \|^{2} \, \delta_{3}(\rups) .
\label{efxutdvfpxutw}
\end{EQA}
We conclude that
\begin{EQA}
	h''(t) 
	&=&
	\detHP
	\int_{\UV} \Bigl\{ \bigl\langle \av, \HP^{2}(\uv + t \av) \bigr\rangle^{2}
		 - \bigl\langle \av, \HP^{2} \av \bigr\rangle + \tau(\av,\uv,t)
		 \Bigr\} \, \ex^{f(\xv + \uv + t \av) - f(\xv) - f'(\xv,\uv + t \av)} \, d\uv
\label{hpptq0efxvtafpfpp}
\end{EQA}
where \( |\tau(\av,\uv,t)| \lesssim \| \HP \av \|^{2} \, \delta_{3}(\rups) \).
Lemma~\ref{Lextftumu} helps to bound 
\begin{EQA}
	\bigl| h''(t) - \hhh''(t) \bigr|
	& \lesssim &
	\| \HP \av \|^{2} \, \delta_{3}(\rups)
\label{Dw2iRpDwg222}
\end{EQA}
uniformly in \( |t| \leq 1 \).
This yields with some \( \rho \in [0,1] \)
\begin{EQA}
	\bigl| h(1) - \hhh(1) \bigr|
	& \leq &
	\bigl| h(0) - \hhh(0) \bigr|
	+ \bigl| h'(0) - \hhh'(0) \bigr| + \bigl| h''(\rho) - \hhh''(\rho) \bigr|/2
	\\
	& \lesssim &
	\err(\rups) 
	+ \bigl( \| \HP \av \| + \| \HP \av \|^{2} \bigr) \delta_{3}(\rups)
\label{Had3ha2hpp}
\end{EQA}
which completes the proof.
\end{proof}


\subsection{Concavity and tail bounds}
\label{SGaussintegr}

Let \( f(\xv) \) be a function on \( \R^{\dimp} \).
Previous results describe the local behavior of \( f(\xv + \uv) \) for \( \uv \in \UV \)
under local smoothness conditions. 
Now we derive some upper bounds on \( f(\xv + \uv) \) for \( \uv \) large using that \( f \) is concave.
More precisely, we fix \( \xv \) and \( \uv \) and bound the values
\(f(\xv + t \uv) - f(\xv) - t f'(\xv,\uv)\) for \(\uv \in \UV\) and large \(t\).

\begin{lemma}
\label{LffpfppEL}
Suppose \eqref{1mffmxum34} with \(\delta_{m} \leq 1\) for \(m=3,4\).
Let \(\xv + \UV \subset \Xs\). 
Let the function \(f(\xv + t \uv)\) be concave in \(t\). 
Then it holds for any \(\uv \in \UV\) and for \(t > 1\)
\begin{EQA}
    f(\xv + t \uv) - f(\xv) - \bigl\langle \nabla f(\xv), \uv \bigr\rangle \, t
    & \leq &
    \Bigl( t - \frac{1}{2} \Bigr) 
    \Bigl\{ \langle \nabla^{2} f(\xv) \uv, \uv \rangle - 3 \delta_{3} \Bigr\} .
\label{tm12g1g026}
\end{EQA}
\end{lemma}

\begin{proof}
The  Taylor expansion of the third order for \(g(t) = f(\xv + t\uv)\) at \(t=0\) yields
\begin{EQA}
    \biggl| g(1) - g(0) - g'(0) - \frac{1}{2} g''(0) \biggr|
    & \leq &
    \delta_{3} .
\label{d3Ad3g10}
\end{EQA}
Similarly one obtains
\begin{EQA}
    g'(1) - g'(0)
    &=&
    g'(1) - g'(0) - g''(0) + g''(0)
    \leq 
    g''(0) + 3 \delta_{3} \, .
\label{gpp063ut}
\end{EQA}
Concavity of \(g(\cdot)\) implies
\begin{EQA}
    g(t) - g(1)
    & \leq &
    (t-1) g'(1) .
\label{gtg1tm1gp1}
\end{EQA}
We summarize that
\begin{EQA}
    g(t) - g(0) - t g'(0)
    & = &
    g(t) - g(1) - (t-1) g'(1) + (t-1) \bigl\{ g'(1) - g'(0) \bigr\} + g(1) - g(0) - g'(0)
    \\
    & \leq &
    (t-1) \bigl\{ g''(0) + 3 \delta_{3} \bigr\}
    + \frac{1}{2} g''(0) + \delta_{3}
    \\
    & \leq &
    (t-1/2) \bigl\{ g''(0) + 3 \delta_{3} \bigr\} .
\label{t12gpp6d3uvt}
\end{EQA}
This implies the assertion in view of \(g''(0) = \langle \nabla^{2} f(\xv) \uv, \uv \rangle\).
\end{proof}

Now we specify the result of Lemma~\ref{LffpfppEL} for the elliptic set 
\(\UV(\rups)\) defined by the condition \(- \langle \nabla^{2} f(\xv) \uv, \uv \rangle \leq \rups^{2}\).
We write \(\delta_{3}(\rups)\) in place of 
\(\delta_{3}(\Xs,\UV(\rups))\).
We aim at bounding from above 
the value \(f(\xv + \uv) - f(\xv)\) for
\(\uv\) with \(- \langle \nabla^{2} f(\xv) \uv, \uv \rangle = \rr^{2} > \rups^{2}\).

\begin{lemma}
\label{CELLDtuvLt}
Consider \(\xv \in \Xs\) and 
\(\UV = \UV(\rups) = \bigl\{ \uv \colon - \langle \nabla^{2} f(\xv) \uv, \uv \rangle \leq \rups^{2} \bigr\}\).
Let \(f(\xv + \uv)\) be concave in \(\uv\).
Then for any \(\uv\) with \(- \langle \nabla^{2} f(\xv) \uv, \uv \rangle = \rr^{2} > \rups^{2}\), it holds
\begin{EQA}
  f(\xv + \uv) - f(\xv) - \bigl\langle \nabla f(\xv), \uv \bigr\rangle
  & \leq &
  - (\rr \rups - \rups^{2}/2) 
  \bigl\{ 1 - 3 \rups^{-2} \delta_{3}(\rups) \bigr\} .
\label{LGtuLGtr221622}
\end{EQA}
\end{lemma}

\begin{proof}
Define \(t = \rr/\rups\) and \(\uvd = \uv \rups/\rr\), so that 
\( - \langle \nabla^{2} f(\xv) \uv, \uv \rangle = \rups^{2} \) and \(\uvd \in \UV(\rups)\).
Then it holds by \eqref{tm12g1g026}
\begin{EQA} 
    f(\xv + \uv) - f(\xv) - \bigl\langle \nabla f(\xv), \uv \bigr\rangle
    &=&
    f(\xv + t \uvd) - f(\xv) - \bigl\langle \nabla f(\xv), \uv \bigr\rangle
    \\
    & \leq &
    - (\rr/\rups - 1/2) \bigl\{ \rups^{2} - 3 \delta_{3}(\rups) \bigr\}
    \\
    &=&
    - (\rr \rups - \rups^{2}/2) 
    \bigl\{ 1 - 3 \rups^{-2} \delta_{3}(\rups) \bigr\} 
\label{rr0mr0216r2}
\end{EQA}
and the result follows.
\end{proof}

The result is meaningful if \(3 \rups^{-2} \delta_{3}(\rups) < 1\).
Then with \(\CONSTru = 1 - 3 \rups^{-2} \delta_{3}(\rups)\), we obtain
for any \(\uv\) with \(- \langle \nabla^{2} f(\xv) \uv, \uv \rangle = \rr^{2} > \rups^{2}\)
\begin{EQA}
    f(\xv + \uv) - f(\xv) - \bigl\langle \nabla f(\xv), \uv \bigr\rangle
    & \leq &
    - \CONSTru (\rr \rups - \rups^{2}/2) .
\label{LGtuLGuuT22r22}
\end{EQA}

\subsection{Gaussian integrals}
Let \( \TAU \) be a linear operator in \(\R^{\dimp}\), \( \dimp \leq \infty \), 
with \(\| \TAU \|_{\oper} \leq 1\).
By \( \TAU^{\T} \) we denote the adjoint operator for \( \TAU \).
Given positive \(\rups\) and \(\CONSTru\), consider the following ratio 
\begin{EQA}[c]
    \frac{\int_{\| \TAU \uv \| > \rups}
    \exp\bigl(
      - \CONSTru \| \TAU \uv \| + \frac{1}{2} \CONSTru \rups^{2}
      + \frac{1}{2} \| \TAU \uv \|^{2} 
      - \frac{1}{2} \| \uv \|^{2} 
    \bigr) d\uv}
    {\int_{\| \TAU \uv \| \leq \rups}
    \exp\bigl( - \frac{1}{2}\| \uv \|^{2} \bigr) d\uv} \, .
\label{AAuvrupsc}
\end{EQA}
Obviously, one can rewrite this value as ratio of two expectations
\begin{EQA}[c]
	\frac{\E \Bigl\{ \exp\bigl(	- \CONSTru \rups \| \TAU \gaussv \| 
			+ \frac{1}{2} \CONSTru \rups^{2}
  			+ \frac{1}{2} \| \TAU \gaussv \|^{2} \bigr) 
			\Ind\bigl( \| \TAU \gaussv \| > \rups \bigr) \Bigr\}}
		 {\P\bigl( \| \TAU \gaussv \| \leq \rups \bigr) } \, ,
\label{EindAAgelerr}
\end{EQA}
where \( \gaussv \sim \ND(0,\Id_{\dimp}) \).
Note that without the linear term \(- \CONSTru \| \TAU \gaussv \|\) in the exponent, 
the expectation in the numerator can be infinite. 
We aim at describing \(\rups\) and \(\CONSTru\)-values which ensure that 
the probability in denominator is close to one while the expectation in the numerator is small.

\begin{lemma}
\label{TGaussintext}
Let \(\TAU\) be a linear operator in \(\R^{\dimp}\) with 
\(\| \TAU \|_{\oper} \leq 1\).
Define \(\dimAA = \tr(\TAU^{\T} \TAU)\).
For any \(\CONSTru, \rups\) with 
\(1/2 < \CONSTru \leq 1\) and 
\(\CONSTru \rups = 2 \sqrt{\dimAA} + \sqrt{\xx} \) for \( \xx > 0 \) 
\begin{EQA}
  	\E \Bigl\{ \exp\Bigl(
    		- \CONSTru \rups \| \TAU \gaussv \| + \frac{\CONSTru \rups^{2}}{2}
    		+ \frac{1}{2} \| \TAU \gaussv \|^{2} 
  	\Bigr) \Ind\bigl( \| \TAU \gaussv \| > \rups \bigr) \Bigr\}
  	& \leq &
	\CONST \ex^{ - (\dimAA + \xx)/2}
\label{EexCruAg12Ag2}
\end{EQA}
and 
\begin{EQA}
    \P\bigl( \| \TAU \gaussv \| \leq \rups \bigr)
    & \geq &
    1 - \exp \Bigl\{ - \frac{1}{2} (\rups - \sqrt{\dimAA})^{2} \Bigr\}
    \geq 
    1 - \ex^{- (\dimAA + \xx)/2}.
\label{EeAg12Ag2}
\end{EQA}
\end{lemma}

\begin{remark}
The result applies even if the full dimension \(\dimp\) is infinite and \(\gaussv\) is a Gaussian element in a Hilbert space, provided that \(\dimAA = \tr(\TAU^{\T} \TAU)\) is finite, that is, \(\TAU^{\T} \TAU\) is a trace operator.
\end{remark}

\begin{proof}
Define
\begin{EQA}
	\PAAr(\rr)
	& \eqdef &
	\P\bigl( \| \TAU \gaussv \| \geq \rr \bigr),
	\\
	f(\rr)
	& \eqdef &
	\exp \Bigl( 
		- \CONSTru \rups \rr + \frac{\CONSTru \rups^{2}}{2} + \frac{\rr^{2}}{2} 
	\Bigr) .
\label{PhirPAgrfrC0}
\end{EQA}
Then
\begin{EQA}
	&& \nquad
	\E \Bigl\{ \exp\Bigl(
        - \CONSTru \rups \| \TAU \gaussv \| + \frac{\CONSTru \rups^{2}}{2}
        + \frac{1}{2} \| \TAU \gaussv \|^{2} 
      \Bigr) \Ind\bigl( \| \TAU \gaussv \| > \rups \bigr) \Bigr\}
    \\
    &=&
    - \int_{\rups}^{\infty} f(\rr) d\PAAr(\rr)
    =
    f(\rups) \PAAr(\rups)
    + \int_{\rups}^{\infty} f'(\rr) \PAAr(\rr) \, d\rr \, .
\label{irifpPhrdrf}
\end{EQA}
Now we use that
\( \PAAr\bigl( \sqrt{\dimAA} + \sqrt{2\xx} \bigr) \leq \ex^{-\xx} \)
for any \( \xx > 0 \).
This can be rewritten as
\begin{EQA}
	\PAAr(\rr)
	& \leq &
	\exp \Bigl\{ - \frac{1}{2} \bigl( \rr - \sqrt{\dimAA} \bigr)^{2} \Bigr\}
\label{Pr12rsA2mr}
\end{EQA}
for \( \rr > \sqrt{\dimAA} \).
In particular, in view of \( \rups \geq 2 \sqrt{\dimAA} + \sqrt{\xx} \)
\begin{EQA}
	f(\rups) \PAAr(\rups)
	& \leq &
	\PAAr(\rups)
	\leq 
	\exp \Bigl\{ - \frac{1}{2} \bigl( \rups - \sqrt{\dimAA} \bigr)^{2} \Bigr\}
	\leq 
	\exp\bigl\{ - \frac{1}{2} \bigl( \dimAA + \xx \bigr) \bigr\} .
\label{?}
\end{EQA}
Now we use that \( f'(\rr) = (\rr - \CONSTru \rups) f(\rr) \) and 
\begin{EQA}
	\int_{\rups}^{\infty} f'(\rr) \PAAr(\rr) \, d\rr
	&=&
	\int_{\rups}^{\infty} (\rr - \CONSTru \rups) f(\rr) \PAAr(\rr) \, d\rr
	\\
	& \leq &
	\int_{\rups}^{\infty} (\rr - \CONSTru \rups) 
		\exp\Bigl\{ - \frac{1}{2} \bigl( \rr - \sqrt{\dimAA} \bigr)^{2} 
			- \CONSTru \rups \rr + \frac{\CONSTru \rups^{2}}{2} + \frac{\rr^{2}}{2}
		\Bigr\} \, d\rr
	\\
	&=&
	\int_{\rups}^{\infty} (\rr - \CONSTru \rups) 
		\exp\Bigl\{ 
			- \bigl( \CONSTru \rups - \sqrt{\dimAA} \bigr) \rr 
			+ \frac{\CONSTru \rups^{2}}{2} - \frac{\dimAA}{2} 
		\Bigr\} \, d\rr
	\\
	&=&
	\int_{0}^{\infty} (x + \rups - \CONSTru \rups) 
		\exp\Bigl\{ 
			- \bigl( \CONSTru \rups - \sqrt{\dimAA} \bigr) (x + \rups) 
			+ \frac{\CONSTru \rups^{2}}{2} - \frac{\dimAA}{2} 
		\Bigr\} \, dx .
\label{i0iex0ixexd}
\end{EQA}
The use of \( \int_{0}^{\infty} \ex^{-x} dx = \int_{0}^{\infty} x \ex^{-x} dx = 1 \)
yields
\begin{EQA}
	\int_{\rups}^{\infty} f'(\rr) \PAAr(\rr) \, d\rr
	& \leq &
	\left( 
		\frac{\rups - \CONSTru \rups}{\CONSTru \rups - \sqrt{\dimAA}} 
		+ \frac{1}{(\CONSTru \rups - \sqrt{\dimAA})^{2}} 
	\right)
	\exp\Bigl\{ \rups \sqrt{\dimAA} 
			- \frac{\CONSTru \rups^{2}}{2} - \frac{\dimAA}{2} \Bigr\} .
\label{irifprPr222}
\end{EQA}
It remains to check that for \( \CONSTru \in (1/2,1) \)
and \( \CONSTru \rups = 2\sqrt{\dimAA} + \sqrt{\xx} \)
\begin{EQA}
	- \rups \sqrt{\dimAA} 
			+ \frac{\CONSTru \rups^{2}}{2} + \frac{\dimAA}{2}
	& \geq &
	\frac{\xx + \dimAA}{2} \, .
\label{xp2r0sQfr2}
\end{EQA}
The result follows.
\end{proof}

Now we consider Gaussian integrals with an additional quadratic multiplier.
\begin{lemma}
\label{TGaussintextquad}
Let \( \TAU \) be a linear operator in \(\R^{\dimp}\) with \(\| \TAU \|_{\oper} \leq 1\).
Let \(\zev \in \R^{\infty}\) be a unit norm vector: \(\|\zev\| = 1\).
Define \(\dimAA = \tr(\TAU^{\T} \TAU)\).
For any positive \(\CONSTru, \rups\) with 
\(1/2 < \CONSTru \leq 1\) and 
\(\CONSTru \rups > 2 \sqrt{\dimAA+1} + \sqrt{\xx}\) 
\begin{EQA}
  	&& \nquad
  	\E \Bigl\{|\langle \zev, \gaussv \rangle|^2 \exp\Bigl(
    - \CONSTru \rups \| \TAU \gaussv \| + \frac{\CONSTru \rups^{2}}{2}
    + \frac{1}{2} \| \TAU \gaussv \|^{2} 
  	\Bigr) \Ind\bigl( \| \TAU \gaussv \| > \rups \bigr) \Bigr\} 
   	\\
	& \leq &
  	\CONST \ex^{ - (\dimAA + \xx)/2}.
\label{lgr0TgCr02Se2}
\end{EQA}
\end{lemma}
  
\begin{proof}
Define \( \AAla \) by \( \AAla^{\T} \AAla = \TAU^{\T} \TAU + \zev \otimes \zev \).
Obviously \( \| \AAla \gaussv \| \geq \| \TAU \gaussv \| \),
\( |\langle \zev, \gaussv \rangle| \leq \| \AAla \gaussv \| \).
Further, \( \rr^{2}/2 - \CONSTru \rups \rr \) grows in \( \rr \geq \rups \)
in view of \( \CONSTru \leq 1 \).
Therefore,
\begin{EQA}
	&& \nquad
	|\langle \zev, \gaussv \rangle|^{2} \exp\Bigl(
        - \CONSTru \rups \| \TAU \gaussv \| + \frac{\CONSTru \rups^{2}}{2}
        + \frac{1}{2} \| \TAU \gaussv \|^{2} 
      \Bigr) \Ind\bigl( \| \TAU \gaussv \| > \rups \bigr)  
    \\
    & \leq &
    \| \AAla \gaussv \|^{2} \exp\Bigl(
        - \CONSTru \rups \| \AAla \gaussv \| + \frac{\CONSTru \rups^{2}}{2}
        + \frac{1}{2} \| \AAla \gaussv \|^{2} 
      \Bigr) \Ind\bigl( \| \AAla \gaussv \| > \rups \bigr)   
\label{lzg2eCTgrC22}
\end{EQA}
Now we can follow the line of the proof of Lemma~\ref{TGaussintext}.
Consider
\begin{EQA}
	\PAArg(\rr)
	& = &
	\P\bigl( \| \AAla \gaussv \| \geq \rr \bigr)
	\leq 
	\exp\bigl\{ - \frac{1}{2} (\rr - \sqrt{\dimAla}) \bigr\},
	\\
	f(\rr)
	& \eqdef &
	\rr^{2} \exp \Bigl( 
		- \CONSTru \rups \rr + \frac{\CONSTru \rups^{2}}{2} + \frac{\rr^{2}}{2} 
	\Bigr) 
\label{PhirElagAgrfrC0r2}
\end{EQA}
with \( \dimAla \eqdef \tr \AAla^{\T} \AAla = \dimA + 1 \).
Then
\begin{EQA}
	&& \nquad
	\E \Bigl\{ \bigl( \langle \zev, \gaussv \rangle \bigr)^{2} \exp\Bigl(
        - \CONSTru \rups \| \TAU \gaussv \| + \frac{\CONSTru \rups^{2}}{2}
        + \frac{1}{2} \| \TAU \gaussv \|^{2}  
      \Bigr) \Ind\bigl( \| \TAU \gaussv \| > \rups \bigr) \Bigr\}
    \\
    & \leq &
    - \int_{\rups}^{\infty} f(\rr) d\PAArg(\rr)
    =
    f(\rups) \PAArg(\rups)
    + \int_{\rups}^{\infty} f'(\rr) \PAArg(\rr) \, d\rr \, .
\label{irifplaPhrdrf}
\end{EQA}
Now we can continue as in the proof of Lemma~\ref{TGaussintext}.
\end{proof}

The bound \eqref{lgr0TgCr02Se2} can be easily extended to the case of a more general functional 
\( \QP \gammav \) in place of \( \langle \lambda, \gammav \rangle \).

\begin{lemma}
\label{TGaussintextquadQ}
Let \( \TAU \) be a linear operator in \(\R^{\dimp}\) with \(\| \TAU \|_{\oper} \leq 1\)
and \(\dimAA = \tr(\TAU^{\T} \TAU) < \infty \).
Let \( \AP \) be a bounded linear operator with  
\( \tr (\AP^{\T} \AP) < \infty \).
For any positive \(\CONSTru, \rups\) with 
\(1/2 < \CONSTru \leq 1\) and 
\(\CONSTru \rups > 2 \sqrt{\dimAA+1} + \sqrt{\xx}\) 
\begin{EQA}
  	\E \Bigl\{\| \AP \gaussv \|^2 \exp\Bigl(
    - \CONSTru \rups \| \TAU \gaussv \| + \frac{\CONSTru \rups^{2}}{2}
    + \frac{1}{2} \| \TAU \gaussv \|^{2} 
  	\Bigr) \Ind\bigl( \| \TAU \gaussv \| > \rups \bigr) \Bigr\} 
   	& \lesssim &
  	\tr(\AP^{\T} \AP) \ex^{ - (\dimAA + \xx)/2}.
\label{lgr0TgCr02Se2Q}
\end{EQA}
\end{lemma}

\begin{proof}
We use the Karhunen-Loeve decomposition of \( \AP^{\T} \AP \):
\begin{EQA}
	\| \AP \gammav \|^{2}
	&=&
	\sum_{j} \mu_{j} \langle \zev_{j}, \gammav \rangle^{2}
\label{Agsjmjzjg2}
\end{EQA}
with orthogonal unit vectors \( \zev_{j} \) and \( \sum_{j} \mu_{j} = \tr(\AP^{\T}\AP) \),
and apply the result of Lemma~\ref{TGaussintextquad} to each term of this decomposition. 
\end{proof}

\section{Proofs of the main results}
This section collects the proofs of the main theorems.

\paragraph{Proof of Proposition~\ref{TconstGrDG}}
The idea of the proof is to show that for each \( \uv \) with \( \| \DPGP \uv \| = \rr_{\GP} \), 
the derivative of the function \( \LGP(\thetavs_{\GP} + t \uv) \) in \( t \) is negative for 
\( |t| \geq 1 \). 
This yields that the point of maximum of \( \LGP(\thetav) \) cannot be outside of 
\( \CA_{\GP}(\rr_{\GP}) \).
Let us fix any \( \uv \) with \( \| \DPGP \uv \| \leq \rr_{\GP} \).
We use the decomposition
\begin{EQA}
	\LGP(\thetavs_{\GP} + t \uv) - \LGP(\thetavs_{\GP})
	&=&
	\bigl\langle \nabla \zeta, \uv \bigr\rangle \, t + \E \LGP(\thetavs_{\GP} + t \uv) - \E \LGP(\thetavs_{\GP}) .
\label{LGtsGtuLGts}
\end{EQA}
With \( f(t) = \E \LGP(\thetavs_{\GP} + t \uv) \), it holds
\begin{EQA}
	\frac{d}{dt} \LGP(\thetavs_{\GP} + t \uv)
	&=&
	\bigl\langle \nabla \zeta, \uv \bigr\rangle + f'(t) .
\label{frddtLtGstu}
\end{EQA}
The bound \eqref{uTzDGvTxm12z} implies on \( \Omega(\xx) \)
\begin{EQA}
	\bigl| \bigl\langle \nabla \zeta, \uv \bigr\rangle \bigr|
	&=&
	\bigl| \langle \DPGP^{-1} \nabla \zeta, \DPGP \uv \rangle \bigr|
	\leq 
	\rr_{\GP} \, \zq(\BB_{\VPGP},\xx) .
\label{uTnzrzBGx}
\end{EQA}
By definition of \( \thetavs_{\GP} \), it also holds \( f'(0) = 0 \).
Condition \nameref{LL0ref} implies 
\begin{EQA}
	\bigl| f'(t) - t f''(0) \bigr|
	&=&
	\bigl| f'(t) - f'(0) - t f''(0) \bigr|
	\leq 
	3 t^{2} \, \rr_{\GP}^{3} \, \dltwu_{3,\HPc} .
\label{fptfp0fpttfpp13}
\end{EQA}
For \( t = 1 \), we obtain
\begin{EQA}
	f'(1) 
	& \leq &
	f''(0) + 3 \rr_{\GP}^{3} \, \dltwu_{3,\HPc}
	=
	- \langle \DPGP^{2} \uv, \uv \rangle + 3 \rr_{\GP}^{3} \, \dltwu_{3,\HPc}
	=
	- \rr_{\GP}^{2} + 3 \rr_{\GP}^{3} \, \dltwu_{3,\HPc} .
\label{fp1fpp13d3rG}
\end{EQA}
If \( 3 \dltwu_{3,\HPc} \leq \rho \) for \( \rho < 1 \), 
then \( f'(1) < 0 \).
Concavity of \( f(t) \) and \( f'(0) = 0 \) imply that \( f'(t) \) decreases in 
\( t \) for \( t > 1 \).
Further, on \( \Omega(\xx) \) by \eqref{uTnzrzBGx}
\begin{EQA}
	\frac{d}{dt} \LGP(\thetavs_{\GP} + t \uv) \cond_{t=1}
	& \leq &
	\bigl\langle \nabla \zeta, \uv \bigr\rangle - \rr_{\GP}^{2} + 3 \rr_{\GP}^{3} \, \dltwu_{3,\HPc}
	\\
	& \leq &
	\rr_{\GP} \, \zq(\BB_{\GP},\xx) - \rr_{\GP}^{2} + 3 \rr_{\GP}^{3} \, \dltwu_{3,\HPc}
	\leq 
	\rr_{\GP} \, \zq(\BB_{\GP},\xx) - (1 - \rho) \rr_{\GP}^{2}
	< 0
\label{ddtLGtstu33}
\end{EQA}
for \( \rr_{\GP} > (1 - \rho)^{-1} \zq(\BB_{\GP},\xx) \).
As \( \frac{d}{dt} \LGP(\thetavs_{\GP} + t \uv) \) decreases with \( t \geq 1 \) together with \( f'(t) \) due to \eqref{frddtLtGstu}, the same applies to all such \( t \).
This implies the assertion.

\paragraph{Proof of Theorem~\ref{TFiWititG}}
To show \eqref{3d3Af12DGttt}, we use that \( \tilde{\thetav}_{\GP} \in \CA_{\GP}(\rr_{\GP}) \) and
\( \nabla \LGP(\tilde{\thetav}_{\GP}) = 0 \).
Therefore,
\begin{EQA}
	\LGP(\tilde{\thetav}_{\GP} + \uv) - \LGP(\tilde{\thetav}_{\GP})
	&=&
	\LGP(\tilde{\thetav}_{\GP} + \uv) - \LGP(\tilde{\thetav}_{\GP})
	 - \langle \nabla \LGP(\tilde{\thetav}_{\GP}), \uv \rangle .
\label{LGtttGLGtTtGW}
\end{EQA}
Let us fix any \( \thetav \in \CA_{\GP}(\rr_{\GP}) \) and \( \uv \) with 
\( \| \DPGP \uv \| \leq \rr_{\GP} \), and consider
\begin{EQA}
	f(t)
	=
	f(t,\uv)
	& \eqdef &
	\LGP(\thetav + t \uv) - \LGP(\thetav) - \langle \nabla \LGP(\thetav), \uv \rangle \, t .
\label{fuELtuELGW}
\end{EQA}
As the stochastic term of \( L(\thetav) \) and thus, of \( \LGP(\thetav) \) is linear in \( \thetav \), it cancels in this expression,
and it suffices to consider the deterministic part \( \E \LGP(\thetav) \). 
Obviously \( f(0) = 0 \), \( f'(0) = 0 \).
Moreover, \( f''(0) = \langle \nabla^{2} \E \LGP(\thetav) \, \uv, \uv \rangle 
= - \langle \DPGP^{2}(\thetav) \, \uv, \uv \rangle < 0 \).
Taylor expansion of the third order implies 
\begin{EQA}
	\bigl| f(1) - \frac{1}{2} f''(0) \bigr|  
	& \leq &
	\bigl| \dltwu_{3,\HPc}(\thetavc,\uv) \bigr| \, , 
	\quad
	\thetavc \in [\thetav,\thetav + \uv] .
\label{fu12un2fud2W}
\end{EQA}
In particular, for any \( \thetav \in \CA_{\GP}(\rr_{\GP}) \) 
\begin{EQA}
	\Bigl| 
		\E \LGP(\thetavs_{\GP}) - \E \LGP(\thetav) - \frac{1}{2} \bigl\| \DPGP (\thetav - \thetavs_{\GP}) \bigr\|^{2} 
	\Bigr|
	& \leq &
	\rr_{\GP}^{3} \, \dltwu_{3,\HPc} .
\label{d3GrGELGtsG12}
\end{EQA}
We now use that \( \nabla \LGP(\tilde{\thetav}_{\GP}) = 0 \) and
by Proposition~\ref{TconstGrDG}, \( \uv = \thetavs_{\GP} - \tilde{\thetav}_{\GP} \) fulfills 
\( \| \DPGP \uv \| \leq \rr_{\GP} \) on \( \Omega(\xx) \).
Therefore, for \( \thetav \in \CA_{\GP}(\rr_{\GP}) \)
\begin{EQA}
	&& \nquad
	\Bigl| 
		\LGP(\thetav) - \LGP(\tilde{\thetav}_{\GP}) 
		- \frac{1}{2} \| \DPGPt (\thetav - \tilde{\thetav}_{\GP}) \|^{2}
	\Bigr|
	\\
	&=&
	\Bigl| 
		\LGP(\thetav) - \LGP(\tilde{\thetav}_{\GP}) 
		- \langle \nabla \LGP(\tilde{\thetav}_{\GP}), \thetav - \tilde{\thetav}_{\GP} \rangle 
		- \frac{1}{2} \| \DPGPt (\thetav - \tilde{\thetav}_{\GP}) \|^{2}
	\Bigr|
	\leq 
	\rr_{\GP}^{3} \, \dltwu_{3,\HPc} .
\label{23GrG122LGt}
\end{EQA}
The result \eqref{3d3Af12DGttt} follows. 
Further, as \( \tilde{\thetav}_{\GP} \in \CA_{\GP}(\rr_{\GP}) \), it holds
\begin{EQA}
	&& \nquad
	\LGP(\tilde{\thetav}_{\GP}) - \LGP(\thetavs_{\GP}) - \frac{1}{2} \| \DPGP^{-1} \nabla \zeta \|^{2}
	=
	\max_{\thetav \in \CA_{\GP}(\rr_{\GP})} 
	\Bigl\{ 
		\LGP(\thetav) - \LGP(\thetavs_{\GP}) - \frac{1}{2} \| \DPGP^{-1} \nabla \zeta \|^{2} 
	\Bigr\}
	\\
	&=&
	\max_{\thetav \in \CA_{\GP}(\rr_{\GP})} 
	\Bigl\{ \bigl\langle \thetav - \thetavs_{\GP}, \nabla \zeta \bigr\rangle
	+ \E \LGP(\thetav) - \E \LGP(\thetavs_{\GP}) - \frac{1}{2} \| \DPGP^{-1} \nabla \zeta \|^{2} 
	\Bigr\}
	\\
	& \leq &
	\max_{\thetav \in \CA_{\GP}(\rr_{\GP})} 
	\Bigl\{ \bigl\langle \DPGP (\thetav - \thetavs_{\GP}), \DPGP^{-1} \nabla \zeta \bigr\rangle
		- \frac{1}{2} \| \DPGP (\thetav - \thetavs_{\GP}) \|^{2} 
		- \frac{1}{2} \| \DPGP^{-1} \nabla \zeta \|^{2} 
	\Bigr\} + \rr_{\GP}^{3} \, \dltwu_{3,\HPc}
	\\
	& \leq &
	\max_{\thetav \in \CA_{\GP}(\rr_{\GP})} 
	\Bigl\{ 
		- \frac{1}{2} \| \DPGP (\thetav - \thetavs_{\GP}) - \DPGP^{-1} \nabla \zeta \|^{2} 
	\Bigr\} + \rr_{\GP}^{3} \, \dltwu_{3,\HPc}
	=  
	\rr_{\GP}^{3} \, \dltwu_{3,\HPc}
\label{d3G1212222B} 
\end{EQA}
and similarly 
\( \LGP(\tilde{\thetav}_{\GP}) - \LGP(\thetavs_{\GP}) - \frac{1}{2} \| \DPGP^{-1} \nabla \zeta \|^{2} \geq - \rr_{\GP}^{3} \, \dltwu_{3,\HPc} \).
This two-sided bound yields as \eqref{DGttGtsGDGm13rG} as \eqref{3d3Af12DGttG}.

The last statement \eqref{DPGPm1Cd3rG} of the theorem follows directly from Lemma~\ref{LellUVD2w}
with \( \QP = \DPGP \) and \( f(\thetav) = \E \LGP(\thetav) \).

\paragraph{Proof of Proposition~\ref{TbiasGP} and Theorem~\ref{TestlosspMLE}}
The definition of \( \thetavs \) and \( \thetavs_{\GP} \) implies 
\begin{EQA}
	\E \LGP(\thetavs_{\GP})
	\geq 
	\E \LGP(\thetavs),
	&\quad &
	\E L(\thetavs_{\GP})
	\leq 
	\E L(\thetavs) .
\label{ELGttsGELELs}
\end{EQA}
As \( \E \LGP(\thetav) = \E L(\thetav) - \| \GP \thetav \|^{2}/2 \), it follows that
\begin{EQA}
	2 \E \LGP(\thetavs_{\GP}) - 2 \E \LGP(\thetavs)
	& \leq & 
	\bigl\| \GP \thetavs \bigr\|^{2}
	- \bigl\| \GP \thetavs_{\GP} \bigr\|^{2}
	\leq 
	\bigl\| \GP \thetavs \bigr\|^{2} .
\label{12ELGSmGts22}
\end{EQA}
The bound \eqref{d3GrGELGtsG12} with \( \thetav = \thetavs \) implies the first statement of \eqref{fr33GrGd3rELG22}. 

Further we show that \( \| \GP \thetavs \| \leq \rrbias/2 \) implies \( \| \DPGP (\thetavs_{\GP} - \thetavs) \| \leq  \rrbias \).
Indeed, suppose the opposite inequality.
Define \( \uv = \rrbias \DPGP (\thetavs - \thetavs_{\GP}) / \| \DPGP (\thetavs_{\GP} - \thetavs) \| \), so that 
\( \| \uv \| = \rrbias \). 
The function \( f(t) = \E \LGP(\thetavs_{\GP}) - \E \LGP(\thetavs_{\GP} + t \uv) \) is convex in \( t \) and 
\( \thetavs_{\GP} + t \uv \in \Thetad \) for \( |t| \leq 1 \).
Using the approximation \eqref{d3GrGELGtsG12} for \( \thetav = \thetavs_{\GP} + \uv \) implies 
\begin{EQA}
	2 \E \LGP(\thetavs_{\GP}) - 2 \E \LGP(\thetavs_{\GP} + t \uv)
	& \geq & 
	\rrbias^{2} - \rrbias^{3} \, \dltwu_{3,\HPc}(\rrbias) 
	\geq 
	\rrbias^{2}/2
\label{tts24ELGEL}
\end{EQA}
and concavity of \( \E \LGP(\thetav) \) together with \( \nabla \E \LGP(\thetavs_{\GP}) = 0 \) 
implies for \( t \geq 1 \)
\begin{EQA}
	\E \LGP(\thetavs_{\GP}) - \E \LGP(\thetavs_{\GP} + t \uv)
	& \geq &
	\rrbias^{2}/2 .
\label{24t1ELGELGtu}
\end{EQA}
This contradicts to the bounds \eqref{12ELGSmGts22} and \( \| \GP \thetavs \|^{2} \leq \rrbias^{2}/2 \).

Now for any \( \thetav \) with \( \| \DPGP (\thetavs_{\GP} - \thetav) \| \leq \rrbias \)
\begin{EQA}
	\Bigl| 
		\E \LGP(\thetavs_{\GP}) - \E \LGP(\thetav) - \frac{1}{2} \bigl\| \DPGP (\thetav - \thetavs_{\GP}) \bigr\|^{2}
	\Bigr|
	& \leq &
	\rrbias^{3} \, \dltwu_{3,\HPc}(\rrbias) .
\label{BB3GrGLGD2}
\end{EQA}
Further we use that \( \thetavs = \argmax \E L(\thetav) \) and 
\( \E \LGP(\thetav) = \E L(\thetav) - \| \GP \thetav \|^{2}/2  \).
By \eqref{BB3GrGLGD2} in view of \( \| \DPGP (\thetavs_{\GP} - \thetavs) \| \leq  \rrbias \)
and \( \DPGP^{2} = \IF(\thetavs_{\GP}) + \GP^{2} = \DP^{2} + \GP^{2} \)
\begin{EQA}
	\E L(\thetavs) - \E \LGP(\thetavs_{\GP})
	&=&
	\max_{\thetav \in \CA_{\GP}(\rrbias)} 
	\bigl\{ \E \LGP(\thetav) + \frac{1}{2} \| \GP \thetav \|^{2} 
		- \E \LGP(\thetavs_{\GP}) \bigr\}
	\\
	& \leq & 
	\max_{\thetav \in \CA_{\GP}(\rrbias)} 
	\Bigl\{ 
		- \frac{1}{2} \bigl\| \DPGP (\thetav - \thetavs_{\GP}) \bigr\|^{2} + \frac{1}{2} \| \GP \thetav \|^{2} 
	\Bigr\} + \rrbias^{3} \, \dltwu_{3,\HPc}(\rrbias)
	\\
	& = & 
	\max_{\thetav \in \CA_{\GP}(\rrbias)} 
	\Bigl\{ 
		- \frac{1}{2} \bigl\| \DP \thetav - \DP^{-1} \DPGP^{2} \, \thetavs_{\GP} \bigr\|^{2} 
		+ \frac{1}{2} \| \DP^{-1} \DPGP^{2} \, \thetavs_{\GP} \|^{2} 
	\Bigr\} + \rrbias^{3} \, \dltwu_{3,\HPc}(\rrbias) .
\label{12d34GDm1DG2}
\end{EQA}
A similar inequality holds from below with opposite sign for \( \dltwu_{3,\HPc} \)-term 
yielding for the maximizer \( \thetavs \) the bound
\begin{EQA}
	\bigl\| \DP \thetavs - \DP^{-1} \DPGP^{2} \, \thetavs_{\GP} \bigr\|^{2}
	& \leq &
	4 \rrbias^{3} \, \dltwu_{3,\HPc}(\rrbias) .
\label{43GrbDtsm1}
\end{EQA}
Equivalently, using again \( \DPGP^{2} = \DP^{2} + \GP^{2} \)
\begin{EQA}
	\bigl\| \DP^{-1} \DPGP^{2} (\thetavs - \thetavs_{\GP}) - \DP^{-1} \GP^{2} \thetavs \bigr\|^{2}
	& \leq &
	4 \rrbias^{3} \, \dltwu_{3,\HPc}(\rrbias) .
\label{Dm1D2ttssGm1}
\end{EQA}
As \( \DP^{2} \leq \DPGP^{2} \), this also implies
\begin{EQA}
	\bigl\| \DPGP (\thetavs - \thetavs_{\GP}) - \DPGP^{-1} \GP^{2} \thetavs \bigr\|^{2}
	& \leq &
	4 \rrbias^{3} \, \dltwu_{3,\HPc}(\rrbias) .
\label{43DGttsttsGDGm122}
\end{EQA}
The statement \eqref{22GtsQDGm2loss} of Theorem~\ref{TestlosspMLE} follows from the bound 
\eqref{22GtsQDGm2} on the bias \( \| \QP (\thetavs_{\GP} - \thetavs) \| \) 
and the deviation bound
\( \bigl\| \QP \DP_{\GP}^{-2} \nabla \zeta \bigr\| \leq \zq(\BB_{\QPGP},\xx) \) on a set 
of probability at least \( 1 - \ex^{-\xx} \) 
by the triangle inequality.

For the second statement of Theorem~\ref{TestlosspMLE}, we apply the decomposition
\begin{EQA}
	\bigl\| \QP (\tilde{\thetav}_{\GP} - \thetavs) \bigr\|
	&=&
	\bigl\| \QP (\tilde{\thetav}_{\GP} - \thetavs_{\GP}) - \av \bigr\|
\label{QttGtsattG}
\end{EQA}
with \( \av = \thetavs - \thetavs_{\GP} \).
Now the result follows from asymptotic normality \eqref{szQDGm2nzzVQG} and 
Gaussian comparison result of Theorem~\ref{Tgaussiancomparison3}
in view of small bias condition \eqref{QDGm2G2ts2V2}. 

\paragraph{Proof of Proposition~\ref{PrhoQPBvM}}
Let \( \tilde{\thetav}_{\GP} = \argmax_{\thetav} \LGP(\thetav) \) be the penalized MLE of the parameter \( \thetav \).
We aim at bounding from above the quantity
\begin{EQA}
    \rho(\rups)
    &=&
    \frac{\int_{\| \HPc \uv \| > \rups} \exp\bigl\{ \LGP(\tilde{\thetav}_{\GP}+\uv) \bigr\} d \uv}
    {\int_{\| \HPc \uv \| \leq \rups} \exp \bigl\{ \LGP(\tilde{\thetav}_{\GP}+\uv) \bigr\} d \uv}  
    \, 
\label{rhopiDGP}
\end{EQA}
with 
\( \HPc^{2} \) from \nameref{ED0ref}.
We suppose in the proof that \( \dimp < \infty \).
The general case can be obtained by taking a limit as \( \dimp \to \infty \).
\paragraph{Step 1}

The use of \( \nabla \LGP(\tilde{\thetav}_{\GP}) = 0 \) allows to represent
\begin{EQA}
    \rho(\rups)
    &=&
    \frac{\int_{\| \HPc \uv \| > \rups} \exp \bigl\{ 
    	  \LGP(\tilde{\thetav}_{\GP} + \uv) - \LGP(\tilde{\thetav}_{\GP}) \bigr\} d \uv}
    	 {\int_{\| \HPc \uv \| \leq \rups} \exp \bigl\{ 
	 	  \LGP(\tilde{\thetav}_{\GP} + \uv) - \LGP(\tilde{\thetav}_{\GP}) \bigr\} d \uv}
    \\
    &=&
    \frac{\int_{\| \HPc \uv \| > \rups} \exp \bigl\{ 
    		\LGP(\tilde{\thetav}_{\GP} + \uv) - \LGP(\tilde{\thetav}_{\GP}) 
			- \bigl\langle \nabla \LGP(\tilde{\thetav}_{\GP}), \uv \bigr\rangle
		  \bigr\} d \uv}
    {\int_{\| \HPc \uv \| \leq \rups} \exp \bigl\{ 
    		\LGP(\tilde{\thetav}_{\GP} + \uv) - \LGP(\tilde{\thetav}_{\GP}) 
    		- \bigl\langle \nabla \LGP(\tilde{\thetav}_{\GP}), \uv \bigr\rangle 
		   \bigr\} d \uv}
    \, .
\label{rhopiDGPt}
\end{EQA}
Now we study this expression for any possible value \( \thetav \) from the concentration set of \( \tilde{\thetav}_{\GP} \).
Consider \( f(\thetav) = \E \LGP(\thetav) \).
As the stochastic term of \( L(\thetav) \) and thus, of \( \LGP(\thetav) \) is linear in 
\( \thetav \), it holds
\begin{EQA}
	\LGP(\thetav + \uv) - \LGP(\thetav) - \bigl\langle \nabla \LGP(\thetav), \uv \bigr\rangle 
	&=&
	f(\thetav + \uv) - f(\uv) - \bigl\langle \nabla f(\thetav), \uv \bigr\rangle.
\label{fuELtuELG}
\end{EQA}
Therefore, it suffices to bound uniformly in \( \thetav \in \Thetad \) the ratio
\begin{EQA}
	\rho(\rups,\thetav)
	& \eqdef &
    \frac{\int \Ind\bigl( \| \HPc \uv \| > \rups \bigr) \exp \bigl\{ 
    		f(\thetav + \uv) - f(\uv) - \bigl\langle \nabla f(\thetav), \uv \bigr\rangle \bigr\} d \uv}
    	 {\int \Ind\bigl( \| \HPc \uv \| \leq \rups \bigr)  \exp \bigl\{ 
	 		f(\thetav + \uv) - f(\uv) - \bigl\langle \nabla f(\thetav), \uv \bigr\rangle \bigr\} d \uv}
	\qquad
\label{rhopifUV}
\end{EQA}
%

\paragraph{Step 2}
First we bound the denominator of \( \rho(\rups,\thetav) \).
Lemma~\ref{Lintfxupp2} yields
\begin{EQA}
	\int_{\| \HPc \uv \| \leq \rups} \ex^{ 
		f(\thetav + \uv) - f(\uv) - \langle \nabla f(\thetav), \uv \rangle
	}
	\, d\uv
	& \geq &
	\bigl( 1 - \err(\rups) \bigr)
	\int_{\| \HPc \uv \| \leq \rups} 
	\ex^{ - \| \DPGP(\thetav) \uv \|^{2}/{2}}
	\, d\uv ,
	\quad
	\\
	\int_{\| \HPc \uv \| \leq \rups} 
	\ex^{	f(\thetav + \uv) - f(\uv) - \langle \nabla f(\thetav), \uv \rangle
	}
	\, d\uv
	& \leq &  
	\bigl( 1 + \err(\rups) \bigr)
	\int_{\| \HPc \uv \| \leq \rups} 
	\ex^{ - \| \DPGP(\thetav) \uv \|^{2}/{2}}
	\, d\uv ,
	\quad
\label{1eGAAiDu22dT31}
\end{EQA}
where \( \DPGP^{2}(\thetav) = \IF_{\GP}(\thetav) = - \nabla^{2} f(\thetav) \) and 
\( \err(\rups) \) is given by \eqref{LmgfquadELGP}.
%
Moreover, after a proper normalization, the integral 
\( \int_{\| \HPc \uv \| \leq \rups} \exp \Bigl( - {\| \DPGP(\thetav) \uv \|^{2}}/{2} \Bigr) \, d\uv \)
can be viewed as the probability of the Gaussian event.
Namely
\begin{EQA}
	\frac{\det \DPGP(\thetav)}{(2\pi)^{\dimp/2}}
	\int_{\| \HPc \uv \| \leq \rups} \exp \Bigl( - \frac{\| \DPGP(\thetav) \uv \|^{2}}{2} \Bigr) \, d\uv
	&=&
	\P\bigl( \bigl\| \HPc \DPGP^{-1}(\thetav) \gammav \bigr\| \leq \rups \bigr)
\label{dtDGt2pip2Dtm1r0}
\end{EQA}
for a standard normal \( \gammav \in \R^{\dimp} \). 
The choice \( \rups \geq \sqrt{\dimA_{\GP}(\thetav)} + \sqrt{2\xx} \) yields by Corollary~\ref{Cchi2p}
\begin{EQA}
	\P\bigl( \bigl\| \HPc \DPGP^{-1}(\thetav) \gammav \bigr\| \leq \rups \bigr) 
	& \geq &
	1 - \ex^{-\xx} .
\label{PDtDGm1tgar0}
\end{EQA}
If the error term \( \err(\rups) \) is small, we obtain a sharp 
bound for the integral in the denominator of \( \rho(\rups,\thetav) \) from 
\eqref{rhopifUV}.

\paragraph{Step 3}
Now we bound the integral on the exterior of \( \UVd = \bigl\{ \uv \colon \| \HPc \uv \| \leq \rups \bigr\} \).
Linearity of stochastic term in 
\( \LGP(\thetav) = L(\thetav) - \| \GP \thetav \|^{2}/2 \) 
and quadraticity of the penalty imply
\begin{EQA}
	\LGP(\thetav + \uv) - \LGP(\thetav) - \bigl\langle \nabla \LGP(\thetav), \uv \bigr\rangle
	&=&
	\E L(\thetav + \uv) - \E L(\thetav) - \bigl\langle \nabla \E L(\thetav), \uv \bigr\rangle
	- \frac{1}{2} \| \GP \uv \|^{2} \, .
\label{t22Gu2uTvnLGt}
\end{EQA}
Now we apply Lemma~\ref{CELLDtuvLt} with 
\( f(\thetav + \uv) = \E L(\thetav + \uv) - \bigl\langle (\IF(\thetav) + \HPc^{2}) \uv, \uv \bigr\rangle / 2  \). 
This function is concave and it holds \( - \bigl\langle \nabla^{2} f(\thetav) \uv, \uv \bigr\rangle = \| \HPc \uv \|^{2} \).
The bound \eqref{LGtuLGtr221622} yields for any \( \uv \) with 
\( \| \HPc \uv \| = \rr > \rups \) and
\( \CONSTru = 1 - 3 \rups \, \dltwu_{3,\HPc} \geq 1/2 \)
\begin{EQA}
	&& \nquad
	\LGP(\thetav + \uv) - \LGP(\thetav) - \bigl\langle \nabla \LGP(\thetav), \uv \bigr\rangle 
	\\
	&=&
	f(\thetav + \uv) - f(\thetav) - \bigl\langle \nabla f(\thetav), \uv \bigr\rangle 
	- \bigl\langle (\DPGP^{2} - \HPc^{2}) \uv , \uv \bigr\rangle / 2
	\\
	& \leq &
	- \CONSTru (\| \HPc \uv \| \rups - \rups^{2}/2) 
	- \bigl\langle (\DPGP^{2} - \HPc^{2}) \uv , \uv \bigr\rangle / 2
	\\
	&=&
	- \CONSTru (\| \HPc \uv \| \rups - \rups^{2}/2) - \| \DP_{\GP|\HPc}(\thetav) \uv \|^{2}/2 ,
\label{LGttLGT22r22}
\end{EQA}
where \( \DP_{\GP|\HPc}^{2} = \DPGP^{2} - \HPc^{2} \).
%
Now we can use the result about Gaussian integrals from Section~\ref{SGaussintegr}.
With \( \TAU = \HPc \DP_{\GP|\HPc}^{-1}(\thetav) \), it holds by Lemma~\ref{TGaussintext}
\begin{EQA}
	&& \nquad
	\frac{\det \DPGP(\thetav)}{(2\pi)^{\dimp/2}}
	\int \Ind\bigl( {\| \HPc \uv \| > \rups} \bigr) \exp \bigl\{ 
		\LGP(\thetav + \uv) - \LGP(\thetav) - \bigl\langle \nabla \LGP(\thetav), \uv \bigr\rangle 
	\bigr\} \, d\uv
	\\
	& \leq &
	     \E \Bigl\{ \exp\Bigl(
        - \CONSTru \rups \| \TAU \gaussv \| + \frac{\CONSTru \rups^{2}}{2}
        + \frac{1}{2} \| \TAU \gaussv \|^{2} 
      \Bigr) \Ind\bigl( \| \TAU \gaussv \| > \rups \bigr) \Bigr\}
    \leq 
    \CONST \ex^{ - (\dimA_{\GP}(\thetav) + \xx)/2} .
\label{dDG2pp2LttG2x2}
\end{EQA}
Putting together of Step 1 through Step 3 yields the statement
about \( \rho(\rups) \).

\paragraph{Proof of Theorem~\ref{TnonparBvMm} and Corollary~\ref{CnonparBvMm}} 
We proceed similarly to the proof of Theorem~\ref{PrhoQPBvM}.
Fix any centrally symmetric set \( A \).
First we restrict the posterior probability to the set 
\( \CAt(\rups) = \{ \uv \colon \| \HPc \uv \| \leq \rups \} \).
Then we apply the quadratic approximation of the log-likelihood function \( L(\thetav) \).
Denote \( A(\rups) = A \cap \CAt(\rups) \).
Obviously, \( A(\rups) \) is centrally symmetric as well. 
Further, 
\begin{EQA}
	\P\bigl( \vthetav_{\GP} - \tilde{\thetav}_{\GP} \in A \cond \Yv \bigr)
	&=& 
	\frac{\int_{A} \exp\bigl\{ \LGP(\tilde{\thetav}_{\GP} + \uv) \bigr\} d\uv}
		{\int_{\R^{\dimp}} \exp\bigl\{ \LGP( \tilde{\thetav}_{\GP} + \uv) \bigr\}d\uv}
	\\
	& \leq &
	\frac{\int_{A(\rups)} 
	\exp\bigl\{ \LGP( \tilde{\thetav}_{\GP} + \uv) - \LGP(\tilde{\thetav}_{\GP}) 
				- \bigl\langle \nabla \LGP(\tilde{\thetav}_{\GP}), \uv \bigr\rangle \bigr\} d\uv}
		{\int_{\| \HPc \uv \| \leq \rups} 
	\exp\bigl\{ \LGP( \tilde{\thetav}_{\GP} + \uv) - \LGP(\tilde{\thetav}_{\GP}) 
				- \bigl\langle \nabla \LGP(\tilde{\thetav}_{\GP}), \uv \bigr\rangle \bigr\}d\uv}
	+ \rho(\rups) .
\label{rorLGttGuuG}
\end{EQA}
Now we apply the bounds from the proof of Theorem~\ref{PrhoQPBvM} yielding the upper bound
\begin{EQA}
	\P\bigl( \vthetav_{\GP} - \tilde{\thetav}_{\GP} \in A \cond \Yv \bigr)
	& \leq & 
	\frac{\bigl\{ 1 + \err(\rups) \bigr\} 
		  	\int_{A(\rups)} \exp\bigl\{ - \| \DPGPt \uv \|^{2}/2 \bigr\} d\uv}
		 {\bigl\{ 1 - \err(\rups) \bigr\}
			\int_{\| \HPc \uv \| \leq \rups} \exp\bigl\{ - \| \DPGPt \uv \|^{2}/2 \bigr\}d\uv}	
	+ \rho(\rups) 
	\\
	& \leq &
	\frac{\bigl\{ 1 + \err(\rups) \bigr\} 
			\P\bigl( \DPGPt^{-1} \gammav \in A \bigr)}
		 {\bigl\{ 1 - \err(\rups) \bigr\}
		 \P\bigl( \| \HPc \DPGPt^{-1} \gammav \| \leq \rups \bigr)}	
	+ \rho(\rups) .
\label{1mdi1pdimDGu22}
\end{EQA}
Now we prove the lower bound.
It obviously holds  
\begin{EQA}
	\P\bigl( \vthetav_{\GP} - \tilde{\thetav}_{\GP} \in A \cond \Yv \bigr)
	&=& 
	\frac{\int_{A} \exp\bigl\{ \LGP(\tilde{\thetav}_{\GP} + \uv) \bigr\} d\uv}
		{\int_{\R^{\dimp}} \exp\bigl\{ \LGP( \tilde{\thetav}_{\GP} + \uv) \bigr\}d\uv}
	\\
	& \geq &
	\frac{\int_{A(\rups)} 
			\exp\bigl\{ \LGP( \tilde{\thetav}_{\GP} + \uv) - \LGP(\tilde{\thetav}_{\GP}) 
				- \bigl\langle \nabla \LGP(\tilde{\thetav}_{\GP}), \uv \bigr\rangle \bigr\} d\uv}
		 {\bigl( \int_{\| \HPc \uv \| \leq \rups} + \int_{\| \HPc \uv \| > \rups} \bigr)
	\exp\bigl\{ \LGP( \tilde{\thetav}_{\GP} + \uv) - \LGP(\tilde{\thetav}_{\GP}) 
				- \bigl\langle \nabla \LGP(\tilde{\thetav}_{\GP}), \uv \bigr\rangle \bigr\}d\uv}
\label{rorLGttGuuGl}
\end{EQA}
and in a similar way as above
\begin{EQA}
	\P\bigl( \vthetav_{\GP} - \tilde{\thetav}_{\GP} \in A \cond \Yv \bigr)
	& \geq & 
	\frac{\bigl\{ 1 - \err(\rups) \bigr\} 
			\P\bigl( \DPGPt^{-1} \gammav \in A(\rups) \bigr)}
		{\bigl\{ 1 + \err(\rups) \bigr\} \P\bigl( \| \HPc \DPGPt^{-1} \gammav \| \leq \rups \bigr) +
			\CONST \ex^{ - (\dimt_{\GP} + \xx)/2}
		} 
	\\
	& \geq & 
	\frac{\bigl\{ 1 - \err(\rups) \bigr\} 
			\bigl\{ \P\bigl( \DPGPt^{-1} \gammav \in A \bigr) - \rho(\rups) \bigr\}}
		{\bigl\{ 1 + \err(\rups) \bigr\} \P\bigl( \| \HPc \DPGPt^{-1} \gammav \| \leq \rups \bigr) +
			\CONST \ex^{ - (\dimt_{\GP} + \xx)/2}
		} 	\, .
\label{rorLGttGuuGl2}
\end{EQA}
For the case of an arbitrary possibly non-symmetric \( A \), the proof is similar with the use of 
\eqref{d324d4efppm3} instead of \eqref{4d324d4efppmz2}.
%

\paragraph{Proof of Theorem~\ref{Ccontactionrate} and of Theorem~\ref{CThonestCS}}

The difference \( \vthetav_{\GP} - \thetavs \) can be decomposed as
\begin{EQA}
	\vthetav_{\GP} - \thetavs 
	& = &
	\bigl( \vthetav_{\GP} - \tilde{\thetav}_{\GP} \bigr) + \bigl( \tilde{\thetav}_{\GP} - \thetavs \bigr)
	= 
	\bigl( \vthetav_{\GP} - \tilde{\thetav}_{\GP} \bigr) + \bigl( \tilde{\thetav}_{\GP} - \thetavs_{\GP} \bigr) 
	+ \bigl( \thetavs_{\GP} - \thetavs \bigr) .
\label{vtvGtsGGtss}
\end{EQA}
Theorem~\ref{TestlosspMLE} with \( \xx = \log n \) allows to bound with high probability
\begin{EQA}
	\bigl\| \QP (\tilde{\thetav}_{\GP} - \thetavs) \bigr\|^{2}
	& \lesssim &
	\| \QP \DPGP^{-2} \GP^{2} \thetavs \|^{2}
	+ \zq^{2}(\BB_{\QPGP}^{2},\xx) 
	\leq 
	\| \QP \DPGP^{-2} \GP^{2} \thetavs \|^{2}
	+ \tr(\BB_{\QPGP}^{2}) + \log n ,
\label{QttGttsG2tr22}
\end{EQA}
where \( \BB_{\QPGP}^{2} = \QP \DPGP^{-2} \HPc^{2} \DPGP^{-2} \QP^{\T} \).
Moreover, as \( \HPc^{2} \lesssim \DPGP^{2} \), we bound 
\( \tr \bigl( \BB_{\QPGP}^{2} \bigr) \lesssim \tr \bigl( \QP \DPGP^{-2} \QP^{\T} \bigr) \).
Further, Theorem~\ref{TnonparBvMm} yields on a random set \( \Omega(\xx) \) for \( \xx = \log n \)
\begin{EQA}
	\P\bigl( \| \QP (\vthetav_{\GP} - \tilde{\thetav}_{\GP}) \| \geq \rr \cond \Yv \bigr)
	& \lesssim &
	\PG\bigl( \| \QP \DPGPt^{-1} \gammav \| \geq \rr \bigr) + 1/n  .
\label{PQvtGttGrQr2}
\end{EQA}
Now we apply Theorem~\ref{TexpbLGA} with \( \rr = \rr_{\QP} = \zq(\QP \DPGPt^{-2} \QP^{\T},\xx) 
\leq \sqrt{\tr (\QP \DPGPt^{-2} \QP^{\T})} + \sqrt{2\xx} \)
to the Gaussian quadratic form \( \| \QP \DPGPt^{-1} \gammav \|^{2} \).
The desired result \eqref{PDvtGtsCrGY} follows by \eqref{DPGPm1Cd3rG} of Theorem~\ref{TFiWititG} and by the bias bound \eqref{trQDGm2G22Q2DG2}.


To check \eqref{Palpo1AGra}, note first that by definition,
it holds for the true parameter \( \thetavs \):
\begin{EQA}
	\P\bigl( \thetavs \in \CA_{\QPGP}(\rr) \bigr)
	&=&
	\P\bigl( \bigl\| \QP (\tilde{\thetav}_{\GP} - \thetavs) \bigr\| \leq \rr \bigr) .
\label{PtsiSAQrQtGtsr}
\end{EQA}
The Fisher expansion \eqref{DGttGtsGDGm13rG} \( \tilde{\thetav}_{\GP} - \thetavs_{\GP} \approx \DPGP^{-2} \nabla \zeta \) of Theorem~\ref{TFiWititG} combined with the CLT 
\( \VP^{-1} \nabla \zeta \tow \gammav \) for a standard normal \( \gammav \) reduces the latter question to Gaussian probability
\begin{EQA}
	\P\bigl( \bigl\| \QP (\tilde{\thetav}_{\GP} - \thetavs) \bigr\| \leq \rr \bigr)
	& \approx &
	\P\Bigl( \bigl\| \QP \bigl(\DPGP^{-2} \VP \gammav + \thetavs_{\GP} - \thetavs\bigr) \bigr\| \leq \rr \Bigr) .
\label{PPQQDGm2ttsGtts}
\end{EQA}
By Gaussian comparison Theorem~\ref{Tgaussiancomparison3}, 
the impact of the bias \( \thetavs_{\GP} - \thetavs \)
is negligible under the undersmoothing condition 
\( \| \QP (\thetavs_{\GP} - \thetavs) \|^{2} \ll \tr \bigl( \QP \DPGP^{-2} \VP^{2} \DPGP^{-2} \QP^{\T} \bigr) \). 
Combining with Theorem~\ref{TestlosspMLE} yields 
in view of \( \DPGP^{-2} \VP^{2} \DPGP^{-2} \leq \DPGP^{-2} \)
\begin{EQA}
	1 - \alp
	&=&
	\PG\bigl( \bigl\| \QP \DPGPt^{-1} \gammav \bigr\| \leq \rr_{\alp} \bigr)
	\approx
	\P\bigl( \bigl\| \QP \DPGP^{-1} \gammav \bigr\| \leq \rr_{\alp} \bigr)
	\leq 
	\P\Bigl( \bigl\| \QP \DPGP^{-2} \VP \gammav \bigr\| \leq \rr_{\alp} \Bigr)
	\\
	& \approx &
	\P\bigl( \bigl\| \QP (\tilde{\thetav}_{\GP} - \thetavs) \bigr\| \leq \rr_{\alp} \bigr),
\label{1maPGPQm2111}
\end{EQA}
that is, the credible set \( \CA_{\QPGP}(\rr_{\alp}) \) is an asymptotically valid confidence set.

\paragraph{Proof of Theorems~\ref{Tconcdens} through \ref{TGLMBvM}}
It suffices to check the conditions of the general results from Section~\ref{SnonBvM}.
We start with the log-density model.
First we show \nameref{LL0ref}.
Remind that \( \E L(\thetav) = n \bigl\{ \langle \Psimean, \thetav \rangle - \cdens(\thetav) \bigr\} \) and \( \IF(\thetav) = - \nabla^{2} \E L(\thetav) = n \nabla^{2} \cdens(\thetav) \).
%
%
Fix any \( \thetav \in \Thetad \) and
denote also \( \Psid = E_{\thetav} \Psi(X_{1}) \) and define for \( \uv \in \R^{\mms} \) 
with \( \| \HPc \uv \| = \rr \)
and any \( t \)
\begin{EQA}
	\qq(t)
	& \eqdef &
	\int \exp\bigl\{ \bigl\langle \Psi(x) - \Psid ,\thetav + t \uv \bigr\rangle - \cdens(\thetav) \bigr\} \Pdom(dx) 
	=
	\int \exp\bigl\{ t \bigl\langle \Psi(x) - \Psid , \uv \bigr\rangle \bigr\} P_{\thetav}(dx),
	\\
	\qq_{k}(t)
	& \eqdef &
	\frac{d^{k}\qq(t)}{dt^{k}} 
	=
	\int \bigl\langle \Psi(x) - \Psid ,\uv \bigr\rangle^{k} \, 
	\exp\bigl\{ t \bigl\langle \Psi(x) - \Psid , \uv \bigr\rangle \bigr\} P_{\thetav}(dx) ,
	\quad
	k \geq 1.
\label{q0tuietPxpbuPd}
\end{EQA}
Due to \nameref{Psithetaref},
all these quantities are well defined.
Moreover, \( \qq(0) = 1 \) and does not depend on \( \thetav, \uv \) while \( \qq_{1}(0) = 0 \).
Also
\begin{EQA}
	\qq_{2}(0)
	&=&
	\langle \nabla^{2} \cdens(\thetav) \uv, \uv \rangle .
\label{q20ln2cdtvvvr}
\end{EQA}
Further, define
\begin{EQA}
	h(t)
	& \eqdef &
	\log \qq(t)  
	=
	\cdens(\thetav + t \uv) - \cdens(\thetav) - t \bigl\langle \Psid, \uv \bigr\rangle .
\label{htedelqtcdtuv}
\end{EQA}
Then
\begin{EQA}
	\dltw_{k}(\thetav,\uv)
	&=&
	- n 
    \frac{1}{k!}
      \frac{d^{k}}{dt^{k}} \cdens(\thetav + t \uv) \biggr|_{t=0} 
    =
	- n 
    \frac{1}{k!}
      \frac{d^{k}}{dt^{k}} h(t) \biggr|_{t=0} 
    \, ,
    \quad 
    k = 3, 4,
\label{dmtudenst0}
\end{EQA}
Straightforward calculus yields
\begin{EQ}[rcl]
    h^{(3)}(0)
    &=&
    - \qq_{3}(0) + 3 \qq_{2}(0) \, \qq_{1}(0) - 2 \qq_{1}^{3}(0) \, ,
    \\
    h^{(4)}(0)
    &=&
    - \qq_{4}(0) 
    + 4 \qq_{3}(0) \, \qq_{1}(0)
    + 3 \qq_{2}^{2}(0) 
    - 12 \qq_{2}(0) \, \qq_{1}^{2}(0) 
    + 6 \qq_{1}^{4}(0) \, .
\label{6h14th04t}
\end{EQ}
Now \nameref{Psiuref} implies
\begin{EQA}
	\bigl| \qq_{k}(0) \bigr|
	& \leq &
	\bigl\{ \CONSTfour^{2} \, \qq_{2}(0) \bigr\}^{k/2}
	\leq 
	\bigl\{ \CONSTfour^{2} \, \langle \nabla^{2} \cdens(\thetav) \uv, \uv \rangle \bigr\}^{k/2}
	, 
	\quad
	k = 3,4 .
\label{k1qktCpskqt}
\end{EQA}
As \( \HPc^{2} \geq \nsize \nabla^{2} \cdens(\thetav) \) and \( \| \HPc \uv \|^{2} = \rr^{2} \),
this yields for some absolute constant \( \CONST_{3}, \CONST_{4} \)
\begin{EQA}
    \dltw_{3}(\thetav,\uv)
    & \leq &
    \CONST_{3} \, \CONSTfour^{3} \, \nsize (\rr^{2}/\nsize)^{3/2},
    \quad
    \dltw_{4}(\thetav,\uv)
    \leq 
    \CONST_{4} \, \CONSTfour^{4} \, \nsize (\rr^{2}/\nsize)^{2} . 
\label{d3tur3n12r4n2}
\end{EQA}
%
Now we check \nameref{ED0ref} for 
\( \nabla \zeta = \Spsi - \E \Spsi \) and \( \VP^{2} = n \, \nabla^{2} \cdens(\thetavs) \). 
%
%
Let \( \| \HPc \uv \| = \lambda \).
It holds from \eqref{gtdelointTPxm0} due to the i.i.d. structure of the data 
in view of 
\( \nabla \cdens(\thetavs) = \Psimean \) 
\begin{EQA}
	&& \nquad
	\log \E \exp \bigl\{ \langle \nabla \zeta , \uv \rangle \bigr\}
	=
	n \log E \exp \bigl\{ \langle \Psi(X_{1}) - \Psimean \, , \uv \rangle \bigr\} 
	\\
	&=&
	n \bigl\{ \cdens(\thetavs + \uv) 
    	- \langle \nabla \cdens(\thetavs) \, , \uv \rangle 
    		\bigr\}
    =
    \frac{n}{2} 
    	\bigl\langle \nabla^{2} \cdens(\thetavs + t \uv) \, \uv, \uv \bigr\rangle 
\label{EtTSPXngt}
\end{EQA}
for \( t \in [0,1] \). 
By \nameref{phiref}, 
\( \| \uv \| \leq \| \HPc^{-1} \| \, \| \HPc \uv \| \leq \CONSTi_{\phi,1} \nsize^{-1/2} \lambda \).
By the Cauchy-Schwarz inequality and \nameref{Psithetaref}
\begin{EQA}
	\bigl| \bigl\langle \Psi(x), \uv \bigr\rangle \bigr|^{2}
	& \leq &
	q_{\mms}^{-2} \| \uv \|^{2} \, \sum_{j = 1}^{\mms} \psi_{j}^{2}(x) q_{j}^{2} 
	\leq 
	\CONSTi_{\phi,1}^{2} \, \CONSTpsi^{2} \, \frac{\lambda^{2} \, q_{\mms}^{-2}}{\nsize} 
	\leq 
	\CONSTi_{\phi,1}^{2} \, \CONSTpsi^{2} \, \frac{\lambda^{2} \log^{2}(\nsize)}{\nsize^{2/3}} \, .
\label{PxPdVPvvrCu}
\end{EQA}
If the latter value is smaller than a constant \( \CONST \) then by \eqref{q0tuietPxpbuPd}
\begin{EQA}
	\bigl\langle \nabla^{2} \cdens(\thetavs + t \uv) \, \uv, \uv \bigr\rangle
	& \leq &
	\ex^{\CONST} \bigl\langle \nabla^{2} \cdens(\thetavs) \, \uv, \uv \bigr\rangle .
\label{Cl2Cmnna2detnl}
\end{EQA}
This yields \nameref{ED0ref} with \( \gmb \lesssim \nsize^{1/3} / \log (\nsize) \lesssim \mms \).
%
%
Theorem~\ref{LLbrevelocroB} ensure 
for \( \zq_{\GP} = \sqrt{\dimA_{\GP}} \, + \sqrt{2 \log n} \)
the probability bound \eqref{uTzDGvTxm12z}
\( \bigl\| \DPGP^{-1} \nabla \zeta \bigr\| \leq \zq_{\GP} \)
on a random set \( \Omega_{n} \) with \( \P\bigl( \Omega_{n} \bigr) \geq 1 - 3/n \).
Now all the statements of Theorem~\ref{Tconcdens} and \ref{Tconcdensex} follow directly from the general 
results of Section~\ref{SnonBvM}.

Now we check the general conditions for the GLM starting with \nameref{LL0ref}.
Let \( \uv \in \R^{\mms} \) with \( \| \HPc \uv \| \leq \rr \).
It holds by \( |\langle \Psi_{i},\thetav \rangle| \leq \CONSTPsi \), \nameref{phi4ref},
and \nameref{PsiuGRref} for \( k=3,4 \)
\begin{EQA}
	\bigl| \dltw_{k}(\thetav,\uv) \bigr|
	&=&
	\biggl| 
		\sum \cdens^{(k)}\bigl( \langle \Psi_{i},\thetav \rangle \bigr) \langle \Psi_{i},\uv \rangle^{k} 
	\biggr|
	\leq 
	\| \cdens^{(k)} \|_{\infty} \biggl( \sum \langle \Psi_{i},\uv \rangle^{4} \biggl)^{k/4}
	\\
	& \leq &
	\nsize^{k/2 - 1}
	\| \cdens^{(k)} \|_{\infty} \biggl( \sum \langle \Psi_{i},\uv \rangle^{2} \biggr)^{k/2}
	\leq 
	\CONST \rr^{k} \nsize^{-k/2+1}
\label{Crknk2m1fkislPi}
\end{EQA}
yielding \eqref{dmtuCmnumdd}.
Further, by independence of the \( Y_{i} \)'s,
\begin{EQA}
	\log \E \exp\bigl\{ \lambda \langle \nabla\zeta, \uv \rangle \bigr\}
	&=&
	\sum \log \E \exp\bigl\{ \lambda \eps_{i} \langle \Psi_{i}, \uv \rangle \bigr\} .
\label{YilEnzeiPiu}
\end{EQA}
Under \nameref{PsithetaGRref} and \nameref{IFtref}, one can bound for any \( \uv \)
with \( \| \HPc \uv \| = 1 \)
\begin{EQA}
	\bigl| \langle \Psi_{i}, \uv \rangle \bigr|
	& \leq &
	\CONSTpsi \sqrt{\mms} \, \| \uv \|
	\leq 
	\CONSTpsi \, \CONSTIF \sqrt{\mms} \, \| \HPc \uv \| \nsize^{-1/2} 
	\leq 
	\CONSTpsi \, \CONSTIF \sqrt{\mms} \, \nsize^{-1/2} \, .
\label{CfCIn1m2m12uPi}
\end{EQA}
Thus, \( \bigl| \lambda \langle \Psi_{i}, \uv \rangle \bigr| \leq 
\lambda \CONSTIF \sqrt{\mms} \, \nsize^{-1/2} \), and, for 
\( |\lambda| \leq \varrho / \bigl( \CONST \sqrt{\mms} \, \nsize^{-1/2} \bigr) \),
it follows by \nameref{epsiref}
\begin{EQA}
	\sum \log \E \exp\bigl\{ \lambda \eps_{i} \langle \Psi_{i}, \uv \rangle \bigr\}
	& \leq &
	\sum \frac{\nunu^{2} \lambda^{2} \sigma_{i}^{2}}{2} \langle \Psi_{i}, \uv \rangle^{2}
	\leq 
	\frac{\nunu^{2} \lambda^{2}}{2} \| \HPc \uv \|^{2}
	=
	\frac{\nunu^{2} \lambda^{2}}{2}
\label{22n02l22Hu2si}
\end{EQA}
yielding \nameref{ED0ref} with \( \lambda \asymp \nsize^{1/3} \).
We complete the proof as in the log-density case.

\subsection{Proof of Theorem~\ref{Tpostmeanv}}
It holds
\begin{EQA}
	\vthetavb_{\GP} - \tilde{\thetav}_{\GP}
	&=&
	\frac{\int \bigl( \thetav - \tilde{\thetav}_{\GP} \bigr) \exp \LGP(\thetav) d\thetav}{\int \exp \LGP(\thetav) d\thetav} \, .
\label{vtbGttGfttGdt}
\end{EQA}
The use of \( \nabla \LGP(\tilde{\thetav}_{\GP}) = 0 \) helps to represent
with \( \CAt(\rups) = \{ \uv \colon \| \HPc \uv \| \leq \rups \} \)
\begin{EQA}
	\QP \bigl( \vthetavb_{\GP} - \tilde{\thetav}_{\GP} \bigr) 
	& = &
	\frac{\left( \int_{\| \HPc \uv \| \leq \rups} + \int_{\| \HPc \uv \| > \rups} \right) 
		\QP \uv \,
		\exp\bigl\{ \LGP( \tilde{\thetav}_{\GP} + \uv) - \LGP(\tilde{\thetav}_{\GP}) 
				- \bigl\langle \nabla \LGP(\tilde{\thetav}_{\GP}), \uv \bigr\rangle \bigr\} d\uv}
		{\int \exp\bigl\{ \LGP( \tilde{\thetav}_{\GP} + \uv) - \LGP(\tilde{\thetav}_{\GP}) 
				- \bigl\langle \nabla \LGP(\tilde{\thetav}_{\GP}), \uv \bigr\rangle \bigr\}d\uv} \, .
\label{QtbGttGfCrCRcLGRG}
\end{EQA}
Now, with \( f_{\thetav}(\uv) = \E \LGP(\thetav + \uv) \), define \( f(\uv) \) by using 
\( \thetav = \tilde{\thetav}_{\GP} \), that is, 
\( f(\uv) = f_{\tilde{\thetav}_{\GP}}(\uv) \). 
Linearity of the stochastic part of \( \LGP(\thetav) \) implies 
\begin{EQA}
	\LGP( \tilde{\thetav}_{\GP} + \uv) - \LGP(\tilde{\thetav}_{\GP}) 
				- \bigl\langle \nabla \LGP(\tilde{\thetav}_{\GP}), \uv \bigr\rangle
	&=&
	f(\uv) - f(0) - f'(0,\uv), 
\label{fuf0fpuLGt}
\end{EQA} 
and it holds 
\begin{EQA}
	\| \QP \bigl( \vthetavb_{\GP} - \tilde{\thetav}_{\GP} \bigr) \| 
	& \leq &
	\rho_{0}(\rups) + \rho_{1}(\rups)
\label{QtbGttGfCrCRcLGRG}
\end{EQA}
with
\begin{EQA}
	\rho_{0}(\rups)
	& \eqdef &
	\biggl\| \frac{\int_{\| \HPc \uv \| \leq \rups}  
		\QP \uv
		\exp\bigl\{ f(\uv) - f(0) - f'(0,\uv) \bigr\} d\uv }
		{\int_{\| \HPc \uv \| \leq \rups} 
	\exp\bigl\{ f(\uv) - f(0) - f'(0,\uv) \bigr\}d\uv} \biggr\|\, ,
	\\
	\rho_{1}(\rups)
	& \eqdef &
	\biggl\| \frac{\int_{\| \HPc \uv \| > \rups}  
		\QP \uv
		\exp\bigl\{ f(\uv) - f(0) - f'(0,\uv) \bigr\} d\uv}
		{\int_{\| \HPc \uv \| \leq \rups} 
	\exp\bigl\{ f(\uv) - f(0) - f'(0,\uv) \bigr\}d\uv} \biggr\|\, ,
\label{r0r1r0fuv0uG}
\end{EQA}
As \( - \nabla^{2} f(0) = \DPGPt^{2} \), Lemma~\ref{Lintfxupp2r}
and Theorem~\ref{PrhoQPBvM} imply
\begin{EQA}
	\rho_{0}(\rups)
	& \lesssim &
	\frac{\delta_{3} \E \| \QP \DPGPt^{-1} \gammav \|}{\P\bigl( \| \HPc \DPGPt^{-1} \gammav \| \leq \rups \bigr)}
	\lesssim
	\delta_{3} \sqrt{\dimt_{\QPGP}} 
\label{r0d3stQ2Dm2}
\end{EQA}
with \( \dimt_{\QPGP} = \tr (\QP \DPGPt^{-2} \QP^{\T}) \).
For bounding the term \( \rho_{1}(\rups) \), we apply the bound from Lemma~\ref{CELLDtuvLt} and then 
Lemma~\ref{TGaussintext} and \ref{TGaussintextquad} with \( \TAU = \HPc \DPGPt^{-1} \)
and \( \dimAA = \tr (\TAU^{\T} \TAU) = \dimt_{\GP} \).
The use of \( 2 \| \QP \uv \| \leq 1 + \| \QP \uv \|^{2} \) yields on \( \Omega(\xx) \)
\begin{EQA}
	\rho_{1}(\rups)
	& \lesssim &
	\dimt_{\QPGP} \, \exp\{ - (\dimt_{\GP} + \xx)/2 \} .
\label{r1r0ltQ2Dm2epA}
\end{EQA}
The second moment of the expression 
\( \bigl\langle \uv, \DPGPt \bigl( \vthetav_{\GP} - \tilde{\thetav}_{\GP} \bigr) \bigr\rangle \)
given \( \Yv \) and a unit vector \( \zv \) is evaluated similarly. 
One gets
\begin{EQA}
	&& \nquad
	\E \Bigl[ 
          \bigl\langle \zv, \DPGPt \bigl( \vthetav_{\GP} - \tilde{\thetav}_{\GP} \bigr) \bigr\rangle^{2} 
        \cond \Yv 
    \Bigr] - 1 
    =
    \frac{\int \bigl\langle \zv, \DPGPt \bigl( \vthetav_{\GP} - \tilde{\thetav}_{\GP} \bigr)\bigr\rangle^{2} 
     		\exp \LGP(\thetav) \, d\thetav}
		  {\int \exp \LGP(\thetav) \, d\thetav} - 1
	\\
	&=&
	\frac{\left( \int_{\| \HPc \uv \| \leq \rups} + \int_{\| \HPc \uv \| > \rups} \right)
			\bigl[ \bigl\langle \zv, \DPGPt \uv \bigr\rangle^{2} - 1 \bigr]
			\exp\bigl\{ f(\uv) - f(0) - f'(0,\uv) \bigr\} d\uv }
		 {\int \exp\bigl\{ f(\uv) - f(0) - f'(0,\uv) \bigr\}d\uv}
	\\
	&=&
	\rho_{2} + \rho_{3} \, .
\label{ro2ro3EcY1fuf0}
\end{EQA}
Similarly to \eqref{r1r0ltQ2Dm2epA}, one can get 
\( |\rho_{3}| \leq \dimt_{\QPGP} \, \exp\{ - (\dimt_{\GP} + \xx)/2 \} \).
For the value \( |\rho_{2}| \), we use symmetricity of 
\( \UVd = \CAt(\rups) = \bigl\{ \uv \colon \| \HPc \uv \| \leq \rups \bigr\} \) and Lemma~\ref{Lintfxupp2}
yielding \( |\rho_{2}| \lesssim \err \).

\subsection{Proof of Theorem~\ref{TbarthetavG}}
We start with an extension of Theorem~\ref{PrhoQPBvM} to some non-symmetric sets \( A \).
Let \( \av \) be a possibly random vector in \( \R^{\dimp} \). 
The next result assumes that \( \av \) is sufficiently small so that 
\( \| \DPGPt \av \| \leq 1 \). 

\begin{proposition}
\label{TGARav}
    Let the conditions of Theorem~\ref{PrhoQPBvM} hold.
Suppose that \( \rups \) satisfies the conditions \eqref{LmgfquadELGP} and \eqref{CONSTruAxx1} with 
\( \xx = 2 \log n \).
Let a random vector \( \av \) satisfy \( \| \DPGPt \av \| \leq 1 \) 
on the set \( \Omega(\xx) \) from \eqref{uTzDGvTxm12z}.
Then for any set \( A \in \cc{B}_{s}(\R^{\dimp}) \), it holds on \( \Omega(\xx) \)
\begin{EQA}
	&& \nquad
		\left| 
		\P\bigl( \vthetav_{\GP} - \tilde{\thetav}_{\GP} - \av \in A \cond \Yv \bigr)
		- \PG\bigl( \DPGPt^{-1} \gammav - \av \in A \bigr)  
	\right|
	\\
	& \lesssim &
	\Bigl\{ \err(\rups) + \delta_{3}(\rups) \| \DPGPt \av \| \Bigr\}
	\PG\bigl( \DPGPt^{-1} \gammav \in A \bigr) + 1/n \, .
\label{erd3DtaDtm1ACn}
\end{EQA}
\end{proposition}

This result can be proved in the same line as Theorem~\ref{TnonparBvMm} using Lemma~\ref{Lintavshift}.
Taking into account the Gaussian approximation result from Corollary~\ref{CnonparBvMm},
we only have to compare the posterior probability of 
\( \| \QP (\vthetav_{\GP} - \bar{\vthetav}_{\GP}) \| \leq \rr \) with 
\( \PG\bigl( \| \QP \DPGPt^{-1} \gammav \| \leq \rr \bigr) \).
Let \( \av \) be defined as
\begin{EQA}
	\av 
	&=&
	\tilde{\thetav}_{\GP} - \bar{\vthetav}_{\GP} .
\label{avQttGbtvtG}
\end{EQA}  
As \( \QP = \QP \Pi \) for a projector \( \Pi \), it also holds with \( \av_{0} = \Pi \av \)
\begin{EQA}
	\| \QP (\vthetav_{\GP} - \bar{\vthetav}_{\GP}) \|
	&=&
	\| \QP (\vthetav_{\GP} - \tilde{\thetav}_{\GP} - \av_{0}) \|
\label{QvtGbtGQvtav}
\end{EQA}
and
\begin{EQA}
	\P\bigl( \vthetav_{\GP} - \bar{\vthetav}_{\GP} \in \CS_{\QP}(\rr) \cond \Yv \bigr)
	&=&
	\P\bigl( \vthetav_{\GP} - \tilde{\thetav}_{\GP} - \av_{0} \in \CS_{\QP}(\rr) \cond \Yv \bigr)
\label{PvGbGEGYttGY}
\end{EQA}
Now Theorem~\ref{TGARav} implies
\begin{EQA}
	&& \nquad
	\left| \P\bigl( \vthetav_{\GP} - \tilde{\thetav}_{\GP} - \av_{0} \in \CS_{\QP}(\rr) \cond \Yv \bigr)
	- \PG\bigl( \DPGPt^{-1} \gammav - \av_{0}  \in \CS_{\QP}(\rr) \bigr) \right|
	\\
	& \lesssim &
	\CONST \Bigl\{ \err + \delta_{3} \| \DPGPt \av_{0} \| + n^{-1} \Bigr\} \, .
\label{PvGbGEGYttGY}
\end{EQA}
Theorem~\ref{Tpostmeanv} yields that the norm of \( \DPGPt \av_{0} \) can be bounded on \( \Omega(\xx) \) as
\begin{EQA}
	\| \DPGPt \av_{0} \|
	=
	\| \DPGPt \Pi \av \|
	& \lesssim &
	\delta_{3} \sqrt{\dimt_{\Pi}} 
	+ n^{-1} \dimt_{\Pi} \, 
\label{aled3sptQG}
\end{EQA}
with \( \dimt_{\Pi} = \tr (\Pi \DPGPt^{-2} \Pi \DPGPt^{2} \Pi) \).
It remains to compare two Gaussian probabilities of
\( \| \QP \DPGPt^{-1} \gammav \| \leq \rr \) and of 
\( \| \QP (\DPGPt^{-1} \gammav - \av_{0}) \| \leq \rr \).
For this one can apply the Pinsker inequality. 
However, the Gaussian comparison result of Theorem~\ref{Tgaussiancomparison3} provides a more precise bound 
in view of the elliptic shape of the considered credible sets:
\begin{EQA}
	\bigl| \P\bigl( \| \QP (\DPGPt^{-1} \gammav - \av_{0}) \| \leq \rr \bigr) 
		- \P\bigl( \| \QP \DPGPt^{-1} \gammav \| \leq \rr \bigr)
	\bigr|
	& \leq &
	\frac{\| \QP \av_{0} \|^{2}}{\| \QP \DPGPt^{-2} \QP^{\T} \|_{\Fr}} \, .
\label{PQDGa0Qa02Fr}
\end{EQA}
Now the assertion follows by one more application of Theorem~\ref{Tpostmeanv} 
in view of \( \dimt_{\QPGP} \leq  \| \QP \DPGPt^{-2} \QP^{\T} \|_{\Fr}^{2} \, \dimt_{\Pi} \).

\subsection{Proof of Theorem~\ref{TnonparamBvM}}
We assume that all the conditions are fulfilled for the smaller prior covariance \( \GP^{2} \),
and all error terms correspond to that prior. 
Corollary~\ref{CnonparBvMm} implies on a set of probability at least \( 1 - 1/n \)
for any measurable set \( A \)
\begin{EQA}
	\left| \P\bigl( \vthetav_{\GP} - \tilde{\thetav}_{\GP} \in A \cond \Yv \bigr)
    - \PG\bigl( \DPGPt^{-1} \gammav \in A \bigr) 
    \right|
    & \lesssim &
    \delta_{3} + n^{-1} \, .
\label{ernm1QtDrlQ}
\end{EQA}
Similarly, again on a set of probability at least \( 1 - 1/n \)
\begin{EQA}
	\left| \P\bigl( \vthetav_{\GPm} - \tilde{\thetav}_{\GPm} \in A \cond \Yv \bigr)
    - \PG\bigl( \DPt_{\GPm}^{-1} \gammav \in A \bigr) 
    \right|
    & \lesssim &
    \delta_{3} + n^{-1} \, .
\label{ernm1QtDrlQ11}
\end{EQA}
Define 
\( \av \eqdef \tilde{\thetav}_{\GP} - \tilde{\thetav}_{\GP_{1}} \).
The bound of Theorem~\ref{Tgaussiancomparison3} yields
\begin{EQA}
  	&& \nquad
  	\left| 
    	\PG\bigl( \| \QP (\DPGPt^{-1} \gammav - \av) \| \leq \rr \bigr)
    	-
    	\PG\bigl( \| \QP \DPt_{\GPm}^{-1} \gammav \| \leq \rr \bigr) 
  	\right|
  	\\
  	& \lesssim &
  	\frac{1}{\| \QP \DPt_{\GP}^{-2} \QP^{\T} \|_{\Fr}} 
  	\bigl\{ \tr\bigl(\QP (\DPGPt^{-2} - \DPt_{\GP_{1}}^{-2}) \QP^{\T}\bigr) 
  	+ \bigl\| \QP \av \bigr\|^{2} \bigr\}.
\end{EQA}
Putting all bounds together completes the proof.

\subsection{Proof of Lemma~\ref{Ltruncprior} through \ref{Lsmoothnessts}}
The condition \( \GP_{\mm}^{-2} \bigl( \Id - \Proj_{\VV_{m}} \bigr) = 0 \) effectively means
that the prior is limited to the subspace \( \VV_{\mm} \) and moreover, its norm is bounded
by \( \gp_{\mm} \) on this subspace.
This implies \eqref{fC1nC2ngm2mpGm} in view of \eqref{uVmGu2C1C2}.

Next we show the statement \eqref{C4mC3mpGAtm} of Lemma~\ref{Lsmoothprior}.
For \( J \geq 1 \), by \eqref{uVmGu2C1C2} and \eqref{uniVmgm2u2uiVm}
\begin{EQA}
	\tr\bigl[ \HPc^{2} \bigl\{ \IF(\thetav) + \GP^{2} \bigr\}^{-1} \Proj_{\VV_{J}} \bigr]
	& \leq &
	\tr \Proj_{\VV_{J}}
	=
	J,
	\\
	\tr\bigl[ \HPc^{2} \bigl\{ \IF(\thetav) + \GP^{2} \bigr\}^{-1} \Proj_{\VV_{J}} \bigr]
	& \geq &
	\tr \bigl[ \IF(\thetav) \bigl\{ \IF(\thetav) + \gp_{J}^{2} \Id \bigr\}^{-1} \Proj_{\VV_{J}} \bigr]
	\\
	& \geq &
	\frac{\CONST_{1,\IF} \, n}{\CONST_{1,\IF} \, n + \gp_{J}^{2}} \tr \Proj_{\VV_{J}}
	=
	\frac{J \, \CONST_{1,\IF} \, n}{\CONST_{1,\IF} \, n + \gp_{J}^{2}} \, .
\label{JC1FnC1FngJ2JV2}
\end{EQA}
Similarly
\begin{EQA}
	\tr\bigl[ \HPc^{2} \bigl\{ \IF(\thetav) + \GP^{2} \bigr\}^{-1} (\Id - \Proj_{\VV_{J}}) 
	\bigr]
	& \leq &
	\sum_{j \geq J} \frac{\CONST_{2,\IF} \, n}{\CONST_{2,\IF} \, n + \gp_{j}^{2}} \, .
\label{jJC2FnCngj2}
\end{EQA}
Further, by \eqref{sumjJgjm2C}
\begin{EQA}
	\sum_{j \geq J} \frac{\CONST_{2,\IF} \, n}{\CONST_{2,\IF} \, n + \gp_{j}^{2}}
	& \leq &
	\CONST_{2,\IF} \, n \sum_{j \geq J} \frac{1}{\gp_{j}^{2}} 
	\leq 
	\CONST_{2,\IF} \, n \, \, \CONST_{2,\gp} \, J \gp_{J}^{-2} \, .
\label{2InC2gJmJjJ}
\end{EQA}
Therefore,
\begin{EQA}
	J \, \frac{\CONST_{1,\IF} \, n}{\CONST_{1,\IF} \, n + \gp_{J}^{2}}
	\leq 
	\tr\bigl[ \HPc^{2} \bigl\{ \IF(\thetav) + \GP^{2} \bigr\}^{-1} \bigr]
	& \leq &
	J \bigl( 1 + \CONST_{2,\IF} \, \CONST_{2,\gp} \, n \gp_{J}^{-2} \bigr) .
\label{J1C2IC2gngJn2}
\end{EQA}
If \( J \) be such that \( \gp_{J}^{2} \asymp n \), then \( \dimA_{\GP}(\thetav) \asymp J \). 

Lemma~\ref{Lsmoothnessts} is a special case of Lemma~\ref{Lsmoothprior}
with \( \gp_{j}^{2} = \CGP j^{2s} \).

\bibliography{exp_ts,listpubm-with-url,semi_bvm}

\end{document}